
\input amstex \documentstyle{amsppt}
\magnification 1200
\hsize=16.5truecm \vsize=23 truecm \voffset=.4truecm
\refstyle{A}
\TagsOnRight
\NoBlackBoxes
\topmatter
\define\ppt#1{#1}
\ppt{}

\define\secNo{00}
\newcount\eqNo
\newcount\pcNo
\define\secno#1. {#1. \global\edef\secNo{#1}\global\eqNo=0\global\pcNo=0}
\define\ssecno#1. {#1. \global\edef\secNo{#1}\global\eqNo=0\global\pcNo=0}
\outer\def\iproclaim#1/#2/ #3. {\bgroup\csname proclaim\endcsname
                          {#3\ \dpcl#1/#2/}}
\define\eproclaim{\endproclaim\egroup}
\define\proof#1:{\demo{Proof#1}}
\define\QED {\qed\enddemo}
\define\examp#1/#2/#3:{\par\noindent{\bf Example \dpcl#1/#2/:\
                                        #3.}\hfill\break}
\define\subex#1 #2\Take#3\par{\noindent{{\bf #1. }{\sl#2}\hfill\break
                           Take\ #3\par}}

\define\Br{\overline}
\define\abs #1{{\left\vert#1\right\vert}}
\define\bra #1>{\langle #1\rangle}
\define\bracks #1{\lbrack #1\rbrack}
\define\ket #1{\vert#1\rangle}
\define\ketbra #1#2{{\vert#1\rangle\langle#2\vert}}
\define\lin{\mathop{{\fam0 lin\,}}\nolimits}
\define\norm #1{\left\Vert #1\right\Vert}
\define\Norm #1#2{#1\Vert #2#1\Vert}
\define\rank{\mathop{{\fam0 rank}}\nolimits}
\define\set #1{\left\lbrace#1\right\rbrace}
\define\stt{\,\vrule\ }
\define\th {\hbox{${}^{\text{th}}$}\ }  
\define\tr {\mathop{{\fam0 tr}}\nolimits}
\define\up#1{^{(#1)}}
\define\clinsp{\overline{\lin}}
\define\iprod{\bra\cdot,\cdot>}
\define\me<#1,#2,#3>{\langle#1\vert#2\vert#3\rangle} 
\define\idty{{\leavevmode{\roman1\mkern -5.4mu\roman I}}}
\define\Nl{{\Bbb N}} \define\Cx{{\Bbb C}}
\define\Rl{{\Bbb R}} \define\Ir{{\Bbb Z}}
\define\pfi{\varphi}
\define\eps{\varepsilon}
\define\em{\it} 
\chardef\hash=`\#
\define\1{^\dagger}
\define\wick{{\Cal W}} 
\define\wicc{{\bold W}} 
\define\A{{\Cal A}}
\define\BB{{\Cal B}}   
\define\HH{{\Cal H}}

\define\LL{{\Cal L}}   
\define\Mat{{\Cal M}}  
\define\No{{\Cal N}}   
\define\NN{{\Cal N}}   
\define\R{{\Cal R}}    
\define\SnU#1{{\roman S}_\nu{\roman U}(#1)}
\define\Tt{\widetilde T}
\define\ZZ{{\Cal Z}}
\define\cEx{{\Bbb E}}   
\define\cV{{\Cal V}}    
\define\crep#1{\lambda_{#1}} 
\define\creps#1{\widehat\lambda_{#1}} 
\define\cst#1{\omega_{#1}}   

\define\ideal{{\Cal J}}
\define\nTt{\Vert\Tt\Vert}   
\define\tA{{\Cal T}}         
\define\twist{\Theta}
\define\gauge{\alpha}        
\define\degree{\deg}
\define\kms#1{\tau_{#1}}   
\define\kmsl{\tau_\lambda} 
\define\Dff{\Omega}        
\define\DffT{\Dff_T}
\define\Dffwick{\Dff\wick}
\define\df{{\roman d}}         
\define\dx#1{{\df x_{#1}}}     
\define\dxs#1{{\df x_{#1}\1}}
\define\mo[#1,#2,#3]{a_{#3}a_{#2}a_{#1}}
\define\Qprod#1#2{\prod_{\textstyle #1 \atop \textstyle#2}}%
\define\aaa#1#2{A_{#1#2}A\1_{#1#2}}
\define\Bozejko{Bo{\accent 95 z}ejko}
\define\Wa{Wick algebra}
\define\HWa{Hermitian Wick algebra}
\define\cbo{collectively bounded}
\define\rep{representation}
\define\irrep{irreducible \rep}
\define\mCAR{$\mu$CAR}
\catcode`@=11
\def\ifundefined#1{\expandafter\ifx\csname
                        \expandafter\eat\string#1\endcsname\relax}
\def\atdef#1{\expandafter\def\csname #1\endcsname}
\def\atedef#1{\expandafter\edef\csname #1\endcsname}
\def\atname#1{\csname #1\endcsname}
\def\ifempty#1{\ifx\@mp#1\@mp}
\def\ifatundef#1#2#3{\expandafter\ifx\csname#1\endcsname\relax
                                  #2\else#3\fi}
\def\eat#1{}
\def\@cite #1,#2\@ver
   {\eachcite{#1}\ifx#2\@ut\else,\@cite#2\@ver\fi}
\def\cit#1{\cite{\@cite#1,\@ut\@ver}}
\def\eachcite#1{\atname{R@#1}}
\def\defonereftag#1=#2,{\atdef{R@#1}{#2}}
\def\defreftags#1, {\ifx\relax#1\relax \let\next\relax \else
           \expandafter\defonereftag#1,\let\next\defreftags\fi\next }
\def\lasteq{\secNo.\number\eqNo}
\def\deq#1(#2){{\ifempty{#1}\global\advance\eqNo by1
       \edef\n@@{\lasteq}\else\edef\n@@{#1}\fi
       \ifempty{#2}\else\global\atedef{E@#2}{\n@@}\fi\n@@}}
\def\eq#1(#2){\edef\n@@{#1}\ifempty{#2}\else
       \ifatundef{E@#2}{\global\atedef{E@#2}{#1}}%
                       {\edef\n@@{\atname{E@#2}}}\fi
       {\rm(\n@@)}}
\def\eqback#1{{(\advance\eqNo by -#1 \lasteq)}}

\def\eqgroup(#1){{\global\advance\eqNo by1
       \edef\n@@{\lasteq}\global\atedef{E@#1}{\n@@}}}
\def\lastpc{\secNo.\number\pcNo}
\def\dpcl#1/#2/{\ifempty{#1}\global\advance\pcNo by1
       \edef\n@@{\lastpc}\else\edef\n@@{#1}\fi
       \ifempty{#2}\else\global\atedef{P@#2}{\n@@}\fi\n@@}
\def\pcl#1/#2/{\edef\n@@{#1}%
       \ifempty{#2}\else
       \ifatundef{P@#2}{\global\atedef{P@#2}{#1}}%
                       {\edef\n@@{\atname{P@#2}}}\fi
       \n@@}
\def\Def#1/#2/{Definition~\pcl#1/#2/}
\def\Thm#1/#2/{Theorem~\pcl#1/#2/}
\def\Lem#1/#2/{Lemma~\pcl#1/#2/}
\def\Prp#1/#2/{Proposition~\pcl#1/#2/}
\def\Cor#1/#2/{Corollary~\pcl#1/#2/}
\def\Exa#1/#2/{Example~\pcl#1/#2/}
\catcode`@=12
\defreftags
AKLT=AKLT, Baez=Bae, Bargmann=Bar, Bergmann=Ber, BoSpei=BS1,
BoSpein=BS2, Blackadar=Bla, JOb=BElGJ, JOa=BEvGJ, BraRo=BR,
Connes=Con, Cuntz=Cu1, Cuntzb=Cu2, Daska=Das, Donohu=Don, Dykema=DN,
Faddeev=FRT, FCS=FNW, Fivel=Fiv, Folland=Fol, GoodWall=GW,
Green=Gre, Jones=Jo1, Jonesb=Jo2, Pallebook=Jor, QCS=JW, QCR=JSW1,
GPOTS=JSW2, KacRaina=KR, Klimek=KL, Worfuncal=KW, Manin=Man,
Mathematica=Mat, Meister=Mei, Foias=NF, Nica=NN, Phillips=Phi,
Pusz=Pus, PWor=PW, Radjavi=RR, RSimon=RS, Rehren=Reh, Rosenberg=Ros,
Rudin=Rud, FNC=SW, Shale=Sha, Schmuedgen=Sm\"u, Slowikow=Slo,
Speix=Sp2, Temperley=TL, vDaele=vDa, QTD=We1, FCL=We2, WessZ=WZ,
WoRims=Wo1, Woron=Wo2, Worona=Wo3, Woronb=Wo4, Zagier=Zag, ,
\def\sIn{0}       
\def\sH{1}        
\def\sHi{\sH.1}   
\def\sHc{\sH.3}   
\def\sHw{\sH.4}   
\def\sE {\sH.2}   
\def\sB{2}        
\def\sBt{\sB.1}   
\def\sBe{\sB.2}   
\def\sBs{\sB.3}   
\def\sBu{\sB.4}   
\def\sBp{\sB.5}   
\def\sBq{\sB.6}   
\def\sI{3}        
\def\sIq{\sI.1}   
\def\sIt{\sI.2}   
\def\sD{4}        
\def\sDf{\sD.1}   
\def\sDd{\sD.2}   
\def\sDb{\sD.3}   
\def\sG{5}        
\title Positive representations of    \\
       general commutation relations  \\
       allowing Wick ordering
\endtitle

\rightheadtext{Wick algebras }

\author P.E.T. J\o rgensen,
        L.M. Schmitt,
        and R.F. Werner
\endauthor

\leftheadtext{J\o rgensen--Schmitt--Werner}

\address
P.E.T. J\o rgensen,
Department of Mathematics, The University of Iowa, Iowa City,
Iowa 52242, U.S.A.
\endaddress
\email jorgen\@math.uiowa.edu\endemail
\address
L.M. Schmitt,
University of Aizu, Aizu-Wakamatsu, Fukushima Prefecture 965, Japan
\endaddress
\email lothar\@rsc.u-aizu.ac.jp \endemail
\address
R.F. Werner,
F.B. Physik, Universit\"at Osnabr\"uck, 49069 Osnabr\"uck,
Germany \endaddress
\email reinwer\@dosuni1.rz.Uni-Osnabrueck.de\endemail



\abstract
We consider the problem of representing in Hilbert space commutation
relations of the form
$$ a_ia_j^*=\delta_{ij}{\bold1} +
           \sum_{k\ell}T_{ij}^{k\ell} a_\ell^*a_k \quad,$$
where the $T_{ij}^{k\ell}$ are essentially arbitrary scalar
coefficients. Examples comprise the $q$-canonical commutation
relations introduced by Greenberg, Bozejko, and Speicher, and the
twisted canonical (anti-)commutation relations studied by Pusz and
Woronowicz, as well as the quantum group S$_\nu$U$(2)$.
Using these relations, any polynomial in the generators $a_i$ and
their adjoints can uniquely be written in ``Wick ordered form'' in
which all starred generators are to the left of all unstarred ones.
In this general framework we define the Fock representation, as well
as coherent representations. We develop criteria for the natural
scalar product in the associated representation spaces to be
positive definite, and for the relations to have representations by
bounded operators in a Hilbert space.
We characterize the relations between the generators $a_i$ (not
involving $a_i^*$) which are compatible with the basic relations.
The relations may also be interpreted as defining a non-commutative
differential calculus. For generic coefficients $T_{ij}^{k\ell}$,
however, all differential forms of degree 2 and higher vanish. We
exhibit conditions for this not to be the case, and relate them to
the ideal structure of the Wick algebra, and conditions of
positivity. We show that the differential calculus is compatible
with the involution iff the coefficients $T$ define a representation
of the braid group. This condition is also shown to imply improved
bounds for the positivity of the Fock representation.
Finally, we study the KMS states of the group of gauge
transformations defined by $a_j\mapsto\exp(it)a_j$.
\endabstract
\endtopmatter
\toc
\head0. Introduction
\page{3}\endhead
\subhead0.1. Main results
\page{3}\endsubhead
\subhead0.2. General background and related work
\page{5}\endsubhead
\head1. Basic Definitions and Examples
\page{8}\endhead
\subhead1.1. Hermitian Wick algebras
\page{8}\endsubhead
\subhead1.2. Examples
\page{11}\endsubhead
\subhead1.3. Coherent representations
\page{18}\endsubhead
\subhead1.4. A characterization of the Fock representation
\page{21}\endsubhead
\head2. Boundedness and Positivity
\page{22}\endhead
\subhead2.1. The coefficients $T_{ij}^{k\ell }$ as an operator
\page{22}\endsubhead
\subhead2.2. The operator $T$ in the examples
\page{24}\endsubhead
\subhead2.3. Bounds for small $T$
\page{27}\endsubhead
\subhead2.4. The universal bounded representation
\page{29}\endsubhead
\subhead2.5. Positivity of the Fock representation\ for positive $T$
\page{31}\endsubhead
\subhead2.6. $T$ satisfying the braid relations
\page{32}\endsubhead
\head3. Ideals in Wick algebras
\page{36}\endhead
\subhead3.1. Wick ideals and quadratic Wick ideals
\page{37}\endsubhead
\subhead3.2. Twisted canonical (anti-)commutation relations
\page{41}\endsubhead
\head4. Wick algebra relations as differential calculus
\page{48}\endhead
\subhead4.1. The algebra of differential operators
\page{48}\endsubhead
\subhead4.2. The algebra of differential forms
\page{50}\endsubhead
\subhead4.3. Differential calculus with braid relations
\page{52}\endsubhead
\head5. Gauge automorphisms and their KMS states
\page{56}\endhead
\head{} Acknowledgement
\page{59}\endhead
\head{} References
\page{59}\endhead
\endtoc
\ppt{\vfill}
\catcode`\@=11 \unvbox\tocbox@ \catcode`\@=12
\ppt{\vfill\eject}

\document

\head    \secno\sIn. Introduction   \endhead
\subhead \ssecno\sIn.1. Main results \endsubhead

This paper is about \rep s in Hilbert space of a wide class of
involutive algebras, which are obtained from -- often only finitely
many -- generators $a_i$ with the relations
$$ a_ia_j^*=\delta_{ij}{\idty} +
           \sum_{k\ell}T_{ij}^{k\ell} a_\ell^*a_k
\quad,\tag\deq(In.1)$$
where the $T_{ij}^{k\ell}$ are complex coefficients constrained only
by a hermiticity condition so that the above relation respects the
involution. The simplest case of such an algebra are certain
deformations of the unit disk \cit{Nica}. Extensions to relations of
the form $aa^*=f(a^*a)$ for a single generator, with more general
function $f$ have also been studied \cit{Daska}, and have
been shown to be analyzed completely in terms of a weighted shift
\cit{Klimek}. In this paper we look at generalizations to more than
one generator, in which no simple weighted shift \rep\ is possible.
Our paper does subsume many examples found in the literature,
however, and unifies, as well as substantially extends previously
known results.

In a basis free language the ``structure constants''
$T_{ij}^{k\ell}$ determine an operator
$\Tt:\HH\1\otimes\HH\to\HH\otimes\HH\1$, where $\HH$ is a Hilbert
space with a basis $\set{e_i}$ labelled as the generators, and
$\HH\1$ is the conjugate of $\HH$. The relations
$$ f\1\otimes g= \bra f,g>\idty + \Tt(f\1\otimes g)
\quad,\quad f\1\in\HH\1,\ g\in\HH \tag\deq(In.1')$$
then define an ideal in the algebra of all tensors over $\HH$ and
$\HH\1$ with tensor multiplication as product. The quotient of the
tensor algebra by this ideal is then an abstract involutive algebra,
which we denote by $\wick(T)$. A fundamental algebraic property of
$\wick(T)$ is that by using \eq(In.1) every
expression can be rearranged uniquely in ``Wick ordered form'', such
that, in every monomial, all starred generators are to the left of all
unstarred ones. To emphasize this property, $\wick(T)$ will be called
a ``\Wa''. We are mainly interested in ``positive \rep s'' of
$\wick(T)$, i.e.\ \rep s of the $a_i$ as operators on a Hilbert space
such that $a_i^*$ is a restriction of the operator adjoint of $a_i$.
Moreover, we would like to characterize the solutions of \eq(In.1)
by bounded operators.

For any choice of $T_{ij}^{k\ell}$, $\wick(T)$ has a unique \rep\
constructed from a cyclic vector $\Omega$ with the property that
$a_i\Omega=0$ for all generators $a_i$. This \rep, called the Fock
\rep, carries a natural hermitian scalar product, which, however,
often fails to be positive semi-definite. One of the main problems
addressed in this paper is to find criteria on $T_{ij}^{k\ell}$
implying the positivity of the Fock \rep. We generalize the Fock
\rep\ to the so-called coherent \rep s, the \rep s generated from a
cyclic vector $\Omega$ such that $a_i\Omega=\pfi_i\Omega$, with
$\pfi_i\in\Cx$. We will demonstrate the usefulness of this concept
for constructing \rep s.

We have shown in earlier work \cit{QCR} that, provided
the coefficients $T_{ij}^{k\ell}$ are sufficiently small, the Fock
\rep\ is positive and bounded, and is, moreover, the universal bounded
solution of the relations (i.e.\ every bounded solution is in a
quotient of the C*-algebra generated by the $a_i$ in the Fock
\rep). Moreover, the universal C*-algebra is independent of (small)
$T$, and is isomorphic, as a C*-algebra, to an extension of the
Cuntz algebra \cit{Cuntz} by the compact operators.
In \cit{QCR} we did not give explicit bounds on $T$. It might
appear that such bounds were best given in terms of the norm of the
operator $\Tt$ appearing naturally in \eq(In.1'), and the definition
of the algebra $\wick(T)$.
It is one of the main conclusions of this paper that this is
misleading, and that another operator, which we will simply denote
by $T$, plays a fundamental r\^ole. $T$ is a partial adjoint of
$\Tt$ with respect to only one tensor factor. Explicitly, we set
$T:\HH\otimes\HH\to\HH\otimes\HH$ with
$$ T e_i\otimes e_\ell =\sum_{} T_{ij}^{k\ell} e_j\otimes e_k
\quad.\tag\deq(In.2)$$
Our new, explicit criterion for the isomorphism of the universal
bounded \rep\ of  $\wick(T)$ with the extended Cuntz algebra
(\Thm\sBu.4/Bs.6/ below) becomes $\norm{T}<\sqrt2-1$, with further
refinements if $T\geq0$ or $T\leq0$. Our main results regarding the
Fock \rep\ are summarized in the following Theorem. The third item
was shown recently by \Bozejko\ and Speicher \cit{BoSpein} using
techniques rather different from those employed in this paper. In its
statement, $T_1=T\otimes\idty$, and $T_2=\idty\otimes T$, acting on
$\HH\otimes\HH\otimes\HH$.

\iproclaim/In.1/ Theorem.
Let $T_{ij}^{k\ell}$ denote the coefficients defining a \Wa\ on
finitely many generators, and let $T$ be the operator defined in
\eq(In.2). Then the Fock \rep\ of the \Wa\ $\wick(T)$ has
positive semi-definite scalar product, if any one of the following
conditions holds:
\roster
\item""{\kern-20pt(\Thm\sBs.2/Bs.2/)}\quad
        $\norm{T}\leq1/2$.
\item""{\kern-20pt(\Thm\sBp.1/Bp.1/)}\quad
        $T\geq0$.
\item""{\kern-20pt(\Thm\sBq.3/Bq.3/)}\quad
        $T$ satisfies the braid relation
          $T_1T_2T_1=T_2T_1T_2$, and $\norm{T}\leq1$.
\endroster
Moreover, all three bounds ($1/2,\ 0$, and $1$) are best possible.
\eproclaim

In the definition of \Wa s, we have imposed no relations between the
generators $a_i$. The idea of Wick ordering, however, makes sense also
in the presence of such relations, provided these relations are
compatible with the coefficients $T_{ij}^{k\ell}$ in a suitable sense.
Ideals in the tensor algebra over $\HH$ satisfying this condition will
be called {\em Wick ideals} (see Section \sI\ for details). An
interesting observation is that in some cases such an ideal is
annihilated by every bounded \rep.

As a case study of the interplay between the structure of Wick
ideals and the positivity of \rep s, we consider three structures
defined by Woronowicz and Pusz, namely their twisted canonical
commutation relations (\cit{PWor}, see \eq\sE.3(E.3)), the twisted
canonical anti-\-commutation relations (\cit{Pusz}, see
\eq\sE.4(E.4)), and the quantum group $\SnU2$ (\cit{WoRims}, see
\eq\sE.5(E.5)). In each case the original definition uses relations of
the form \eq(In.1), as well as relations between the $a_i$ alone. In
this paper we consider, in each case, the \Wa\ defined by retaining
only the \Wa\ relations, but discarding any further relations.

\iproclaim/In.2/ Theorem.
Consider an \irrep\ $\pi$ of the \Wa\ $\wick(T)$ by bounded
operators on a Hilbert space $\HH$, where $T$ is given by one of the
following:
\roster
\item
the twisted canonical commutation relations
(cf.\ \eq(E.3); \Thm\sIt.1/It.1/ )
\item
the twisted canonical anti-commutation relations
(cf.\ \eq(E.4); \Thm\sIt.3/It.3/ )
\item
$\SnU2$
(cf.\ \eq(E.5); \Prp\sE.4/E.4'/ )
\endroster
Then $\pi$ is coherent. More specifically there are a cyclic
vector $\Omega\in\HH$, $k\in I$, and $z\in\Cx$ such that
$a_k\Omega=z\Omega$, and $a_i\Omega=0$ for $i\neq k$.
\eproclaim

In the case of the twisted canonical commutation relations (resp.\
$\SnU2$) this result implies that the additional relations automatically
follow from the boundedness (resp.\ positivity) of the \rep; in the
case of the twisted canonical anti-commutation relations we also
find \rep s in which these relations are {\em not} satisfied.
In all cases, the parameters $k$ and $z$ defining the coherent \rep\
provide an efficient parametrization of the equivalence classes of
\irrep s, and we obtain a neat summary of the classification work
done in \cit{PWor,Pusz,Meister}, correcting, in passing, an omission
of some cases in the classification of \cit{PWor}.

Relations of the form \eq(In.1) appear naturally as the commutation
rules between coordinates and partial derivatives in {\em
non-commutative differential geometry}. We indicate the connections,
and point out a characterization of the braid relation for $T$ in
terms of differential geometry (\Thm\sDb.2/Db.2/). Finally, we study
the automorphism group acting on the generators as $a_j\mapsto
e^{-it}a_j$, and its associated KMS-states.

\subhead \ssecno\sIn.2. General background and related work \endsubhead

Many of the structures that have been studied in recent years as
deformations of ``classical'' objects, such as quantum groups
\cit{Woron,Faddeev}, and quantum planes \cit{Manin,WessZ}, as
well as various deformations of the canonical commutation relations
\cit{BoSpei,PWor} are defined in terms of algebraic relations
between the generators of an algebra. Consequently, much of the work
on these structures has been in the purely algebraic category. An
aspect that has often been neglected is the question under what
conditions a given set of relations can be realized by (possibly
unbounded) operators in a Hilbert space. This is unfortunate, both
from a physical and a mathematical point of view.

{}From the physical point of view, it is desirable to be able to
interpret the algebra under consideration as the algebra of
observables of a physical quantum system. Then the statistical
interpretation of quantum theory gives special importance to the
positive elements of the algebra, since only these can describe
probabilities. The positive elements are in turn defined as those of
the form $X^*X$, in terms of an involution $X\mapsto X^*$. Thus
observable algebras are always involutive algebras with a generating
positive cone. This in turn implies that the algebra has a faithful
\rep\  by Hilbert space operators. We do not claim, of
course, that these are the only algebras that may have ``physical''
significance. A good example to the contrary is the superselection
structure in low dimensional quantum field theory, often described
in terms of the \rep\ theory of some quantum group \cit{Rehren}, which
itself is not interpreted as an observable algebra. Nevertheless, we
feel that the observable interpretation of quantum theory is too
fundamental to be ignored lightly in the development of new algebraic
structures for quantum theory.

On the mathematical side the relevance of Hilbert space \rep s
(``positive definite'', or ``positive energy'' \rep s for emphasis
\cit{KacRaina,Pallebook,JOb}) is obvious from the example of Lie
algebra theory, where a rich structure unfolds even for finite
dimensional algebras. Picking a setting with \rep s realized in
Hilbert space demands positivity, and in a related way positivity
plays a basic r\^ole in \cit{Jones}. The algebras we consider in this
paper are all infinite dimensional (although typically with a finite
number of generators), so an additional aspect is the necessity of
topological concepts to tame the various infinite sums arising in
the theory. An example from quantum group theory for the need for
topological completions is the identification of the tensor product
of two function algebras with the algebra of functions on the
Cartesian products. This is usually taken to motivate the definition
of the coproduct, but generally fails for the {\em algebraic} tensor
product of algebras of continuous functions on non-discrete spaces.
The topological notions needed for such purposes in the
non-commutative situation are often imported most naturally from a
collection of Hilbert space \rep s.

In Sections \sB\ and \sI\ we show that there is a surprising
interplay between the purely algebraic side of the theory and the
constraints imposed by demanding a Hilbert space \rep. We show that
often there will be additional relations, which do not follow {\em
algebraically} from the given ones, which are nevertheless valid in
every Hilbert space \rep. This phenomenon is part of a general
theory of C*-algebras defined in terms of generators and relations,
which has been studied recently by several authors
\cit{Worfuncal,Blackadar,FNC,Phillips}.

The choice of the relations \eq(In.1) was made since several
examples of such structures are currently under investigation in the
literature. Our initial motivation was to understand the
relationship between two different deformations of the canonical
commutation relations, namely the ``$q$-deformed'' commutation
relations \cit{BoSpei,Green,QCS}, and the ``twisted''
canonical commutation relations \cit{PWor,Pusz}. The latter are
given as two sets of relations, one set of the form \eq(In.1), and
another set involving only the generators $a_i$. This is in contrast
to the situation for the $q$-deformed relations, where any such
additional relations lead to inconsistency. In order to understand
the relationship between these structures we will treat both cases
from the point of view of \eq(In.1) alone. Admissible relations
between the $a_i$ then correspond precisely to Wick ideals studied in
Section \sI.

The Wick order viewpoint is also a natural ingredient in the Wess-Zumino
\cit{WessZ}, and Woronowicz \cit{Worona}, approach to
non-commutative differential geometry, and we shall take up this
point in Section \sD\ below; see also Baez \cit{Baez} where the
classical Poincar\'e lemma from de Rham theory is proved in the
quantum theoretic setting. Here the adjoints $a_j$ of the generators
$a_j^*$ are interpreted as partial derivatives on the algebra of
``functions'' generated by the $a_j$ alone. This is motivated by the
Segal-Bargmann \cit{Folland} formalism for the complex formulation
of the Schr\"odinger \rep. The fundamental idea in
\cit{Bargmann} was to pick a positive definite inner product on
the functions of several complex variables $(z_1,\ldots,z_d)$ such
that the multiplication operator $f(z)\longmapsto z_j f(z)$ becomes
adjoint to ${\partial/\partial z_j}$, in the sense
$$\Big\langle z_j f(z),\, h(z)\Big\rangle
    =\Big\langle f(z),\, {\partial\over\partial z_j} h(z)\Big\rangle
\tag\deq(In.3)$$
for functions $f(z), h(z)$ in the corresponding Hilbert space, or
Bargmann space.  Specifically, Bargmann finds his (now familiar)
Gaussian Hilbert space as a solution $\bra\cdot,\cdot>$ to the
ansatz \eq(In.3). Our viewpoint here is that for the general Wick
relation \eq(In.1), the requirement of hermiticity, i.e.\ that
\rep s take the involution in the abstract \Wa\ into the
operator theoretic adjoint, serves as a generalization of \eq(In.3).
It turns out that this approach is equivalent to defining an algebra
of differential forms, with an additional commutation rule
transforming expressions of the form ``$x_i\dx j$'' into
``$\dx k\,x_\ell$''. For generic coefficients, however, all
differential forms of degree $2$ and higher vanish. In the presence
of Wick ideals non-trivial higher order forms can exist. In the
classical case the anticommutativity of the differentials implies
that any form of degree higher than the number of generators
vanishes. This is not true in the non-commutative setup discussed
here, and we will give an example of a calculus with forms of all
orders, using recent results from the theory of spin chains
\cit{FCS}.

For further examples of \Wa s, and more details on those mentioned
here, we refer to Section \sE.

The term ``{\em Wick ordering}'' for the process of ordering all
creation operators to the left of all annihilation operators is
standard usage in physics for the canonical (anti-) commutation
relations. We therefore use the term ``{\em \Wa}'' for an algebra
obtained form relations of the form \eq(In.1). The term has been
used before by Slowikowski \cit{Slowikow} in a much more special
(commutative) context.
The terms ``{\em coherent \rep}'' for the ``highest weight'' \rep s
with a cyclic vector $\Omega$ satisfying $a_i\Omega=\alpha_i\Omega$
with $\alpha_i\in\Cx$, and ``{\em Fock \rep}'' for the special case
$\alpha_i\equiv0$ are also in agreement with standard terminology
for the canonical (anti-) commutation relations. We give an
algebraic definition of these \rep s, and the canonical hermitian
forms on the respective \rep\ spaces,  for any choice of coefficients
$T_{ij}^{k\ell}$, and eigenvalues $\alpha_i$.
Our main problem, namely deciding the positivity of this hermitian
form is reminiscent of of the positivity question which is central to
the study of the coherent \rep\ of the {\em Virasoro algebras} $\cV$
(see \cit{KacRaina}); in the simplest case these \rep s are
labeled by two parameters $(c, \lambda)$ where $c$ is the central
charge, and $\lambda$ is the spectral value for the Cartan generator
of $\cV$. The coherent \rep\ is then defined by the highest weight
method and the \rep\  space is the corresponding Verma module
(see \cit{GoodWall,KacRaina}). The explicit range of the parameters
$(c,\lambda)$ where the corresponding inner product $\langle\cdot
,\cdot\rangle_{c,\lambda}$ is positive semi-definite is known, see
\cit{KacRaina} and \cit{GoodWall} for details. In a different
context, a similar positivity condition dictates the values of the
Jones index \cit{Jones,Jonesb}.

As \Thm/In.1/ shows, \Wa s with small $\norm{T}$ can be understood in
a rather general way. The universal bounded \rep\ exists by
\Thm\sBu.4/Bs.6/, and is isomorphic to the Cuntz-Toeplitz algebra,
which has only one proper closed two-sided ideal. Hence no relations
between the generators, in addition to \eq(In.1), can be imposed. When
$\norm{T}$ exceeds $1$ this isomorphism breaks down. As a simple
example in which additional relations are not only consistent with
\eq(In.1), but are forced by it in conjunction with the positivity
requirement, we present the C*-algebra of Woronowicz' quantum group
$\SnU2$. More generally, the possibility of further relations,
conveniently expressed as the existence of a non-trivial ``Wick
ideal'' implies that $-1$ is an eigenvalue of $T$, and the Fock scalar
product becomes partly degenerate. These Wick ideals are a typical
feature of the twisted canonical commutation relations
\cit{PWor,Pusz}, and also play a fundamental r\^ole in setting up a
differential calculus \cit{WessZ} (see below). Typically, Wick ideals
satisfy a version of the relations \eq(In.1) without constant term.
This allows us to conclude that the Wick ideal in twisted canonical
commutation relations \cit{PWor} is automatically annihilated in
every bounded \rep, and to describe the structure of \rep s of the
twisted canonical anti-commutation relations not assuming relations
between the $a_i$. The ``untwisted'' case is the theory of Clifford
algebras, and was treated from this point of view in \cit{QCS}.

Finally, we briefly consider the KMS states associated with the
one parameter automorphism group $\gauge$ acting on the generators
as $\gauge_t(a_j)=e^{-it}a_j$. In Fock space this automorphism is
generated by the number operator, and we connect the growth of the
$n$-particle subspaces to the existence of KMS states.

Explicit computations with Wick relations can be quite painful,
since each application of the rule \eq(In.1) may multiply the number
of terms involved. For the present paper the largest computations of
this sort were necessary to verify the Wick ideal property in
\Thm\sIt.4/It.4/. We therefore developed a little {\smc Mathematica}
\cit{Mathematica} package to perform such computations. It takes the
commutation rules as a parameter, and computes Wick ordered forms,
as well as Fock and coherent functionals. We have made it available
by anonymous {\tt ftp} from
``{\tt nostromo.physik.Uni-Osnabrueck.de}''
in the directory ``{\tt pub/Qrelations/Wick}''.

\head    \secno\sH.  Basic Definitions and Examples \endhead
\subhead \ssecno\sHi. Hermitian Wick algebras  \endsubhead

In this section we supply the general facts on the algebra
$\wick(T)$ based on the relation \eq(In.1). This will be the initial
axiom. We found it convenient to write the generators not as indexed
quantities ``$a_i$'', with redundant letters ``$a$'' cluttering every
computation, but to denote the generators directly by ``$i$''.
Moreover, we changed the position of the stars relative to \eq(In.1).
The connection with the notation of \eq(In.1) will be made explicit
once more in \eq\sHi.6(H.5) below.

\iproclaim/H.1/ Definition.
Let $I$ be a set, and let $T_{ij}^{k\ell}\in\Cx$ for
$i,j,k,\ell\in I$ such that for each pair $i,j$ only finitely many
$T_{ij}^{k\ell}\neq0$. Then the {\bf\Wa} on the coefficients $T$,
denoted by $\wick(T)$, is the algebra of polynomials (possibly with
constant term) in the symbols
$i, i\1$ for $i\in I$, divided by the relation
$$ i\1j=\delta_{ij}\idty
         + \sum_{k,\ell\in I} T_{ij}^{k\ell}\ \ell k\1
\quad,\tag\deq(H.1)$$
A {\bf\HWa} is a \Wa\  with involution ``$\ \1$'' extending
the map $i\mapsto i\1$. The existence of this involution is
equivalent to the condition
$$ T_{ji}^{\ell k}=\Br{T_{ij}^{k\ell}}
\quad.\tag\deq(H.2)$$
\eproclaim

The distinguishing feature of the relations \eq(H.1) is
that they allow us to order any polynomial in the generators $i,i\1$
such that in every monomial all $i$ are to the left of all $i\1$. The
result of this transformation is called the {\em Wick ordered form} of
the given polynomial. It is unique, i.e.\ independent of the order in
which the rules are applied to different parts of the expression.
The proof given in the following Lemma is based on a simple graph
theoretical principle known as the Diamond Lemma \cit{Bergmann}.

\iproclaim/H.2/ Lemma.
Let $\wick(T)$ be a \Wa. Then the Wick ordered monomials form a basis
of $\wick(T)$.
\eproclaim

\proof:
In order to apply the theory of algebraic reduction presented in
\cit{Bergmann} we have to verify two things: firstly, the
successive application of Wick ordering substitutions must
terminate, which is obvious. Secondly, the so-called diamond
condition has to be satisfied, i.e.\ whenever there are two
different possibilities for applying the substitution rules to the
same expression, there must be further substitutions on the
respective results leading to the same final result. In a \Wa, a
monomial in which two different substitutions are possible is of the
form $Xi\1jYk\1\ell Z$, where $i,j,k,\ell$ are generators, and
$X,Y,Z$ arbitrary (possibly empty) monomials. Since the two possible
substitutions in this situation commute, the diamond condition is
also obviously satisfied.
\QED

By $\HH$ we will denote the free vector space over $I$, i.e.\ the
set of finite linear combinations $f=\sum_{i\in I}f_i\,i$. Then
$f\1=\sum_{i\in I}\Br{f_i}\,i\1$ is considered in a natural way as
an element of the complex conjugate space $\HH\1$ of $\HH$.
The canonical inner product of $\HH$, denoted by
$$ \bra f,g>=\sum_i \Br{f_i}g_i
\tag\deq(H.3)$$
makes $\HH$ into a pre-Hilbert space. The free algebra over the
generators $i,i\1$ is the same as the tensor algebra
$$\tA(\HH\1,\HH)=\Cx\idty\oplus\HH\oplus\HH\1\oplus(\HH\otimes\HH)
               \oplus(\HH\otimes\HH\1)\oplus(\HH\1\otimes\HH)
              \oplus(\HH\1\otimes\HH\1)\oplus\cdots
\tag\deq(H.tA)$$
over $\HH$ and $\HH\1$, i.e.\ the direct sum of all tensor products
with factors $\HH$ or $\HH\1$, with tensor product as multiplication.
The empty tensor product is identified with $\Cx$, and this summand
of $\tA(\HH\1,\HH)$ contains the unit $\idty$ of the algebra.
We will often omit the tensor product signs, using them only for
punctuation improving the legibility of a formula.
The \Wa\ $\wick(T)$ is then the quotient of $\tA(\HH\1,\HH)$ by the
ideal in $\tA(\HH\1,\HH)$, generated by
$f\1\otimes g-\bra f,g>\idty - \Tt(f\1\otimes g)$, where $\Tt$
denotes the linear operator
$$ \Tt:\HH\1\otimes\HH\to\HH\otimes\HH\1\quad
      :f\1\otimes g\ \mapsto
         \sum_{i,j,k,\ell} T_{ij}^{k\ell}\ \Br{f_i}g_j\
             \ell\otimes k\1
\quad.\tag\deq(H.4)$$
Note that using the operator $\Tt$ we can state the relations
\eq(H.1) in basis-free form. This will be useful for extending the
above structures to suitable completions of $\HH$ and $\wick(T)$
(compare \Exa\sE.8/E.7/).

We are mainly interested in \rep s of a \Wa\ by linear
operators on a complex vector space $\R$. In many examples the
representatives of $i$ and $i\1$ will have an interpretation as
creation and annihilation operators. In such cases we sometimes
adhere to the curious, but universally accepted convention to make
$f\mapsto a\1(f)$ a linear operator, and $f\mapsto a(f)$ antilinear.
In other words, we denote the \rep\ by $a\1:\wick(T)\to\R$,
and use the notational conventions
$$\matrix
   \hfill                a\1_i &=& a\1(i)
\hfill\quad&,\quad&\hfill
                           a_i&=& a\1(i\1)
\hfill\\ \vtop to2pt{}\\ \hfill
                        a\1(f)&=&\sum_if_ia\1_i
\hfill\quad&,\quad&\hfill
                          a(f)&=&a\1(f\1)=\sum_i\Br{f_i}a_i
\quad.\hfill\\ \endmatrix\tag\deq(H.5)$$
The following Definiton lists the properties of \rep s in
which our main interest lies.

\iproclaim/H.3/ Definition.
A {\bf Hermitian \rep} of a \HWa\ $\wick(T)$ is a
\rep\ $\pi:\wick(T)\to\LL(\R)$ by linear operators on a
complex vector space with Hermitian form
$\iprod:\R\times\R\to\Cx$ such that for all
$\xi,\eta\in\R$, and $i\in I$:
$$ \bra \xi,\pi(i)\eta>=\bra \pi(i\1)\xi,\eta>
\quad.\tag\deq(H.6)$$

A \rep\ $\pi$ is called {\bf positive}, if the form on the
\rep\ space $\R$ is positive definite. It is called {\bf
bounded}, if it is positive, and $\pi(i)$ is bounded for all $i\in
I$. It is called {\bf \cbo} with bound $\beta$, if it is positive, and
for all $\xi\in\R$:
$$ \sum _i  \norm{\pi(i\1)\xi}^2
      \leq \beta\ \bra\xi,\xi>
\quad.\tag\deq(H.7)$$
\eproclaim

Note that even when the hermitian form of $\R$ is positive definite
we have not assumed $\R$ to be complete, i.e.\ a Hilbert space. If
$\widetilde\R$ denotes the completion of $(\R,\bra\cdot,\cdot>)$, we
may consider $\pi(X)$ for $X\in\wick(T)$ as an operator from the
domain $\R\subset\widetilde\R$ into $\widetilde\R$. Each $\pi(X)$
will be closable with adjoint $\pi(X)^*$ relative to the scalar
product of $\widetilde\R$, and and we have the operator inclusion
$\pi(X\1)\subset \pi(X)^*$. Note that $\R$ is a common, dense,
invariant domain for all the $\pi(X)$.
Given a hermitian linear functional $\omega:\wick(T)\to\Cx$, i.e.\ a
functional with $\omega(X\1)=\Br{\omega(X)}$, we get a hermitian
form on $\wick(T)$ by $\bra X,Y>_\omega=\omega(X\1Y)$. The \rep\
acting on $\wick(T)$ by left multiplication becomes hermitian with
respect to this form. We call it the {\em cyclic}, or GNS- {\em\rep}
$\pi_\omega$ associated with $\omega$. Obviously, this \rep\ is
positive if and only if $\omega$ is, i.e.\ iff $\omega(X\1X)\geq0$.
Conversely, every vector in a positive \rep\ defines a positive
functional. Non-zero positive functionals need not exist on a Wick
algebra, since the cone of elements $\sum X_i\1X_i$ may contain all
elements of the form ``$-Y\1Y$''.
Other related such ``singular features'' are at the heart of the \rep\
theoretic problems for non-compact quantum groups, see e.g.\
\cit{Woronb} and \cit{vDaele}.

When a \rep\ $\pi$ is bounded, we may study the C*-algebra generated
by $\pi\bigl(\wick(T)\bigr)$ in $\BB(\R)$. It is clear that \rep s of
this algebra in turn produce bounded \rep s  of $\wick(T)$.
Therefore we may ask whether every bounded \rep\ of $\wick(T)$
arises in this way from a single \rep. The following definition
describes this situation.

\iproclaim/H.4/ Definition.
The {\bf universal bounded \rep} of a \Wa\  $\wick(T)$ is a
C*-algebra, denoted $\wicc(T)$, which is generated by elements
$i\in I$ satisfying the relations \eq(H.1), with the following
universal property: If $\pi:\wick(T)\to\BB(\R)$ is a bounded
\rep\  on a Hilbert space $\R$, there is a unique C*-algebra
\rep\  $\hat\pi:\wicc(T)\to\BB(\R)$ such that
$\pi(i)=\hat\pi(i)$. The universal {\bf\cbo\ \rep} is defined
analogously.
\eproclaim

The definition fits the pattern of a C*-algebra presented by
generators and relations \cit{Blackadar,Worfuncal,GPOTS,FNC}. The
existence of a universal bounded \rep\ is equivalent to the
following two conditions: firstly, there must be {\em some} bounded
\rep, and secondly, for each generator $i\in I$, there must be a
number $\beta_i<\infty$ such that the uniform bound
$\norm{\pi(i)}\leq\beta_i$ holds in {\em every} bounded \rep. If
only the first condition holds, one can still find a universal \rep\
in the class of ``inverse limits'' of C*-algebras \cit{Phillips}, or
LMC*-algebras \cit{Schmuedgen}. Since in all the examples we study
below either both conditions hold or none of them, we will not
elaborate on this structure.

\subhead \ssecno\sE. Examples \endsubhead

In this section we describe the principal examples, and classes of
examples from the literature which motivated our study. We will
repeatedly refer to these examples in the main parts of the paper,
supplementing them by further examples illustrating more technical
points. The descriptions given in this section are accompanied by
brief indications of the new results we have obtained in each case,
and where to find them in this paper. The coherent \rep s
$\crep\pfi$, which will be defined explicitly in the next section,
are defined by $\crep\pfi(i\1)\Omega=\pfi_i\Omega$, where
$\pfi_i\in\Cx$, and $\Omega$ is a cyclic vector. For $\pfi_i\equiv0$
we get the Fock \rep\ $\crep0$.

\examp/E.1/ $q$-canonical commutation relations:
We fix a real constant $q$, and set
$$    i\1j=\delta_{ij}\idty + q\, ji\1
\quad,\tag\deq(E.1)$$
where the index set $I\ni i,j$ is arbitrary.
Usually these relations are written in terms of ``creation'' operators
$a\1_j$ and ``annihilation operators'' $a_i$ as
$$ a_ia_j\1=\delta_{ij}\idty + q\, a_j\1a_i
\quad,\tag\deq\lasteq'(E.1')$$
For $q=1$ we have the canonical commutation relations for Bosons,
for $q=-1$, we have the canonical anticommutation relations for
Fermions. The relations \eq(E.1') have been suggested by Greenberg
\cit{Green} as an interpolation between Bose and Fermi statistics.
Independently they were introduced by \Bozejko\ and Speicher
\cit{BoSpei}, who showed Fock positivity in this case (see also
\cit{Zagier,Fivel,Dykema}). In the range $\abs q<\sqrt2-1$ the
universal bounded \rep\ was determined in \cit{QCR}. Coherent states
were studied in \cit{QCS}, generalizing earlier work \cit{JOa} on
the case $q=0$.

\examp/E.2/ Temperley-Lieb-Wick relations:
Superficially these relations look like the $q$-canonical
commutation relations, but differ in the positioning of the indices
in the term containing $q$:
$$    i\1j=\delta_{ij}\idty + q\, ij\1
\quad,\tag\deq(E.2)$$
for $i,j=1,\ldots,d$.
These relations were first studied in \cit{QCR}. We constructed the
universal bounded \rep\ for $\abs qd<(\sqrt5-1)/2$, and showed the
independence of the universal C*-algebra of $q$. Apart from the
range $q\geq-1/d$ positive \rep s can only exist for a discrete
series characterized by $q^{2N}=(1+qd)/\bigl(q(d+q)\bigr)$,
$N\in\Nl$. One of the remarkable features of these relations in
comparison with \Exa/E.1/ is the sensitive dependence of the \rep\
theory on the number $d$ of generators. The reason behind this will
become completely clear in Section \sBe. The connection with the
Temperley-Lieb algebra \cit{Temperley,Jones} will also be explained
there, and we show positivity of the Fock \rep\ for $q>-1/(2d)$.

\examp/E.3/ Twisted canonical (anti-)commutation relations:
These relations have been introduced by Pusz and Woronowicz
\cit{PWor,Pusz}.
The Bosonic version \cit{PWor} depends on a parameter
$0<\mu<1$, and is given by
$$
i\1j=\cases \mu\ ji\1   &\text{ for $i\neq j$}\\
                  \idty+\mu^2\ ji\1
                      -(1-\mu^2)\sum_{k<i}kk\1
                                &\text{ for $i=j$.}\\
\endcases
\tag\deq(E.3)$$
We caution the reader that we have changed the ordering of the
indices relative to \cit{PWor}. This will facilitate the discussion
of the relations with infinitely many generators. A typical feature of
these relations is that, in the Fock \rep\  $\creps0$, one has
$\creps0(ij-\mu ji)=0$ for $i>j$. In the work \cit{PWor} these
relations are additional postulates. Here, in \Thm\sIt.1/It.1/, we
show that the they are true in every bounded \rep, and
moreover, that all bounded \irrep s are coherent.

The Fermionic version of twisted commutation relations in this sense
was discussed by \cit{Pusz}, and is given by
$$
i\1j= \cases -\mu\ ji\1   &\text{ for $i\neq j$}\\
             \idty - ji\1 -(1-\mu^2)\sum_{k<i}kk\1
                          &\text{ for $i=j$.}\\
\endcases
\tag\deq(E.4)$$
Here the additional relations which turn out to be satisfied in the
Fock \rep\  are $\creps0(i)^2=0$, and
$\creps0(ij+\mu ji)=0$ for $i>j$. In Section \sIt\ we investigate to
what extent these relations are also consequences of positivity, like
their Bosonic counterparts. We find that there are \rep s not
satisfying these relations, and parametrize all possibilities by
showing that all \irrep s must be coherent (see \Thm\sIt.3/It.3/).

\examp/E.4/ The quantum group $\SnU2$:
Some of the quantum groups of Woronowicz also can be considered as
Hermitian \Wa s with additional relations. For example,
$\SnU2$ is given by the relations \cit{WoRims}
$$\align
  \left.\aligned
    \alpha\alpha^*&=\idty-\nu^2\gamma^*\gamma   \\
    \gamma\gamma^*&=\idty- \alpha^*\alpha       \\
    \alpha\gamma^*&= \nu \gamma^*\alpha
  \endaligned\right\rbrace \tag{\deq(E.5)}\\
  \left.\aligned
    \gamma\gamma^*&=\gamma^*\gamma    \hskip28pt\\
    \alpha\gamma  &= \nu \gamma\alpha
  \endaligned\right\rbrace \tag{\deq(E.5')}
\endalign$$
The separation into \eq(E.5) and \eq(E.5') is not made in
\cit{WoRims}; it is suggested in our context since only the
relations \eq(E.5) fit the pattern of \Wa\ commutation rules (the
generators are $\alpha^*$ and $\gamma^*$). It turns out that this is
a fruitful way of looking at $\SnU2$, since the two remaining
relations are redundant. The following Proposition also shows how
coherent \rep s may help to obtain quickly all \irrep s of the
C*-algebra $\SnU2$.

\global\advance\pcNo by -1
\iproclaim/E.4'/ Proposition.
Let $\alpha$ and $\gamma$ be operators on a Hilbert space,
defined on a dense domain invariant under $\alpha,\alpha^*,\gamma$,
and $\gamma^*$, on which the relations \eq(E.5) hold. Then $\alpha$
and $\gamma$ are bounded, and if $\nu\neq0$, relations \eq(E.5') also
hold.
Moreover, all \irrep s are coherent in the following sense:
if $\set{\alpha,\gamma,\alpha^*,\gamma^*}$ is an irreducible set of
bounded operators, the \rep\ space contains a
vector $\Omega$ such that $\alpha\Omega=\hat\alpha\Omega$, and
$\gamma\Omega=\hat\gamma\Omega$, where
$\hat\alpha,\hat\gamma\in\Cx$, and
$$ \abs{\hat\alpha}=1, \hat\gamma=0 \qquad,\quad\text{or}\quad
        \hat\alpha=0, \abs{\hat\gamma}=1
\quad.\tag\deq(E.5Prp)$$
\eproclaim

\proof:
Note first that from either of the first two relations we
get that $\alpha$ and $\gamma$ are both bounded in every positive
\rep. Consider the operators
$$\aligned
   C&=\alpha\gamma-\nu\gamma\alpha        \\
\text{and}\qquad
   R&=\idty-\alpha^*\alpha-\gamma^*\gamma
\quad.\endaligned\tag\deq(E.5p)$$
Then we have to show that $C=R=0$.
A straightforward computation using the Wick relations \eq(E.5) gives
\eqgroup()$$\align
      C^*C&= R(\idty-R)
\quad\text{and}     \tag{\deq\lasteq a(E.5c+c)}\\
      CC^*&= -\nu^2 R(\idty+\nu^2 R)
\quad. \tag{\deq\lasteq b(E.5cc+)}\endalign$$
The since $R=R^*$ \eq(E.5c+c) forces $0\leq R\leq1$. Hence, by
\eq(E.5cc+), $CC^*\leq0$, and $C=0$. Then, from \eq(E.5c+c), $R$ has
spectrum $\set{0,1}$, i.e.\ $R$ is a projection. When $\nu\neq0$,
\eq(E.5cc+) implies $R=0$, thus proving the first
statement. In the case $\nu=0$ we can only conclude with \eq(E.5p)
that $\alpha$ and $\gamma$ are partial isometries with
$\idty=\alpha\alpha^* \geq\idty-\alpha^*\alpha
      =\gamma\gamma^* \geq\gamma^*\gamma$.

It follows that the universal bounded \rep\ of \eq(E.5) is precisely
the C*-algebra of $\SnU2$, as defined by Woronowicz.
We compute all positive coherent \rep s of
this \Wa: Suppose that  $\omega(X\alpha)=\hat\alpha\omega(X)$, and
$\omega(X\gamma)=\hat\gamma\omega(X)$ for some complex constants
$\hat\alpha,\hat\gamma$. Then positivity of $\omega$ requires
$\omega(C)=\omega(R)=0$, and hence $\hat\alpha\hat\gamma=0$, and
$\vert\hat\alpha\vert^2 +\vert\hat\gamma\vert^2 =1$.
These conditions are also sufficient for $\omega$ to be positive:
when $\hat\gamma=0$, and $\hat\alpha$ has modulus $1$, the \rep\ is
one dimensional with $\alpha=\hat\alpha\idty$, and $\gamma=0$.
On the other hand, when $\hat\alpha=0$, and $\hat\gamma$ has
modulus $1$, we first construct the Fock \rep\  of the \Wa\
determined by $\alpha\alpha^*=(1-\nu^2)\idty+\nu^2\alpha^*\alpha$,
with a vacuum vector $\Omega$ satisfying $\alpha\Omega=0$.
Then we explicitly define
$$ \gamma\ \alpha^{*n}\Omega
       = \hat\gamma\,\nu^n\ \alpha^{*n}\Omega
\quad.$$
One readily verifies that this is the coherent \rep, and, comparing
with the known \irrep s of the C*-algebra of $\SnU2$
\cit{Meister} we find the result.
\QED

\examp/E.5/ Braid relations:
Examples \pcl/E.1/ and \pcl/E.3/ have a common feature, namely the
validity of the identity
$$ \sum_{efh} T_{be}^{fc}\ T_{ak}^{he}\ T_{h\ell}^{mf}
  =\sum_{deh} T_{ad}^{eb}\ T_{eh}^{mc}\ T_{dk}^{\ell h}
\qquad,\quad\forall\ abc,k\ell m\in I
\quad,\tag\deq(E.6)$$
where the index set $I$ labelling the generators is finite.
In Section \sBe\ we will identify this as the familiar relation
determining the braid group.
This turns out to be a useful assumption on $T_{ij}^{k\ell}$ in two
very different contexts: it allows an improvement of the bounds insuring
Fock positivity (see Section \sBq), and also allows the extension of
the differential calculus of Section \sDf\ from the tensor algebra
$\tA(\HH)$ to the involutive algebra $\wick(T)$ (see Section \sDb).
The simplest case are the relations
$$    i\1j=\delta_{ij}\idty + q_{ij}\, ji\1
\quad,\tag\deq(E.7)$$
with $q_{ij}\in\Cx$.
When $q_{ij}$ takes only values $\pm1$, but the indices cannot be
grouped into ``bosonic'' and ``fermionic'' (which would be
equivalent to $q_{ij}=(1+q_i+q_j-q_iq_j)/2$, $q_i=\pm1$) one speaks
of anomalous statistics. Speicher \cit{Speix} has used the known
Fock positivity in this case to show it also for general real
$q_{ij}$ with $\abs{q_{ij}}\leq1$. Below in \Thm\sBq.2/Bq.2/ we will
extend this Fock positivity result to complex $q_{ij}$ with
$\abs{q_{ij}}\leq1$. An independent proof generalizing to all $T$
satisfying equation \eq(E.6) was given by \Bozejko\ and Speicher
\cit{BoSpein} (see Section \sBq). We also remark that the braid
condition \eq(E.6) is assumed in the work of Baez \cit{Baez}.

\examp/E.6/ Clifford algebras:
Clifford algebras are defined by the relations
$i\1j+ji\1=\delta_{ij}\idty$, $i,j=1,\ldots,d\leq\infty$, and hence
correspond to the limit $q\to-1$ in \eq(E.1), or $\mu\to1$ in
\eq(E.4). The Fock \rep\  of this algebra is identical with
the Fermi relations, i.e.\ the generators anticommute among
themselves, not only with their adjoints. In a coherent \rep\ one
has $\crep\pfi(ij+ji)=2\pfi_i\pfi_j\idty$. In the universal bounded
\rep\ the elements $\theta_{ij}=ij+ji$ are central. Hence, in every
irreducible \rep, $\theta_{ij}$ is a multiple of the identity. It
turns out \cit{QCS} that an irreducible \rep\ with given matrix
$\theta$ exists if and only if $\norm{\theta}\leq1$ as Hilbert space
operator on $\Cx^d$. Moreover, for each $\theta$ there are at most
two different irreducible \rep s, which are finite dimensional if
$d<\infty$.

\examp/E.6A/ A degenerate case:
We consider the relations
$$ i\1j=\delta_{ij}\bigl(\idty-\sum_{k\in I}kk\1\bigr)
\quad,\tag\deq(E.deg)$$
where $I=\set{1,\ldots,d}$ is a finite set. This relation has
special properties with respect to several of the questions we will
discuss below: it is the only choice of coefficients for which the
Fock scalar product becomes completely degenerate for two or more
``particles''; it allows the largest possible quadratic Wick ideal
(see Section \sIq), and it is the only choice for which all products
of coordinate differentials in the calculus of Section \sDf\ are
linearly independent. Moreover, it is in the two parameter family of
commutation relations
$$ i\1j=\delta_{ij}\idty
         +q\,ji\1
         -\lambda \delta_{ij}\sum_{k\in I}kk\1
\quad,\tag\deq(E.usym)$$
which is characterized by the full symmetry with respect to unitary
transformations of $\HH$ \cit{Shale}.
Here we will briefly consider the \rep\ theory of \eq(E.deg), using
direct methods.

\global\advance\pcNo by -1
\iproclaim/E.deg/ Theorem.
Let $\pi$ be an \irrep\ of the relations \eq(E.deg) in a Hilbert
space $\R$, and let $R=\idty-\sum_iii\1\in\wick(T)$.
Then there is a
number $c_\pi$, $0\leq c_\pi\leq1/4$ such that
$$ \pi\bigl(R(\idty-R)\bigr)=c_\pi\idty
\quad.$$
\roster
\item
For $c_\pi=0$, $\pi$ has to be the Fock \rep, for which $\dim\R=d+1$
\item
For $c_\pi=1/4$, we have $\pi(i)=2^{-1/2}\,v_i$, where the $v_i,
i=1,\dots,d$ satisfy the Cuntz relations \cit{Cuntz}
$$ v_i^*v_j=\delta_{ij}\idty
\quad,\quad\text{and}\quad
    \sum_iv_iv_i^*=\idty
\quad.$$
\item
For $0<c_\pi<1/4$, there is an \irrep\ $w_{ij}$, $i,j=1\dots,d$ of
the Cuntz relations on $d^2$ generators in a Hilbert space $\R_0$
such that, up to a unitary equivalence, $\R$ is the direct sum of
$(d+1)$ copies of $\R_0$, and the generators are given by the
following $\BB(\R_0)$-valued $(d+1)\times(d+1)$-matrix, which we
write in terms of the parameters $\alpha,\beta$ with
$\alpha^2=1/2+(1/4-c_\pi)^{1/2}$, and
$\beta^2=1/2-(1/4-c_\pi)^{1/2}$:
$$ \pi(i)= \pmatrix 0&     0&\dots\beta\dots&0\\
        \alpha w_{i1}&     0&    \dots      &0\\
               \vdots&\vdots&   \ddots      &\vdots\\
        \alpha w_{id}&     0&    \dots      &0 \endpmatrix
\quad,\tag\deq()$$
where $\beta$ in the first row is in the $(i+1)$\th\ position.
\endroster
Moreover, the coherent \rep\ $\crep\pfi$ with cyclic vector $\Omega$
satisfying $\crep\pfi(i\1)\Omega=\pfi_i\Omega$ is positive if and only
if
$$ \bigl(\textstyle\sum_i\abs{\pfi_i}^2\bigr)^2
    = c_\pi\leq {1\over4}
\quad.$$
\eproclaim

\proof:
One readily checks that $iR+Ri=i$, and, by taking adjoints
$i\1R+Ri\1=i\1$. Hence by induction on the degree of polynomials we
get $XR=(\idty-R)X$ for all odd polynomials in the generators
$i,i\1$, and $RY=YR$ for all even ones. It follows that $R(\idty-R)$
is in the center of $\wick(T)$. Hence $\pi(R(\idty-R))$ is a
multiple of the identity in every \irrep\ $\pi$.

(1) By Wick ordering we get $\pi(R(\idty-R))=
\sum_{ij}\pi(ij)\pi(ij)^*$. Hence $c_\pi=0$ implies $\pi(ij)=0$ for
all $i,j$. Hence only the Wick ordered monomials of the forms
$\idty, i,i\1, ij\1$ can be non-zero. Hence
$\dim\pi(\wick(T))\leq(d+1)^2$. Moreover,
$$ \sum_i\pi(i\1R)^*\pi(i\1R)=\pi\bigl(R(\idty-R)R\bigr)=0
\quad.$$
Hence $\pi(i\1R)=0$, and $\pi(R)$ is the projection onto the
one-dimensional subspace of Fock vectors.

(2) For $c_\pi=1/4$, we have $\pi(R)=(1/2)\idty$, and the relations
\eq(E.deg) become the Cuntz relations up to a trivial change in
normalization.

(3) $\pi(R)$ has the eigenvalues $\alpha^2$ and $\beta^2$. Let $P$
denote the eigenprojection for eigenvalue $\alpha^2$, and
$\R_0=P\R$. Then, since $iR=(\idty-R)i$, $\pi(i)$ swaps the
eigenspaces of $R$, i.e.\ $\pi(i)=u_i+v_i$, where
$u_i=(\idty-P)\pi(i)P:\R_0\to\R_0^\perp$, and
$v_i=P\pi(i)(\idty-P):\R_0^\perp\to\R_0$.
In terms of these operators the relations \eq(E.deg) become
$u_i^*u_j=\alpha^2\,\delta_{ij}\,P$, and
$v_i^*v_j=\beta^2\,\delta_{ij}\,(\idty-P)$, and the definition of
$R$ becomes $\sum_iu_iu_i^*=\alpha^2(\idty-P)$, and
$\sum_iv_iv_i^*=\beta^2\,P$. Using the partial isometries
$\alpha^{-1}u_j$, we can identify $\R_0^\perp=(\idty-P)\R$ with the
direct sum of $d$ copies of $\R_0$.
$$ \R=\R_0\oplus\R_0^\perp
     = \R_0\oplus\bigoplus_{i=1}^d \alpha^{-1}u_j\R_0
     \cong \bigoplus_{i=0}^d \R_0
\quad.$$
In this structure we have encoded the action of $u_i$. The action of
$v_i$ is given by the isometries
$$ w_{ij}=(\beta^{-1}v_i)(\alpha^{-1}u_j)
         =c_\pi^{-1/2}P\pi(ij)P
\quad.$$
One readily verifies that these now satisfy the Cuntz relations, and
that $\pi(i)=u_i+v_i$ is given by the matrix shown in the
Theorem. To see the irreducibility of $w_{ij}$, suppose that
$X\in\BB(\R_0)$ commutes with the $w_{ij}$ and their adjoints. Then
$\widehat X=PXP+\alpha^{-2}\sum_i u_iXu_i^*$ commutes with all
$\pi(i),\pi(i)^*$. Hence $\widehat X$ is a multiple of $\idty$, and
so is $X$.

The $\R_0$-component of the cyclic vector $\Omega$ has to be
coherent \cit{JOa,QCS} for the Cuntz algebra generators $w_{ij}$.
The other components have to multiples of the $\R_0$-component. It
is known for which eigenvalues such coherent states exist \cit{JOa},
and this is translated into the statement given in the Theorem.
\QED

\examp/E.7/ Infinitely many generators:
\Def/H.1/ allows the set $I$ labelling the generators to be
infinite with the proviso that the sums in the basic relations
\eq(H.1) are always finite. This assumption was sufficient to make
the definitions of the previous section, as well as the definition of
Fock and coherent \rep s meaningful. We saw that a
\Wa\ can also be thought of as being generated by the vector
space $\HH$, of which $I$ is a linear basis. However, in some of the
applications it is desirable to take $\HH$ as a (possibly infinite
dimensional) Hilbert space, with the set $I$ as an orthonormal
basis. In that case we will take the sums in \eq(H.3) to be
convergent in the norm of $\HH$. All tensor products are to be read
as Hilbert space tensor products, which makes the tensor algebra
$\tA(\HH\1,\HH)$ an involutive Banach algebra. We next have to make
sense of the operator $\Tt$ in \eq(H.4). The assumption we have made
on its matrix elements $T_{ij}^{k\ell}$ are {\em not} sufficient to
guarantee that $\Tt$ is bounded. This is illustrated by \Exa/E.3/,
where $\nTt$ grows like $(1-\mu^2)\sqrt d$ for the relations with
$d<\infty$ generators. However, if we assume for the moment that
$\nTt<\infty$, the operator
$f\1\otimes g\mapsto \bra f,g>\idty + \Tt(f\1\otimes g)$ is also
bounded, and its graph generates an ideal $\ideal$, such that
$\wick(T):=\tA(\HH\1,\HH)/\ideal$ is a Banach algebra with
involution. Of course, as for finitely many generators, this algebra
may fail to have non-zero positive linear functionals.

\examp/E.8/ Relations with small coefficients $T$:
In this case the relations can be considered as a small perturbation
\cit{QCR} of the relations with $T=0$. The algebra $\wick(0)$ has a
universal bounded \rep\  $\wicc(0)$, known as the
Cuntz-Toeplitz algebra. It has a faithful \rep\  as the
algebra of creation and annihilation operators in the full Fock space
over $\HH=\ell^2(I)$, i.e.\ in the direct sum of the $n$-fold
(unsymmetrized) tensor products of $\HH$ with itself. We showed in
\cit{QCR} that for sufficiently small $T$, $\wicc(T)$ exists, and is
canonically isomorphic to $\wicc(0)$. From the form of this
isomorphism we also get Fock positivity, and the faithfulness of the
Fock \rep. The precise meaning of ``small $T$'' in this
context will be given below in Section \sBu.

\examp/E.9/ Counterexamples:
Here we collect a few elementary examples to demonstrate some of the
phenomena that can occur in \Wa\ theory.

\subex1 A \Wa\ without any positive \rep s
\Take two generators $i,j$ with $i\1i=j\1j=\idty$, and
$i\1j=\mu ii\1$. The first two relations make $i,j$ non-zero
co-isometries in every positive \rep, i.e.\ $\norm{i}=\norm{j}=1$.
Then the last equation gives a contradiction as soon as $\abs\mu>1$.
One can also see this in a more algebraic way by observing that
$$ i\1(i\1-j\1)(i-j)i + (\mu-1)\idty=0
\quad.$$
Hence, for $\mu>1$, we find that $\idty$ is both positive and
negative, and consequently vanishes in every positive \rep.

\subex2 A \Wa\ with positive, but without bounded \rep s
\Take the canonical commutation relations (\Exa/E.1/ with $q=1$).

\subex3 A \Wa\ with a universal bounded \rep, whose Fock \rep\ is not
positive.
\Take
Woronowicz' $\SnU2$ as described in \Exa/E.4/. The universal
bounded \rep\ is Woronowicz' C*-algebra $\SnU2$. On the other hand,
the element $C=\alpha\gamma-\nu\gamma\alpha$ gives a negative Fock
expectation $\cst0(CC^*)=-\nu^2(1+\nu^2)$.

\subex4 A \Wa\ for which the family of positive coherent  \rep s is
not faithful.
\Take the Clifford algebras of \Exa/E.5/.
With $\theta_{ij}=ij+ji$, we have
$\crep\pfi(\theta_{ij}\theta_{k\ell}-\theta_{ik}\theta_{j\ell})=0$
in any coherent \rep\ $\crep\pfi$. But since there are \irrep s in
which $\theta$ is an arbitrary (but small) scalar symmetric matrix
\cit{QCS}, this identity does not hold in $\wick(T)$.

\subhead \ssecno\sHc. Coherent \rep s\endsubhead

There is an important set of \rep s of a \Wa\ $\wick(T)$ which is
intimately connected with the Wick ordering process:
Let $\tA(\HH)=\Cx\idty\oplus\HH\oplus(\HH\otimes\HH)\cdots$ denote the
tensor algebra over $\HH$.
Let
$\pfi:\HH\to\Cx$ denote a conjugate linear functional on $\HH$. For
example, we could have $\pfi\in\HH$, and consider the functional
$f\mapsto \bra f,\pfi>$. We will use this notation even if
there is no $\pfi\in\HH$ allowing us to read $\bra f,\pfi>$ as a bona
fide inner product.
Then we define a \rep\
$\crep\pfi:\wick(T)\to\LL(\tA(\HH))$ by
$$\aligned
    \crep\pfi(f)X        &= f\otimes X  \\
    \crep\pfi(f\1)\idty  &= \bra f,\pfi>\idty
\endaligned\tag\deq(H.8)$$
for all $f\in\HH$. Thus $\crep\pfi(f)$ acts on the tensor algebra
$\tA(\HH)$ by left multiplication, and according to our notational
conventions we could also write $\crep\pfi(f)X= f\otimes X$.
The crucial point is that we do not need to define
$\crep\pfi(f\1) X$ for $X\neq\idty$. To see this, consider the
following rule, which is automatic from our axiomatic setup:
$$\aligned
 \crep\pfi(f\1)g_1\otimes \cdots\otimes g_n
     &=\crep\pfi(f\1)\crep\pfi(g_1)\
                    g_2\otimes \cdots\otimes g_n\\
     &=\bra g_1,f> g_2\otimes \cdots\otimes g_n\\
     &\qquad + \sum_{ijk\ell} T_{ij}^{k\ell}\
              \crep\pfi(\ell)\  \crep\pfi(k\1)\
                      g_2\otimes\cdots\otimes g_n
\quad. \endaligned\tag\deq(H.9)$$
On the right hand side $\crep\pfi(k\1)$ is applied to a tensor product
of only $n-1$ factors, so this is an inductive formula allowing us to
extend $\crep\pfi$ to tensors in $\tA(\HH)$ of arbitrary length $n$.
Note that we have not yet imposed an hermitian inner product on
$\tA(\HH)$. There is, however, only one way of introducing it:

\iproclaim/H.5/ Proposition.
Let $\wick(T)$ be a \HWa, and let $\pfi:\HH\to\Cx$ be a conjugate
linear functional. Then
\roster
\item[1]
There is a unique hermitian (but not necessarily positive
semidefinite) inner product $\iprod_{T,\pfi}$ on $\tA(\HH)$ with
$\bra\idty,\idty>_{T,\pfi}=1$, making $\crep\pfi$ a hermitian \rep.
\item[2]
There is a unique hermitian linear functional $\cst\pfi$ on
$\wick(T)$ such that $\cst\pfi(\idty)=1$, and
$$ \cst\pfi(f\,X)
     = \bra\pfi,f>\ \cst\pfi(X)
\qquad\hbox{for all $f\in\HH$, and $X\in\wick(T)$.} $$
\item[3]
$$\align
     \bra F,G>_{T,\pfi}&= \cst\pfi(F\1G)
\qquad\hbox{for all $F,G\in\tA(\HH)$, and} \\
          \cst\pfi(X)&=\bra\idty,\crep\pfi(X)\idty>_{T,\pfi}
\qquad\hbox{for all $X\in\wick(T)$.}
\\\endalign$$
\item[4]
The inner product $\iprod_{T,\pfi}$ is positive
semidefinite if and only if $\cst\pfi(X\1X)\geq0$ for all
$X\in\wick(T)$.
\endroster
\eproclaim

\proof:
We begin by showing (2).
Because $\cst\pfi$ is hermitian we must also have
$\cst\pfi(X\,g\1)= \bra g,\pfi>\,\cst\pfi(X)$.
The existence and uniqueness of $\cst\pfi$ follows from the
uniqueness of the Wick ordering process. The Wick ordered
monomials $i_1\cdots i_nj_1\1\cdots j_m\1$ form a basis of $\wick(T)$,
on which $\cst\pfi$ is explicitly defined as
$$ \cst\pfi\bigl(i_1\cdots i_nj_1\1\cdots j_m\1\bigr)
    = \prod_{k=1}^n   \bra\pfi,i_k> \quad
      \prod_{\ell=1}^m \bra j_\ell,\pfi>
\quad.$$
Suppose $\iprod$ is an inner product with the
properties specified in (1). Then
$$\align
   \bra \idty, \crep\pfi\bigl(i_1\cdots i_nj_1\1\cdots j_m\1\bigr)
              \idty>
    &=\bra \crep\pfi\bigl(i_n\1\cdots i_1\1\bigr)\idty,
           \crep\pfi\bigl(j_1\1\cdots j_m\1\bigr)\idty>     \\
    &=\prod_{k=1}^n   \bra\pfi,i_k> \quad
      \prod_{\ell=1}^m \bra j_\ell,\pfi>
\quad.\endalign$$
by inductive application of the formula
$\crep\pfi(f\1)\idty=\bra f,\pfi>\idty$. Since Wick ordered monomials
form a basis, this proves the second part of (3). Hence
$$\align
  \bra i_1\cdots i_n,\ j_1\cdots j_m>
     &=\bra \idty,
             \crep\pfi (i_n\1\cdots i_1\1\ j_1\cdots j_m)\idty>  \\
     &=\cst\pfi\bigl((i_1\cdots i_n)\1\ (j_1\cdots j_m)\bigr)
\quad,\endalign$$
which proves the first formula in (3), and hence (1), the existence and
uniqueness of $\iprod\equiv\iprod_{T,\pfi}$. (4) is immediate from
the first equation in (3).
\QED

It is interesting to note that in all these \rep s $\crep\pfi(i)$ is
the same operator on $\tA(\HH)$. Therefore, we may recover $T$ and
$\pfi$ from the inner product $\iprod_{T,\pfi}$.

We will call $\crep\pfi$ the {\em coherent \rep}, and $\cst\pfi$ the
{\em coherent state} associated with $\pfi$. The coherent state and
\rep\ with $\pfi=0$ are called the {\em Fock state} and the {\em Fock
\rep} of $\wick(T)$. We use our notation $\iprod_{T,\pfi}$
in such a way that zeros in the index are suppressed. Thus
$\iprod_{T}\equiv\iprod_{T,0}$ is the Fock inner
product, and $\iprod\equiv\iprod_{0,0}$ is the
Fock scalar product for the relations with $T=0$, i.e.\
$i\1j=\delta_{ij}$. A moment's reflection shows that this is the same
as the inner product of $\tA(\HH)$ inherited by the canonical inner
product on $\HH$. The completion of $\tA(\HH)$ in the associated norm is
the so-called ``full Fock space''.

In general the inner product $\iprod_{T,\pfi}$ may be
degenerate. For example, for the canonical commutation relations and
the Fock state $\cst0$, $f_1f_2-f_2f_1$ has zero norm. Dividing out
such vectors we obtain a new, usually simpler \rep. It coincides
with the cyclic (GNS-) \rep\ associated with the coherent state
$\cst\pfi$.

\iproclaim/H.6/ Definition.
Let $\pfi:\HH\to\Cx$ be a conjugate linear functional, and let
$\tA(\HH)/(T,\pfi)$ denote the quotient of $\tA(\HH)$ by the null
space of the inner product $\iprod_{T,\pfi}$. Then the \rep\
$\crep\pfi$ of $\wick(T)$ lifts to a \rep\ on $\tA(\HH)/(T,\pfi)$,
and if $\crep\pfi$ is bounded, it also extends to the completion of
$\tA(\HH)/(T,\pfi)$ with respect to $\iprod_{T,\pfi}$.
In either case this \rep\ will be called the {\bf separated coherent
\rep} associated with $\pfi$, and will be denoted by $\creps\pfi$.
\eproclaim

The following simple observation gives one reason why coherent \rep
s are useful.

\iproclaim/H.6a/ Proposition.
Suppose that the coherent \rep\ $\crep\pfi$  of $\wick(T)$ is
bounded. Then $\creps\pfi$ is irreducible.
\eproclaim

\proof:
When the coherent \rep\ is positive, the coherent functional
$\cst\pfi$ is uniquely characterized by
$$ \cst\pfi\bigl((f-\bra\pfi,f>)(f-\bra\pfi,f>)\1\bigr)=0
\quad.$$
This property would also hold for any component in a convex
decomposition of $\cst\pfi$ into positive functionals. Hence
$\cst\pfi$ is pure, and consequently its GNS-\rep\ $\creps\pfi$ is
irreducible by a standard result \cit{BraRo}.
\QED

\subhead \ssecno\sHw. A characterization of the Fock \rep   \endsubhead

The Fock \rep\ plays a distinguished r\^ole, as shown by the
following Theorem. It generalizes the Wold decomposition in
single operator theory \cit{Foias}, which
provides a canonical direct sum decomposition for an arbitrarily
given isometry as a sum of a unilateral shift (with multiplicity)
and a unitary. An isometry is, after all, the very simplest case of a
relation of the form \eq(H.1): we take a single generator, and all
coefficients $T_{ij}^{k\ell}=0$, so that $i\1i=\idty$. The Fock \rep\
of this relation is indeed the unilateral shift.
The classical Wold decomposition now allows us to decompose an
arbitrary \rep\ $\pi$ as $\pi=\pi_0\oplus\pi_1$, where $\pi_0$ is a
multiple of the Fock \rep, and the ``unitary'' sub\rep\ $\pi_1$ is
characterized by the property that every vector in the \rep\ space
can be given an arbitrarily long ``iteration history'' with respect
to the operators $\pi_1(i)$. In the following Theorem this statement
is extended verbatim to the general case of \Wa\ commutation
relations. In the following it is understood that an empty product
$\pi\bigl(i_1\cdots i_n\bigr)$ with $n=0$ is the identity.

\iproclaim/H.7/ Theorem (Wold decomposition).
Let $\pi:\wick(T)\to\BB(\R)$ be a bounded Hermitian \rep\ of a
\HWa\ on a Hilbert space $\R$.
Let
$$\align
  \No       &=\set{\pfi\in\R
              \stt  \forall_{i}\,
                         \pi\bigl(i\1\bigr)\pfi=0 }  \\
  \R_0      &=\clinsp\set{\pi\bigl(i_1\cdots i_n\bigr)\pfi
                     \stt n\geq0,\ i_1,\ldots,i_n\in I,
                          \pfi\in\No }  \\
  \R_1      &=\bigcap_{n\geq0} \clinsp\set
                     {\pi\bigl(i_1\cdots i_n\bigr)\psi
                     \stt i_1,\ldots,i_n\in I,\
                          \psi\in\R}
\quad.\\\endalign$$
Then $\R=\R_0\oplus\R_1$, and $\R_0$ and $\R_1$ are invariant under
$\pi(\wick(T))$ (either subspace may be $\set0$).
Moreover, the \rep\ $\pi$ restricted to $\R_0$ is a multiple of the
separated Fock \rep\  with multiplicity $\dim\No$.
\eproclaim

\proof:
Let $\xi,\eta\in\No$. Then
$$\bra\xi,\pi(X)\eta>=\bra\xi,\eta>\cst0(X)
$$
by definition of the Fock state.
Let $\eta_\alpha$, $\alpha=1,\ldots,\dim\No$ be an orthonormal basis
in $\No$. Then the cyclic subspaces
$\R_{0,\alpha}=\clinsp\pi(\wick(T))\eta_\alpha$ are orthogonal, and
restricted to each $\pi$ is the Fock \rep. The sum of these spaces is
$\R_0$.

Now let
$$\align
  \R_0\up n &=\clinsp\set{\pi\bigl(j_1\cdots j_m\bigr)\pfi
                     \stt m < n,\ j_1,\ldots,j_m\in I,
                          \pfi\in\No }  \\
  \R_1\up n &=\clinsp\set
                     {\pi\bigl(i_1\cdots i_n\bigr)\psi
                     \stt i_1,\ldots,i_n\in I,\ \psi\in\R}
\quad.\\\endalign$$
Then since for $m<n$ Wick ordering of
$ i_n\1\cdots i_1\1 j_1\cdots j_m$ leaves at least one factor $i\1$
to the right of every monomial, we have
$$ \bra \pi\bigl(j_1\cdots j_m\bigr)\pfi,\
        \pi\bigl(i_1\cdots i_n\bigr)\psi>=0
$$
for $\pfi\in\No$. Hence $\R_0\up n\subset\bigl(\R_1\up n\bigr)^\perp$
for all $n$.

We claim that even $\R_0\up n=\bigl(\R_1\up n\bigr)^\perp$. For this
we need to show that $\psi\in\bigl(\R_1\up n\bigr)^\perp$, or,
equivalently,
$$ \pi\bigl(i_n\1\cdots i_1\1\bigr)\psi=0
\eqno(*)$$
for all $i_1,\ldots,i_n\in I$, implies $\psi\in\R_0\up n$. Since the
decomposition
$\R=\R_0\oplus(\R_0)^\perp
   =\Bigl(\bigoplus_{\alpha=1}^{\dim\No} \R_{0,\alpha}\Bigr)
      \, \oplus(\R_0)^\perp$
is invariant under $\pi$, it suffices to show this for $\psi$ in each
summand separately. Assume first that $\psi\in\R_0^\perp$ satisfies
$(*)$. Then for all $i_1,\ldots,i_{n-1}$ we have
$\pi\bigl(i_{n-1}\1\cdots i_1\1\bigr)\psi\in\No\cap\R_0^\perp=\set0$.
Hence $(*)$ also holds for $n-1$. Proceeding by downwards induction we
find $\psi=0$. Assume next that $\psi\in\R_{0,\alpha}$ for one of the
Fock summands, generated from a vector $\eta_\alpha\in\No$.
Then for $m\geq n$ we have
$\pi\bigl(i_1\cdots i_n\bigr)\eta_\alpha\in\R_1\up n$.
By construction of the Fock \rep\ these vectors span the orthogonal
complement of $\R_0\up n$ so  we have $\psi\in\R_0\up n$.

Since $\R_0=\clinsp\bigcup_{n\in\Nl}\R_0\up n$, and
$\R_1=\bigcap_{n\in\Nl}\R_1\up n$, the complementarity of the
spaces $\R_i\up n$ for each $n$ implies $\R_0=\R_1^\perp$.
\QED

\head    \secno\sB.  Boundedness and Positivity           \endhead
\subhead \ssecno\sBt. The coefficients $T_{ij}^{k\ell}$ as an operator
\endsubhead

One might expect that positivity and boundedness of the Fock \rep\ of
a \Wa\ depend on positivity and boundedness of $\Tt$, considered
as an operator between tensor product Hilbert spaces. However, this is
{\em not} the case. To see this, let us compute the Fock inner product
of the ``two-particle'' space $\HH\otimes\HH$, i.e.\
$$ \bra ij,k\ell>_T
      :=\cst0\bigl(j\1i\1k\ell\bigr)
       = \delta_{ik}\delta_{j\ell}
            + T_{ik}^{\ell j}
\quad.\tag\deq(Bt.1)$$
We would like to consider this as a matrix element of an operator in
the $2$-fold Hilbert space tensor product $\HH\otimes\HH$ with its
natural inner product $\iprod$.
Thus we introduce an operator $T:\HH\otimes\HH\to\HH\otimes\HH$ with
the matrix elements
$$ \bra ij\abs T k\ell>=T_{ik}^{\ell j}
\quad,\tag\deq(Bt.2)$$
or, equivalently,
$$ T\,k\otimes\ell= \sum_{ij}T_{ik}^{\ell j}\ i\otimes j
\quad.\tag\deq\lasteq'(Bt.2')$$
Then we can write
$$\bra ij,k\ell>_T
    =\bra ij\abs{\idty+ T} k\ell>
    =\bra ij,P_2 k\ell>
\quad.\tag\deq(Bt.3)$$
Hence the Fock inner product of order $2$ is positive, if and only
if $P_2=\idty+T\geq0$. Notice in equation \eq(Bt.3) that the indices
$k\ell$ are associated with the domain of $T$, and $ij$ with its
range, whereas the pair $\ell j$ corresponds to the domain of $\Tt$.
Thus in contrast to the transposition operation, which would just
swap domain and range indices, only one index has changed sides,
i.e.\ $T$ is obtained from $\Tt$ by transposing only in one tensor
factor. Since the transposition operation is not completely bounded
or completely positive, the norm and positivity properties of $T$
may be quite different from those of $\Tt$. As the above computation
\eq(Bt.2) shows, it is $T$ rather than $\Tt$, which will be relevant for
questions like Fock positivity.

The formula $P_2=\idty+T$ has an extension to higher tensor powers. We
define the operators $P_n$ on the $n$-fold tensor product space
$\tA_n(\HH)=\HH^{\otimes n}$ such that
$$ \bra i_1\cdots i_n,\ j_1\cdots j_n>_T
       =:\bra i_1\cdots i_n,\ \ P_n\ j_1\cdots j_n>
\quad.\tag\deq(Bt.4)$$
Note that with respect to $\iprod_T$ the subspaces
$\tA_n(\HH)$ and $\tA_m(\HH)$, for $n\neq m$, are orthogonal as for
$\iprod$. Positivity of the Fock inner product means that
$P_n\geq0$ for all $n$. In order to state the formula for $P_n$ we
introduce the notation
$$ X_i\bigl(f_1\otimes f_2\cdots f_i\otimes
             f_{i+1}\otimes\cdots f_n\bigr)
    =\bigl(f_1\otimes f_2\cdots \otimes
    X(f_i\otimes f_{i+1})\otimes\cdots \otimes f_n\bigr)
\tag\deq(Bt.5)$$
where $i=1,\ldots,n-1$, and $X$ is any operator on a product Hilbert space
$\HH\otimes\HH$. Thus $X_i$ becomes an operator on the $n$-fold tensor
product $\tA_n(\HH)=\HH\otimes\HH\cdots\otimes\HH$.

\iproclaim/Bt.1/ Lemma.
For any  bounded operator $T$ on $\HH\otimes\HH$,
define a sequence $P_n(T)$ of bounded
operators on $\tA_n(\HH)$ by the inductive formula $P_1=\idty$,
$$ P_{n+1}(T)=(\idty\otimes P_n(T))
            \bigl(\idty+T_1+T_1T_2+\cdots+T_1\cdots T_n)
\quad.\tag\deq(Bt.6)$$
Then if $T$ is given by equation \eq(Bt.2), $P_n$ from equation
\eq(Bt.4) is equal to $P_n(T)$.
\eproclaim

\proof:
$P_n$ describes the difference between the Fock inner products
$\bra\cdot,\cdot>_T$ and $\bra\cdot,\cdot>_0\equiv\bra\cdot,\cdot>$.
Denoting by $\crep0$ and $\mu$ the Fock \rep s of the respective
relations on $\tA(\HH)$, we have from \eq(H.8):
$\crep0(i)X=\mu(i)X=i\otimes X$, and $\mu(i\1)\mu(j)=\delta_{ij}$.
$\crep0(i\1)$ is given by equation \eq(H.9). Our first task will be
to express it in terms of $\mu(i\1)$, in the form
$$ \crep0(i\1)X =\mu(i\1)R_nX
\quad\hbox{for $X\in\tA_n(\HH)$}
\quad,\eqno(*)$$
where $R_n$ is an operator on $\tA_n(\HH)$ to be determined in terms
of $T$ by an inductive formula. With $X\in\tA_n(\HH)$ we have
$$\align
   \mu(i\1)R_{n+1}\mu(j)X&= \crep0(i\1)\crep0(j)X \\
        &= \delta_{ij}X  + \sum_{k\ell} T_{ij}^{k\ell}\
                                 \crep0(\ell)\  \crep0(k\1)\ X \\
        &= \delta_{ij}X + \sum_{k\ell} T_{ij}^{k\ell}\
                                 \mu(\ell)\  \mu(k\1)R_n\ X \\
        &= \mu(i\1)\mu(j)X + \sum_{k\ell,nm} T_{nm}^{k\ell}\
           \mu(i\1)\mu(n) \mu(\ell)\mu(k\1) \mu(m\1)\ \mu(j)R_n\ X   \\
        &= \mu(i\1)\Big( \idty+T(\idty\otimes R_n)
                   \Big)\mu(j)X
\quad,\endalign$$
where at the last step we used the identity
$\mu(j)R_n=(\idty_\HH\otimes R_n)\mu(j)$ for operators on $\tA_n(\HH)$,
and introduced the operator $T$ from \eq(H.2) in the form
$$  T=\sum_{k\ell,nm} T_{nm}^{k\ell}\
           \mu(n) \mu(\ell)\mu(k\1) \mu(m\1)
\quad.$$
Hence we have proved the formula $(*)$ with
$$ R_{n+1}= \idty+ T\, (\idty_\HH\otimes R_n)
\quad,$$
or $R_n=\bigl(\idty+T_1+T_1T_2+\cdots+T_1\cdots T_{n-1})$.

The claim of the Lemma thus reduces to the formula
$P_{n+1}=(\idty_\HH\otimes P_n)R_{n+1}$. But
$$\align
    \bra \mu(i)X,P_{n+1} \mu(j)Y>
       &= \bra X,\crep0(i\1)\crep0(j)Y>_T
        = \bra X,\mu(i\1)R_{n+1}\mu(j)Y>_T                  \\
       &= \bra X,P_n\mu(i\1)R_{n+1}\mu(j)Y>
        = \bra X,\mu(i\1)(\idty\otimes P_n)R_{n+1}\mu(j)Y>  \\
       &= \bra \mu(i)X,(\idty\otimes P_n)R_{n+1}\mu(j)Y>
\quad,\endalign$$
which completes the proof.
\QED

\subhead \ssecno\sBe. The operator $T$ in the examples  \endsubhead

It is instructive to compute the operator $T$ in the examples
of Section \sE. In each example, the properties of $T$ as an operator
on a Hilbert space illuminates special features of the example.

For {\bf \Exa/E.1/}, the $q$-canonical commutation relations, we get
$$   T\vert ij\rangle =q\,\vert ji\rangle
\quad,\tag\deq(Bt.7)$$
i.e.\ $q$ times the flip operator.
Thus in the Bose case ($q=1$), $P_n(T)$, from equation \eq(Bt.4) is
proportional to the symmetrization projection. In the Fermi case
($q=-1$) each term gets a factor equal to the sign of the permutation
so that $P_n$ is proportional to the anti-symmetrization projection.

For {\bf \Exa/E.2/}, the Temperley-Lieb-Wick relations, we get
$$ T\vert ij\rangle = q \delta_{ij} \sum_k \vert kk\rangle
\quad.\tag\deq(Bt.8)$$
In order for this sum to converge, we must assume that the number $d$
of generators is finite. Then $T$ is $qd$ times a one-dimensional
projection, so $\norm{T}=\abs qd$. This readily accounts for the
$d$-dependence of the \rep\ theory, which we remarked in \cit{QCR},
without being able to pinpoint the reason for this dependence in
comparison with \Exa/E.1/.
One easily checks that $T$ satisfies the Temperley-Lieb relations
\cit{Temperley,Jones}
$$\aligned
    T_1T_2T_1&=q^2\,T_1  \\
         T^2 &=qd\ T
\quad.\endaligned\tag\deq(Bt.9)$$
As a generator of the Temperley-Lieb algebra it would be natural to
consider $T'=(qd)^{-1}T$, which is a projection, and satisfies
$T_1'T_2'T_1'=d^{-2}\,T_1'$.
This relation no longer contains the parameter $q$, from which we
conclude that the analysis of the associated Temperley-Lieb algebra
is not sufficient to decide the positivity questions we address in
this paper. Note that while in \Exa/E.1/ $T$ always has both
positive and negative eigenvalues, $T$ in this example is positive
for $q\geq0$. Hence by \Thm/Bp.1/, the Fock \rep\ of these relations
is positive for $q\geq0$. \Thm/Bs.2/ enlarges this range by
$\norm{T}=\abs qd<1/2$, so all together we have a positive Fock
\rep\ for $q>-1/(2d)$.

For {\bf \Exa/E.3/}, the twisted CCR and CAR,  we have
$$ T\ket{ij}=\cases
              \mu\ket{ji}      &i<j\\
              \mu^2\ket{ii}    &i=j\\
              -(1-\mu^2)\ket{ij}+\mu\ket{ji}
                               &i>j\\
\endcases\tag\deq(Bt.10)$$
in the Bosonic case \eq(E.3), and
$$ T\ket{ij}=\cases
              -\mu\ket{ji}      &i<j\\
              -\ket{ii}         &i=j\\
              -(1-\mu^2)\ket{ij}-\mu\ket{ji}
                               &i>j\\
\endcases\tag\deq(Bt.11)$$
in the Fermionic case \eq(E.4). In both these cases $T$ has the
eigenvalues $\mu^2$ and $-1$, and $\norm{T}=1$.

For {\bf \Exa/E.4/}, i.e.\ $\SnU2$ with the relations \eq(E.5), we
get the eigenvalues $0$ and $-(1+\nu^2)$, so $\norm{T}=1+\nu^2$.
As noted above, and as is evident from \eq(Bt.3) and the
eigenvalue $<-1$ of $T$, Fock positivity fails in this example.

For {\bf \Exa/E.5/}, the braided case, we get the rather more
transparent form $T_1T_2T_1=T_2T_1T_2$ of equation \eq(E.6), which
justifies calling it a braid relation. For the special case of
equation \eq(E.7) we find
$$   T\vert ij\rangle =q_{ij}\vert ji\rangle
\quad.\tag\deq(Bt.12)$$
Hence, in this case $\norm{T}=\max_{ij}\abs{q_{ij}}$.

{\bf \Exa/E.6/}, the case of Clifford algebras, is a special case of
\Exa/E.1/, so $T\vert ij\rangle =-\vert ji\rangle$.

The ``degenerate'' {\bf \Exa/E.6A/} is simply $T=-\idty$.
This is equivalent to $P_2(T)=0$, and $P_n(T)=0$ for all $n\geq2$ by
\Lem/Bt.1/. Note also that the braid relation is satisfied.

In {\bf \Exa/E.7/}, the case of infinitely many generators, we see
that $T$ may turn out to be unbounded: in \Exa/E.2/ we have seen
that $\norm{T}=\abs qd$ for the relations on $d$ generators, hence
for infinitely many generators $T$ is not bounded. It cannot even be
densely defined as an unbounded operator, as is evident from
\eq(Bt.8). On the other hand, $\nTt=1$ in this example. At the
other extreme we have \Exa/E.3/, where $\nTt$ diverges, but
$T$ is a well-defined bounded operator. It is not clear whether
boundedness of $T$, without the finiteness condition of \Def/H.1/,
would allow us to give a definition of the abstract Wick
algebra $\wick(T)$ on infinitely many generators. However, it
suffices to define the Fock \rep, assuming that it happens to be
positive: from \Lem/Bt.1/ we get the operators $P_n(T)$, which are
non-negative definite by assumption.  These define the inner product
on $\tA(\HH)$, and its completion, the Fock space. By construction,
this space has a natural orthogonal decomposition into
``$n$-particle'' spaces $\HH^n$. As shown in the proof of \Prp\sBs.1/Bs.1/
below, the operator $X\mapsto f\otimes X$ is bounded with
respect to the Fock scalar products as an operator taking $\HH^n$ to
$\HH^{n+1}$ (The bound is $\sum_{i=1}^n\norm{T}^n$). Hence we get a
well-defined closed operator $\creps0(f)$ as the orthogonal sum of
these bounded operators. The operators $\creps0(f)$ and
$\creps0(g)^*$ have the space of vectors with finite particle number
as a common dense and invariant domain. Clearly, the algebra
generated by these operators coincides with
$\creps0\bigl(\wick(T)\bigr)$, whenever the latter can be defined
with the methods of Section \sH.

In {\bf \Exa/E.8/}, i.e.\ the case of small coefficients
$T_{ij}^{k\ell}$, the operator $T$ will also be small. The next
section is devoted to showing that $\norm{T}$ is, in fact, the
appropriate measure of ``smallness'' for many questions.

\ppt{\vfill\eject}
\subhead \ssecno\sBs. Bounds for small $T$\endsubhead

In this section we discuss conditions (in the general case) under
which one can obtain Fock positivity or boundedness of the
Fock \rep. A crucial parameter in these problems is $\norm{T}$,
where $T$ is given in terms of the coefficients $T_{ij}^{k\ell}$ by
\eq(Bt.2).

The following Proposition generalizes Lemma 4 in \cit{BoSpei}.

\iproclaim/Bs.1/ Proposition.
Let $\norm{T}<1$, and suppose that  $P_n\geq0$ for all $n$.
Then in the Fock \rep\  the operators $\crep0(i)$ are
bounded with
$$ \norm{\crep0(i)}^2\leq {1\over 1-\norm{T}}
\quad.$$
\eproclaim

\proof:
For bounded operators $A,B,C$ with $A,B\geq0$ and $A=BC$, we have
$A^2=AA^*=BCC^*B^*\leq\norm{C}^2B^2$, and, by the operator
monotonicity of the square root \cit{RSimon,Donohu}, we have
$A\leq\norm{C}B$.
Applying this to the inductive formula \eq(Bt.6)
with $A=P_{n+1}$, $B=\idty\otimes P_{n}$, and using the norm
estimate
$\norm{C}=\norm{\idty+T_1+T_1T_2+\cdots+T_1\cdots T_{n-1}}
         \leq
$\ppt{\break}$
          (1-\norm{T}^n)/(1-\norm{T})\leq(1-\norm{T})^{-1}$,
we get $P_{n+1}\leq(1-\norm{T})^{-1}\idty\otimes P_{n}$.
Hence, for any vector $\Phi\in\HH^{\otimes n}$ we have
$$\align
    \bra\Phi, \crep0(i)^*\crep0(i) \Phi>_T
      &=    \bra i\otimes\Phi,\ P_{n+1}\  (i\otimes\Phi)> \\
      &\leq (1-\norm{T})^{-1}
            \bra i\otimes\Phi,\ (\idty\otimes P_n)\  i\otimes\Phi> \\
      &= (1-\norm{T})^{-1} \bra\Phi, \Phi>_T
\quad.\\\endalign$$
\QED

We note that the bound on $\norm{T}$ is best possible: for the
Bose canonical commutation relations we have $\norm{T}=1$, and it is
known (see e.g.\ \cit{BraRo}) that these have no bounded \rep s. On
the other hand, the criterion $\norm{T}<1$ is not a necessary
condition for the boundedness of the $a_i$, as the example of the
canonical anticommutation relations shows. In the next Theorem we
use the inductive formula to get the positivity of $P_n$ for
sufficiently small $T$. Since the factors in the inductive formula
do not commute, we cannot use it directly to show positivity of
$P_n$. Instead, we show first that the $P_n$ are non-singular, and
then invoke analytic perturbation theory to get positivity. With a
different argument for the non-singularity of $P_n$ this idea has
been used by \cit{Zagier,Fivel} to treat the $q$-canonical
commutation relations, and by \cit{BoSpein} in a much more general
context.

\iproclaim/Bs.2/ Theorem.
For $\norm{T}<1/2$, $P_n\geq\eps_n\idty>0$  for all $n$.
\ppt{\hfill\break}
For $\norm{T}\leq1/2$,  $P_n\geq0$ for all $n$.
\eproclaim

\proof:
Assume  $\norm{T}<1/2$. Then in the inductive formula
$$ P_{n}=P_{n}(T)
        =(\idty\otimes P_{n-1}(T))
            \bigl(\idty+\sum_{r=1}^{n-1}\prod_{i=1}^r T_i\bigr)
$$
the sum in the second factor is bounded by
$\sum_{r=1}^{n-1}\norm{T}^r\leq\norm{T}/(1-\norm{T})<1$.
Hence
$$ \norm{P_n(T)^{-1}}
      \leq {1-\norm{T}\over1-\norm{2T}}\,
                \norm{P_{n-1}(T)^{-1}}
\quad,\tag\deq(Bs.1)$$
and $P_n(T)$ is invertible for all $n$.

Note that $P_n(\lambda T)$ is hermitian and invertible for all
$\lambda\in\bracks{0,1}$. The eigenprojection of $P_n(\lambda T)$ for
the negative half axis can thus be computed by the same Cauchy
integral in the analytic functional calculus for all $\lambda$, and is
hence a norm continuous function of $\lambda$. Since this projection
vanishes for $\lambda=0$, it must be identically zero.
Hence $P_n>0$ for all $n$.

Since $P_n(T)$ is a polynomial in translates of $T$, it is clear that
$\bra\pfi,P_n(\lambda T)\pfi>$ is a continuous function of $\lambda$.
If $\norm{T}=1/2$ this function is positive for $0<\lambda<1$, hence
also for $\lambda=1$.
\QED

We now show that, without further information about $T$, the above
bound on $\norm{T}$ is best possible.

\iproclaim/Bs.ce/ Example.
Consider the \Wa\ $\wick(T)$ on two generators $x,y$ with the
relations
$$\aligned
    x\1y&=y\1x=0 \\
    x\1x&=\idty+\tau(xx\1-yy\1)  \\
    y\1y&=\idty-\tau(xx\1-yy\1)
\quad,\\\endaligned\tag\deq(Bs.2)$$
where $\tau\geq0$ is a real parameter. Then $\norm{T}=\tau$.
For $\tau>1/2$, $\wick(T)$ has no bounded \rep s, and for
$1/2<\tau<1$ the Fock \rep\ is not positive.
\eproclaim

\proof:
The defining coefficients are
$$ T_{ij}^{k\ell} =\tau(-1)^{j+k}\delta_{ij}\delta_{k\ell}
\quad,\quad i,j,k,\ell=1,2\quad.$$
Thus $T$ is diagonal in the tensor product basis, and has eigenvalues
$\pm\tau$. In particular, $\norm{T}=\tau$. Consider the completely
positive map
$$ \Phi(a)={1\over2}\bigl(x\1ax+y\1ay\bigr)
\quad.$$
By the relations it satisfies $\Phi(\idty)=\idty$, and for the
special element $\Delta=(xx\1-yy\1)$ a short computation gives
$\Phi(\Delta)=2\tau\Delta$.

We show next that, for $\norm{T}=\tau>1/2$, the relations \eq(Bs.2)
have no bounded \rep s. If $\pi$ were a bounded \rep, $\pi(\Delta)$
would be bounded, and
$2\tau\norm{\pi(\Delta)}
      =\norm{\pi\Phi(\Delta)}
      \leq\norm{\pi(\Delta)}$.
When $\tau>1/2$, this implies $\pi(\Delta)=0$. Thus
$x:=\pi(x\1x)=\pi(yy\1)$, and $x^2=\pi(x\,x\1y\, y\1)=0$. It follows
that $x=0$, hence $\pi(x)=\pi(y)=0$ in contradiction with the Wick
relations.

Now let $1/2<\tau<1$. Then if the Fock \rep\ were positive, we could
conclude by \Prp/Bs.1/ that $\creps0(x),\creps0(y)$ are
bounded, which is impossible by the preceding paragraph.
\QED

\subhead \ssecno\sBu. The universal bounded \rep \endsubhead

In this section we prove the criterion for the isomorphism between
the universal bounded \rep\ and the Cuntz-Toeplitz algebra for
sufficiently small $T$.
First we need a bound on the right hand side of the relations \eq(H.1).
We will think of it as an operator in block matrix form with respect
to $i,j$. Thus we set
$$ M_{ij}(a,b)=\sum_{k\ell} T_{ij}^{k\ell}\ a^*_\ell b_k
\quad,\tag\deq(Bs.3)$$
where $a_i,b_i$ are arbitrary collections of bounded operators
(usually $b_i=a_i=\creps0(i\1)$). The right hand side can also be
written as a matrix
$$ E_{ij}(a,b)=b_ia_j^*
\quad.\tag\deq(Bs.4)$$
In this notation the relations \eq(H.1) become $E(a,a)=\idty+M(a,a)$,
where $\idty$ now refers to the identity in the $d$-fold direct sum of
the space on which the $a_i$ act. For any $d$-tuple of bounded
operators $a_i$ we define the norm
$$ \norm{a}^2=\Norm\Big{\sum_ka_k^*a_k}
\quad.\tag\deq(Bs.5)$$

\iproclaim/Bs.3/ Proposition.
Suppose that $(a_i)_{i\in I}$ and $(b_i)_{i\in I}$ are tuples of
operators with
\ppt{\break}
$\norm{a},\norm{b}<\infty$. Let $t_+,t_-$ denote the
supremum and infimum of the spectrum of $T$, respectively.
Then
\roster
\item[1]
$\norm{E(a,b)}\leq\norm{a}\norm{b}$,
\item[2]
$\norm{E(a,a)}=\norm{a}^2$,
\item[3]
$ \norm{M(a,b)}\leq\norm{T}\norm{a}\norm{b}$, \quad and
\item[4]
$ t_-\norm{a}^2\idty \leq  M(a,a)\leq t_+\norm{a}^2\idty$.
\endroster
\eproclaim

\proof:
(1) and (2) were proven in \cit{QCR}, Lemma 10.

(3) Let $\R$ denote the Hilbert space on which the $a_i, b_i$ act.
Let $\xi,\eta\in\R^I$. Then
$$\align
  \bra\xi,M(a,b)\eta>
    &=\sum_{ij}\bra\xi_i,M_{ij}(a,b)\eta_j>
     =\sum_{ijk\ell} \bra i\ell \abs T jk>
                   \bra i\ell\abs R jk>         \\
    &=\tr \bigl(T R \bigr)
\quad,\endalign$$
where $R$ denotes the operator with matrix elements
$$ \bra i\ell \abs R jk>
     =\bra a_\ell\xi_i, b_k\eta_j>
\quad.$$
Now by \Lem\sBu.2/Bs.4/ below, $R$ is trace class, with trace norm
$$
   \norm{R}_1^2
      \leq \bigl(\sum_{i\ell}\norm{a_\ell\xi_i}^2\bigr)
            \bigl(\sum_{kj}   \norm{b_k  \eta_j}^2\bigr)
      \leq \norm{\xi}^2\norm{\eta}^2\ \norm{a}^2 \norm{b}^2
\quad.$$
Hence
$$ \abs{\bra\xi,M(a,b)\eta>}
    \leq \norm{T}\norm{a} \norm{b}\ \norm{\xi}\norm{\eta}
\quad.$$

(4) When $a=b$ and $\xi=\eta$, $R$ is a positive trace class operator,
hence $\bra\xi,M(a,a)\xi>$ can be estimated above by the corresponding
expression with $T=t_+\idty$, which is
$t_+\sum_{i\ell}\norm{a_\ell\xi_i}^2\leq t_+\norm{a}^2\norm{\xi}^2$.
The lower bound follows analogously.
\QED

\iproclaim/Bs.4/ Lemma.
Let $d\leq\infty$, and let $\HH$ be a Hilbert space. For $i\in\Nl,
1\leq i\leq d$, let $\psi_i,\pfi_i\in\HH$. Suppose that
$\norm{\psi}^2=\sum_i\norm{\psi_i}^2<\infty$, and
$\norm{\pfi}^2=\sum_i\norm{\pfi_i}^2<\infty$. Then there is a unique
trace-class operator $R$ on $\Cx^d$ with
$\bra i\abs R j>=\bra\psi_i,\pfi_j>$, which satisfies the estimate
$\tr\abs R  \leq \norm{\psi}\norm{\pfi}$.
\eproclaim

\proof:
Consider the operator $V_\pfi:\Cx^d\to\HH$ given by
$V_\pfi\vert i\rangle=\pfi_i$.
Then $V_\pfi V_\pfi^*=\sum_i\abs{\pfi_i\rangle\langle\pfi_i}$. This
positive operator has trace $\norm{\pfi}^2$, hence $V_\pfi$ is
Hilbert-Schmidt class with Hilbert-Schmidt norm
$\norm{V_\pfi}_2^2=\norm{\pfi}^2$. With $R=V_\pfi^*V_\psi$, and
a standard trace norm estimate we get
$\tr\abs R=\norm{R}_1\leq \norm{V_\pfi}_2\norm{V_\psi}_2
          \leq\norm{\pfi}\norm{\psi}$.
\QED

The bounds in \Prp/H.3/ imply the existence of universal
\rep s in the following sense:

\iproclaim/Bs.5/ Proposition.
Let $T\leq t_+\idty$ with $t_+<1$, and suppose that there is a
\cbo\ \rep\ of $\wick(T)$ in the sense of
\Def/H.1/. Then there is a universal \cbo\ \rep\ with
$$ \beta=\Norm\Big{\sum_i\pi(ii\1)}\leq(1-t_+)^{-1}
\quad.$$
\eproclaim

\proof:
The tuple $a_i=\pi(i\1)$, $i=1,\ldots\abs I$ satisfies the relations
$E(a,a)=\idty+M(a,a)$, and $\beta=\norm{a}^2$. From \Prp/Bs.3/(2,4)
we get $\norm{a}^2\leq1+t_+\norm{a}^2$.
\QED

The following is the application of Theorem 9 in \cit{QCR} in the
present context.

\iproclaim/Bs.6/ Theorem.
Let $d<\infty$, and let $t_\pm\in\Rl$ with $t_-\idty\leq T\leq
t_+\idty$. Suppose that
$$  \max\set{\abs{t_+},\abs{t_-}}^2
             < 1-t_++t_-
\quad,\tag\deq(Bs.6)$$
or, more specially, $\norm{T}\leq\sqrt2-1 \approx.414$, or
$\norm{T}\leq(\sqrt5-1)/2\approx.618$ with either $T\geq0$ or
$T\leq0$.
Then the universal bounded \rep\  of $\wick(T)$ exists, and
$\wicc(T)$ is isomorphic to the Cuntz-Toeplitz algebra $\wicc(0)$.
\eproclaim

\proof:
We verify the premises of Theorem 9 in \cit{QCR}. From \Prp/Bs.5/
it is clear that any solution to $E(a,a)=\idty+M(a,a)$ must
satisfy $\norm{a}^2\leq(1-t_+)^{-1}\equiv\mu$. Moreover, this bound
on an arbitrary tuple $a$ implies that $\idty+M(a,a)\leq\mu\idty$. Under
the same hypothesis
$\idty+M(a,a)\geq(1+t_-\mu)\idty\equiv\eps\idty$.
Finally from \Prp/H.3/ we get for
$\norm{a}^2,\norm{b}^2\leq\mu$ the Lipshitz bound
$\norm{M(a,a)-M(b,b)}\leq\norm{M(a-b,a)}+\norm{M(b,a-b)}
  \norm{T}(\norm{a}+\norm{b})\norm{a-b}\leq\lambda\norm{a-b}$
with $\lambda=2\sqrt\mu\norm{T}$. The hypothesis of Theorem 9 in
\cit{QCR} now requires $\lambda<2\sqrt\eps$, or
$\norm{T}^2< 1-t_+-t_-$.
Since $\norm{T}\leq\max\set{\abs{t_+},\abs{t_-}}$, this is
implied by the condition given in the Theorem. The result
then follows from the cited Theorem.
\QED

\subhead \ssecno\sBp. Positivity of the Fock \rep\ for positive $T$
\endsubhead

In this section we show a restricted form of monotonicity of the map
$T\mapsto P_n(T)$, namely that the positivity of the operator $T$ is
also a sufficient condition for the positivity of the Fock \rep. We
will comment on possible stronger versions of monotonicity at the
end of the section.

\iproclaim/Bp.1/ Theorem. Let $T$ be a bounded operator on
$\HH\otimes\HH$ with $T\geq0$, and define $P_n(T)$ as in \Lem/Bt.1/.
Then $P_{n+1}(T)\geq\idty\otimes P_n(T)\geq\idty$ for all $n$.
\eproclaim

\proof:
It suffices to consider the finite dimensional case $d<\infty$,
since $P_n(pTp)$ for $p$ a finite dimensional projection
converges strongly to $P_n(T)$ as $p\nearrow\idty$. So let
$$ \bra i\ell \abs T jk>
       =\sum_\alpha\Br{\Phi^\alpha_{i\ell}}{\Phi^\alpha_{jk}}
\quad.\tag\deq(Bp.1)$$
We now proceed by induction and assume that $P_n\geq0$. We have to
show that the second term on the right in
$$\align
  \bra i_0,i_1,\ldots i_n \mid j_0,j_1,\ldots j_n>_T
    =\delta_{i_0,j_0}& \bra i_1,\ldots i_n \mid j_1,\ldots j_n>_T \\
      &+\sum_{k\ell} \bra i_0\ell \abs T j_0k>\
           \bra i_1,\ldots i_n \abs{\crep0(\ell k\1)}
                j_1,\ldots j_n>_T
\endalign$$
is positive definite. So let $\Psi_{i_0,i_1,\ldots i_n}$ be
arbitrary complex coefficients. Then
$$\align
  \sum \Br{\Psi_{i_0,i_1,\ldots i_n}}\Psi_{j_0,j_1,\ldots j_n}
        &\bra i_0\ell \abs T j_0k>\
        \bra i_1,\ldots i_n \abs{\crep0(\ell k\1)}
             j_1,\ldots j_n>_T    \\
   &=\sum_\alpha\norm{\Psi_{j_0,j_1,\ldots j_n}
                     {\Phi^\alpha_{j_0k}}
                       \crep0(k\1)\ket{j_1,\ldots j_n}}_T^2
\quad.\endalign$$
This norm is evaluated in the $n$-particle Fock space, whose scalar
product is positive definite by inductive hypothesis.
Hence $P_{n+1}=\idty\otimes P_n +$positive terms. The bound
$P_n\geq\idty$ follows from this by induction, since $P_1=\idty$.
\QED

In view of \Thm/Bs.2/ one might be led to believe that the
above Theorem can be improved to hold for all
$T\geq-(1/2)\idty$. However, the following example, a slight
modification of \Exa/Bs.ce/ shows that the
lower bound $T\geq0$ in the Theorem is optimal:

\iproclaim/Bp.ce/ Example.
Consider the \Wa\ $\wick(T)$ on two generators $x,y$ with the
relations
$$\aligned
    x\1y&=y\1x=0 \\
    x\1x&=\idty+  \lambda xx\1 +\eps yy\1  \\
    y\1y&=\idty+ \eps xx\1 +\lambda  yy\1
\quad,\\\endaligned\tag\deq(Bp.2)$$
where $\lambda>\eps$.
Then $\eps$ is the best lower bound for $T$, and, for any
$\eps<0$, we can find $\lambda$ such that the Fock \rep\ of
$\wick(T)$ is not positive.
\eproclaim

\proof:
Here the coefficients are
$$ T\ket{ij}=\bigl( \lambda\delta_{ij}+\eps(1-\delta_{ij})\bigr)
             \ket{ij}
\quad.$$
Then $\idty+T\geq(1+\eps)\idty>0$ is boundedly invertible, and
$P_3=\idty+T_2+(\idty+T_2)T_1(\idty+T_2)$ is positive iff
$(\idty+T_2)^{-1}+T_1$ is positive. But
$$ \bra xyy\abs{(\idty+T_2)^{-1}+T_1}xyy>
    ={1\over1+\lambda}+\eps
\quad.$$
Hence for any $\eps<0$ it suffices to take
$\lambda>1+(-\eps)^{-1}$ to obtain a relation with
$T\geq\eps\idty$ and non-positive Fock space.
\QED

Note that this example also disproves the monotonicity conjecture
that $(T'\leq T)$ implies $(P_n(T')\leq P_n(T))$: for $T'=-1$, i.e.\
\Exa/E.6A/, we have $P_n(T')\geq0$ for all $n$, but in the example
above we have $T\geq\eps\idty>T'$, as long as $\eps\geq-1$, and
$P_n(T)$ is not positive.

\subhead \ssecno\sBq. $T$ satisfying the braid relations \endsubhead

We have seen that without further assumptions the condition
$\norm{T}\leq1/2$ is the optimal sufficient bound for the Fock \rep\
to be positive. In \Exa/E.1/ and \Exa/E.3/, however, it is the bound
$\norm{T}\leq1$ that marks the boundary of Fock positivity. A similar
observation was made by Speicher \cit{Speix}, who generalized the
$q$-relations to the
``$q_{ij}$''-relations
$$    i\1j=\delta_{ij}\idty + q_{ij}\ ji\1
\quad,\tag\deq(Bq.1)$$
where the $q_{ij}$ are real constants. Then
$$   T\vert ij\rangle =q_{ij}\vert ji\rangle
\quad,\tag\deq(Bq.2)$$
and this operator has a natural factorization into $T=Q\pi_1$, where
$$   Q\vert ij\rangle =q_{ji}\vert ij\rangle
\quad,\tag\deq(Bq.3)$$
and $\pi_1\vert ij\rangle=\vert ji\rangle$ is the flip operator.
Then $\norm{T}=\norm{Q}=\max_{ij}\abs{q_{ij}}$, and Speicher's
result \cit{Speix} is again the Fock positivity for $\norm{T}\leq1$.

One feature that makes these relations tractable is the commutativity
of $Q$ with all its permuted versions. In fact, the following
Proposition shows that this feature can be used to
characterize relations of the form \eq(Bq.1). We use the following
notation: when $\HH$ is a Hilbert space, and $Q$ is a bounded
operator on $\HH\otimes\HH$, we denote by $Q_{ij}$ the operator on the
$n$\th tensor power of $\HH$, acting on the tensor slots $i,j$ as $Q$
acts on slots $1,2$, or, more precisely, $Q_{ij}$ is the image of
$Q\otimes\idty^{\otimes(n-2)}$ under the permutation automorphism
taking $(1,2)$ to $(i,j)$.

\iproclaim/Bq.1/ Proposition.
Let $\HH$ be a finite dimensional Hilbert space, and let $Q$ be an
operator on $\HH\otimes\HH$.
Then the following are equivalent:
\roster
\item[1] $Q_{12}=Q_{21}^*$, and the operators $Q_{ij}$ and
$Q_{k\ell}$ commute for all $n$, and $i,j,k,\ell=1,\ldots,n$.
\item[2]
There are a basis $\ket i$ in $\HH$, and complex numbers $q_{ij}$ such
that $Q\ket{ij}=q_{ij}\ket{ij}$.
\endroster
\eproclaim

\proof:
Let $\cEx_1^\omega:\BB(\HH\otimes\HH)\to\BB(\HH)$ denote the conditional
expectation with respect to a functional $\omega$ on the second
factor, i.e.\
$\rho\bigl(\cEx_1^\omega(A\otimes B)\bigr)
  =\rho\otimes\omega(A\otimes B)$, and define
$\cEx_2^\omega$ similarly with respect to a functional on the first
factor. Then because $Q_{12}$ and $Q_{13}$ commute, we find that the
operators of the form $\cEx_1^\omega(Q)$ commute, for all choices of
$\omega$. Similarly, the operators $\cEx_2^\omega(Q)$ commute, and
because $Q_{12}=Q_{21}^*$ we find among them the adjoints of the
first commuting set. Finally, because $Q_{12}$ and $Q_{23}$ commute,
the $\cEx_i^\omega(Q)$, $i=1,2$, also commute with each other, and
hence generate an abelian *-subalgebra $\ZZ\subset\BB(\HH)$. Let
$P_\alpha$ denote the minimal projections of $\ZZ$. Then because
$\cEx_1^\omega(Q)\in\ZZ$ we must have
$\cEx_1^\omega(Q)=\sum_\alpha \omega(X_\alpha)P_\alpha$, and by the
symmetrical condition, $X_\alpha\in\ZZ$. Hence we have
$Q\in\ZZ\otimes\ZZ$, or
$Q=\sum_{\alpha\beta} \tilde q_{\alpha\beta}\
          P_\alpha\otimes P_\beta$.
The result now follows by choosing a basis in $\HH$ in which $\ZZ$ is
diagonal.
\QED

Note that the first condition in \pcl/Bq.1/(1) is just
the hermiticity of
$T$. Moreover, the Proposition shows that it is natural to consider
complex $q_{ij}$ rather than just real ones. The following
Proposition shows that in this case, too, $\norm{T}\leq1$ implies
Fock positivity. Since $T\geq-\idty$ is equivalent to the positivity
of $P_2$, this bound is clearly optimal.

\iproclaim/Bq.2/ Proposition.
For $i,j=1,\ldots,d$, let $q_{ij}\in\Cx$ with $\Br{q_{ij}}=q_{ji}$,
and $\abs{q_{ij}}\leq1$.
Then the \Wa\ defined by  $i\1j=\delta_{ij}\idty + q_{ij}\ ji\1$
has a positive Fock \rep.
\eproclaim

We proved this result in an earlier version of this paper. In the
meantime, however, \Bozejko\ and Speicher \cit{BoSpein} have shown a
strictly stronger result, which has a very natural statement in the
setting of \Wa s. We will therefore describe their result below, and
give only a sketch of our direct proof of the above Proposition. In
the course of this sketch we will also describe the structures
needed to state the result of \cit{BoSpein}.

Consider for the moment the factorization $T=Q\pi_1$ without further
assumptions on $Q$. Then
$$ T_iT_{i+1}T_i=T_{i+1}T_iT_{i+1}
   \iff Q_{ij}Q_{ik}Q_{jk}
      = Q_{jk}Q_{ik}Q_{ij}
\quad,\tag\deq(Bq.4)$$
where we have taken into account that by applying permutation
automorphisms to the equation on the right, the special case
$(i,j,k)=(i,i+1,i+2)$  is transformed into the case of general
$(i,j,k)$.
(In fact, \eq(Bq.4) is nothing but the equivalence of the braid
relation stated in terms of the operations $T_i$ of braiding the
current $i$\th\ and $(i+1)$\th\ strands, and the same relation stated
in terms of the operations $Q_{ij}$ braiding the $i$\th\ and $j$\th\
labelled strand).
Under the assumptions of \Prp/Bq.2/, \Lem/Bq.1/(1) clearly
implies the validity of the identity on the right, and hence the
$T_i$ satisfy the braid relation on the left.
This has an important consequence. Consider the inductive
expression for $P_n$ in \Lem/Bt.1/. Multiplying out the product for
$P_n$ we obtain $n!$ terms. Each one of these terms is labelled in a
unique way by a permutation $\pi$, which is obtained by replacing the
product $T(\pi)=T_{i_1}T_{i_2}\cdots T_{i_k}$ of $T$-operators by
$\pi=\pi_{i_1}\pi_{i_2}\cdots \pi_{i_k}$, where $\pi_i$ is the
transposition exchanging $i$ and $i+1$.
Moreover, $\pi_{i_1}\pi_{i_2}\cdots \pi_{i_k}$ is one
of the \rep s of $\pi$ as a minimal-length product of
transpositions. Now the braid
relation \eq(Bq.4) has the consequence that different minimal-length
factorizations of $\pi$ lead to the same operator $T(\pi)$. To see
this, represent each permutation by a planar diagram with points
labelled $1,\ldots,n$ on the bottom and top lines. Now connect each
point $i$ in the top row with $\pi(i)$ in the bottom line, such that
no connecting line turns upwards anywhere, and such that no three
lines intersect in the same point. This diagram can be read as a
\rep\  of $\pi$ as a product of nearest neighbour
transpositions. The product is of minimal length iff no two lines
intersect twice. Now if the braid relation holds, we can move one
strand in this diagram over the crossing of any other two, and thus
we can transform any two ways of connecting the top and bottom rows
according to $\pi$ into each other. Hence the braid relation implies
that $T(\pi)$ depends only on the permutation $\pi$, and is
``quasi-multiplicative'' \cit{BoSpein}, i.e.\
$T(\pi\sigma)=T(\pi)T(\sigma)$, as long as $\pi$ and $\sigma$ are
sub-products in a minimal length factorization of $\pi\sigma$.

The crucial property to establish to get
\Prp/Bq.2/ and its generalizations is the positive
definiteness of the kernel $(\pi,\sigma)\mapsto T(\pi^{-1}\sigma)$,
i.e.\ the operator inequality
$$  \sum_{\pi\sigma}\Br{a_\pi}a_\sigma T(\pi^{-1}\sigma)
       \geq0
\tag\deq()$$
for arbitrary choices of $a_\pi\in\Cx$. For if this holds, we get
$P_n(T)=\sum_\pi T(\pi)=$ \ppt{\break}
$=(1/N!)\sum_{\pi,\sigma}T(\pi^{-1}\sigma)\geq0$.

\demo{Sketch of proof of \pcl/Bq.2/:}
We want to use the commutativity of the $Q_{ij}$. Hence, for every
permutation $\pi$,  we introduce
$Q(\pi)=T(\pi)\pi^{-1}$.
Then, because
$$ \pi Q_{ij}=Q_{\pi i,\pi j}\ \pi
\quad,$$
we find that $Q(\pi)$ is a product of factors $Q_{ij}$. This is true
even when these factors do not commute. However, if they do, we get
the following simple expression for $Q(\pi)$:
$$ Q(\pi)=\Qprod{i<j}{\pi^{-1}i>\pi^{-1}j} Q_{ij}
\quad.$$
Consider the kernel
$$  \widehat Q(\pi,\sigma)
      :=\pi T(\pi^{-1}\sigma)\sigma^{-1}
       =\pi Q(\pi^{-1}\sigma)\pi^{-1}
       =\Qprod{\pi^{-1}i<\pi^{-1}j}
            {\sigma^{-1} i>\sigma^{-1} j} Q_{ij}
\quad.$$
Then the positivity of this kernel is equivalent to the positivity
of the kernel $T(\pi^{-1}\sigma)$.
Moreover, $\widehat Q$ is an abelian product of simpler kernels, each one
associated with a pair $(i,j), (j,i)$ of indices. These kernels are
readily seen to be positive definite, and hence so is $\widehat Q$.
\QED

It is clear from this outline that the braid property of $T$, as the
necessary condition for defining the operators $T(\pi)$, is a
crucial element of this proof. On the other hand, in \Exa/E.3/ the
braid property is also satisfied. This had suggested to us the
question whether the braid property alone suffices to show Fock
positivity up to $\norm{T}=1$. The affirmative answer was given by
\cit{BoSpein} in the following Theorem, which we cite here in the
terminology of this paper.

\iproclaim/Bq.3/ Theorem (\cit{BoSpein}).
Let $T$ be a bounded operator on $\HH\otimes\HH$, satisfying the braid
relation $T_1T_2T_1=T_2T_1T_2$, and
define the operators $P_n(T)$ as in \Lem/Bt.1/.
Then $\norm{T}\leq1$, implies $P_n(T)\geq0$ for all $n$,
and  $\norm{T}<1$     implies $P_n(T)>0$ for all $n$.
\eproclaim

We have already argued that the first bound is optimal. To see the
optimality of the second, consider the CAR-relations.

\ppt{\vfill\eject}
\head    \secno\sI.  Ideals in Wick algebras            \endhead

This section deals with an important special kind of ideals in a
\Wa, which we will call ``Wick ideals''. Loosely speaking,
they are the ideals which are generated by subsets of $\tA(\HH)$,
i.e.\ by functions of the generators $i\in I$ which do not depend on
the adjoints $i\1$. Their study is motivated by the well-known case
of the canonical (anti-) commutation relations, which are usually
given as {\em two} sets of relations: one set of Wick ordering type,
for commuting $a_i$ with $a_j^*$, and another set for commuting
$a_i$ with $a_j$. Wick ideals describe those additional relations
that can be introduced consistently in a given \Wa.

The consistency problem between Wick relations and the relations
defining a Wick ideal has often been posed in the literature
from the opposite perspective \cit{WessZ,PWor,Manin}. These authors
usually start from an algebra (without involution) defined by
generators and relations. They then ask for the possible choices of
coefficients $T_{ij}^{k\ell}$ such that the corresponding rules
\eq(H.1) are consistent with the already given relations in the
algebra.
In order to put this more formally, recall that in $\wick(T)$
there are no relations between the generators $i$, so that the
canonical embeddings $\tA(\HH)\hookrightarrow\wick(T)$ and
$\tA(\HH\1)\hookrightarrow\wick(T)$ are injective.
However, this injectivity is not essential to the idea of Wick
ordering. A more general approach would be to consider a *-algebra
$\wick$, which is generated by a subalgebra $\A\subset\wick$, and its
adjoints $\A\1=\set{x\1\stt x\in\A}$ in such a way that factors from
$\A$ and $\A\1$ may be disentangled completely. In other words, we
demand that $\wick$ is spanned linearly by $\A\A\1$. If $\A$ is
generated by $i\in I$, we have a homomorphism
$\tA(\HH)\to\A\hookrightarrow\wick$, which is now no longer injective.
Its kernel is the ideal describing the relations between generators.
It is still clear that there must be some rule reordering $i\1j$ into
products $k\1\ell$, but this commutation rule may not be chosen
independently of the ideal. In this picture ``Wick ideals'' are the
ideals in $\tA(\HH)$ which are consistent with the commutation rules
given by $T$.

The importance of Wick ideals for the questions of positivity and
boundedness addressed in this paper is twofold: on the one hand it
is easy to see that certain ideals of this kind are annihilated by
certain coherent \rep s (see \Lem\sIq.2/Iq.co/ below). For example,
every quadratic ideal (in the sense described below) automatically
vanishes in the Fock \rep. The additional relations allow one to
simplify the \rep\ space, and so to decide positivity questions more
easily. An extreme example are the twisted commutation relations of
\Exa/E.3/, where the verification of the positivity of the Fock
\rep\ reduces to a triviality \cit{PWor,Pusz}. On the other hand, we
may consider a Wick ideal, and study its \rep\ theory as an algebra
in its own right. We will find that the generators of a quadratic
Wick ideal satisfy ``Wick relations without constant term''. These
have typically no bounded \rep s, so that we can conclude that the
ideal has to vanish in every bounded \rep\ of the given \Wa. We will
treat \Exa/E.3/ from this point of view in Section \sIt.

\ppt{\vfill\eject}
\subhead \ssecno\sIq. Wick ideals and quadratic Wick ideals
\endsubhead

In the following Lemma we use the following standard notation: when
$X,Y$ are linear subspaces of an algebra, $XY$ denotes the linear span
of the products $xy$ with $x\in X,\ y\in Y$. With this notation, the
characteristic property of \Wa s is
$$\aligned
     \HH\1\HH          &\subset  \Cx\idty + \HH\HH\1     \\
     \tA(\HH\1)\tA(\HH)&\subset   \tA(\HH)\tA(\HH\1)
                      \quad=\quad\wick(T)
\quad.\endaligned\tag\deq(Iq.1)$$

\iproclaim/Iq.1/ Lemma.
Let $\wick(T)$ be a \Wa\  $\ideal\subset\tA(\HH)$ be an
ideal in the tensor algebra over $\HH$, and let
$\ideal_0\subset\ideal$ be a subset generating $\ideal$.
Then the following conditions are equivalent:
\roster
\item[1]
$\ideal\tA(\HH\1)$ is a two-sided ideal in $\wick(T)$.
\item[2]
$\tA(\HH\1)\ideal\subset\ideal\tA(\HH\1)$.
\item[3]
$\HH\1\ideal_0\subset \ideal+\ideal\HH\1$.
\endroster
If these conditions are satisfied, $\ideal$ is called a
{\bf Wick ideal}.
\eproclaim

\proof:
(1)$\iff$(2): (1) trivially implies (2). For the converse note that
the set $\ideal\tA(\HH\1)$ is clearly invariant under
multiplication with elements of $\tA(\HH\1)$ from the right and
$\tA(\HH)$ from the left. By Wick ordering we have
$\tA(\HH\1)\tA(\HH)\subset \tA(\HH)\tA(\HH\1)$, hence
$$\ideal\tA(\HH\1)\tA(\HH)\subset\ideal\tA(\HH)\tA(\HH\1)
                       \subset\ideal\tA(\HH\1)
\quad.$$
The only property missing to make $\ideal\tA(\HH\1)$ an ideal is the
invariance with respect to multiplication by $\tA(\HH\1)$ from the
left, which is (2).

Suppose (3) holds. Then since $\ideal=\HH\ideal_0\HH$ we get
$\HH\1\ideal\subset(\Cx\idty+\HH\HH\1)\ideal_0\HH
           \subset\ideal+\HH(\ideal+\ideal\HH\1)\HH
           \subset\ideal+\ideal(\Cx\idty+\HH\HH\1)
           \subset\ideal+\ideal\HH\1$.
Hence (3) is equivalent to (3) with $\ideal_0$ replaced by $\ideal$.

(2)$\iff$(3):
The right hand side of (3) is contained in
$\ideal\tA(\HH\1)$. Hence (3) implies (2) by induction on the degree of
monomials in $\tA(\HH\1)$. Conversely, if (2) holds, we have, on the
one hand $\HH\1\ideal\subset \tA(\HH) + \tA(\HH)\HH\1$ by Wick ordering,
and, on the other hand
$\HH\1\ideal\subset\ideal\tA(\HH\1)
                  =\ideal + \ideal\HH\1+\ideal\HH\1\HH\1+\cdots$.
Now each term appearing on the right hand side of these expressions
is Wick ordered. Since the Wick ordered polynomials form a linear
basis of $\wick(T)$, we may compare coefficients of monomials in
adjoint generators $i\1$, and obtain
$\HH\1\ideal\subset \ideal+\ideal\HH\1$.
\QED

Note that the Lemma does not assert that $\ideal\tA(\HH\1)$ is a
star-ideal in $\wick(T)$. Hence in order to lift the involution of
$\wick(T)$ to the quotient we have to divide out the larger ideal
$\ideal\tA(\HH\1)+\tA(\HH)\ideal\1$.
The following Lemma gives a simple criterion for deciding which
coherent \rep s lift to the quotient.

\iproclaim/Iq.co/ Lemma.
Let $\ideal$ be a Wick ideal in a \Wa\ $\wick(T)$, and let
$\pfi:\HH\to\Cx:f\mapsto\bra f,\pfi>$ be a conjugate linear
functional. Denote by $\widehat\pfi:\tA(\HH)\to\Cx$ the induced
homomorphism with $\widehat\pfi(\idty)=1$, and
$\widehat\pfi(f)={\bra\pfi,f>}=\Br{\bra f,\pfi>}$, and denote by
$\creps\pfi$ the separated coherent \rep\ associated with $\pfi$
(cf.~\Def/H.6/). \hfill\break
Then $\creps\pfi(\ideal)=\set0$  if and only if
$\widehat\pfi(\ideal)=\set0$.
\eproclaim

\proof:
By definition of the coherent \rep\ $\crep\pfi$, we have, for all
$G\in\tA(\HH)$, and $f\in\HH$:
$\langle\idty,
             \crep\pfi(f)\crep\pfi(G)\idty\rangle_{T,\pfi}
      = \langle\crep\pfi(f\1)\idty,\
                   \crep\pfi(G)\idty\rangle_{T,\pfi}
      = \Br{\bra f,\pfi>}
         \langle\idty,\crep\pfi(G)\idty\rangle_{T,\pfi}
$.
By induction on the degree of $F$, we have
$$ \Bigl\langle\idty,\crep\pfi(F)\crep\pfi(G)\idty\Bigr\rangle_{T,\pfi}
      = \widehat\pfi(F)\
         \Bigl\langle\idty,\crep\pfi(G)\idty\Bigr\rangle_{T,\pfi}
\quad.$$
In particular, for $G=\idty$, we get
$\bra\idty,\creps\pfi(F)\idty>_{T,\pfi}=\widehat\pfi(F)$. Hence
$\creps\pfi(\ideal)=\set0$  implies $\widehat\pfi(\ideal)=\set0$.
Conversely, if $\widehat\pfi(\ideal)=\set0$, we have
$\bra\idty,
      \creps\pfi(\ideal)\creps\pfi(\tA(\HH))\idty>_{T,\pfi}
    =\set0$,
and, by the Wick ideal property:
$$\align
   \Bigl\langle\creps\pfi(\tA(\HH))\idty, &\
       \creps\pfi(\ideal)\creps\pfi(\tA(\HH))\idty\Bigr\rangle_{T,\pfi}
   = \Bigl\langle\idty,\  \creps\pfi\bigl(\tA(\HH\1)\,
               \ideal\bigr)  \idty\Bigr\rangle_{T,\pfi} \\
   &\subset \Bigl\langle\idty,\
       \creps\pfi(\ideal)\creps\pfi(\tA(\HH\1))\idty\Bigr\rangle_{T,\pfi}
    =\set0
\endalign$$
Since the vectors $\creps\pfi(\tA(\HH))\idty$ span the \rep\ space of
$\creps\pfi$, we get $\creps\pfi(\ideal)=\set0$.
\QED

We now consider ideals which are generated by homogeneous quadratic
expressions. The following result, stated in a slightly different
setup, can be found in \cit{WessZ}.

\iproclaim/Iq.2/ Proposition.
Let $\wick(T)$ be a \HWa. Let $P$ be an orthogonal
projection on $\HH\otimes\HH$, and consider the ideal
$\ideal\subset\tA(\HH)$ generated by $P\HH\otimes\HH$.
Then $\ideal$ is a Wick ideal if and only if
\eqgroup(Iq.id)$$\align
        (1+T)P&=0 \tag{\deq\lasteq a(Iq.lin)}\\
\hbox{and}\qquad
    \bigl(\idty\otimes(\idty-P)\bigr)
    \bigl(T\otimes\idty\bigr)\bigl(\idty\otimes T\bigr)
    \bigl(P\otimes\idty\bigr)
    &=0
\quad, \tag{\deq\lasteq b(Iq.qua)}\endalign$$
where $T$ is the operator introduced in equation \eq(Bt.2). Ideals
$\ideal$ of this form will be called {\bf quadratic ideals}.
\eproclaim

\proof:
The ideal is generated by the elements of the form
$$    A_{ij}= \sum_{k\ell} P_{ij}^{k\ell}\ k\otimes\ell
\quad,$$
where $i,j\in I$. (Typically, these are linearly dependent).
A straightforward computation using \eq(H.1) gives
$$ k\1A_{ij}= \Bigl(P_{ij}^{km}+P_{ij}^{nr}T^{rm}_{kn}\Bigr)\ m
            + P_{ij}^{nm}T_{kn}^{rs}T_{rm}^{uv}\ svu\1
\quad,$$
where we have used the convention that every $i\in I$ appearing twice
is automatically summed over. By \Lem/Iq.1/(3) the first term must
be in $\ideal$, and the second in $\ideal\HH\1$. Since
$\ideal\subset\tA(\HH)$ is generated by homogeneous quadratic
expressions, the first term, being only linear, must be zero, and the
second must be of the form $J_uu\1$ with $J_u\in\ideal$ for all $u$.
Using $P_{ij}^{k\ell}=\me<k\ell,P,ij>$, and \eq(Bt.2), the vanishing of
the first degree term is written as
$$ 0=\me<km,P,ij>
     +\me<km,T,nr>\me<nr,P,ij>
    =\me<km,(\idty+T)P,ij>
\quad,$$
for all $k,m,i,j\in I$.

Since $P$ was chosen as a projection, we can express the condition
$J_u\in\ideal$ as
$$\align
  0&= P_{ij}^{nm}T_{kn}^{rs}T_{rm}^{uv}(\idty-P)_{sv}^{pq} \\
   &= \me<pq,\idty-P,sv>\me<ks,T,nr>\me<rv,T,mu>\me<nm,P,ij>\\
   &= \me<kpq,\idty\otimes(\idty-P),k'sv'>
      \me<k'sv',T\otimes\idty,nrv>
      \me<nrv,\idty\otimes T,n'mu'>
      \me<n'mu',P\otimes\idty,iju>\\
  &=\me<kpq,\bigl(\idty\otimes(\idty-P)\bigr)
            \bigl(T\otimes\idty         \bigr)
            \bigl(\idty\otimes T        \bigr)
            \bigl(P\otimes\idty\bigr),iju>
\quad,\endalign$$
for all $k,p,q,i,j,u\in I$.
\QED

The two conditions on $P$ and $T$ are called the ``linear'' and the
``quadratic'' condition in \cit{WessZ}. Note that any quadratic
ideal is annihilated by the Fock \rep. This gives an immediate
interpretation  of the linear condition: since $P_2=\idty+T$
determines the scalar product in $2$-particle Fock space, the norm
of the generators $A_{ij}$, considered as vectors in Fock space,
must vanish. In particular, a necessary condition for $\wick(T)$ to
have a non-trivial quadratic ideal is that $-1$ is an eigenvalue of
$T$. As \Exa/E.2/ at $qd=-1$ shows, this condition is not
sufficient. However, if the braid relation for $T$ holds (compare
\Exa/E.5/, and Sections \sBq, and \sDb), the quadratic condition
does become redundant. We record this simple observation for future
reference.

\iproclaim/Iq.2b/ Corollary.
Suppose $T$ satisfies the braid condition $T_1T_2T_1=T_2T_1T_2$, and
let $P$ be the eigenprojection of $T$ for eigenvalue $-1$. Then the
ideal in $\tA(\HH)$ generated by $P(\HH\otimes\HH)$ is a quadratic Wick
ideal.

\eproclaim

\proof:
\eq(Iq.lin) is obviously satisfied. Moreover,
$T_2(T_1T_2P_1)=(T_1T_2)T_1P_1=-(T_1T_2P_1)$. Hence $T_1T_2P_1$ maps
into the eigenspace of $T_2$ for eigenvalue $-1$, and
$(\idty-P_2)T_1T_2P_1=0$, proving \eq(Iq.lin). Now \Prp/Iq.2/
applies.
\QED

Let $\ideal$ be a Wick ideal. Then
$\ideal\1\ideal\subset\tA(\HH\1)\ideal\subset\ideal\tA(\HH\1)$, and by
taking adjoints,  $\ideal\1\ideal\subset\tA(\HH)\ideal\1$. This
suggests that $\ideal\1\ideal\subset\ideal\ideal\1$. We do not know
whether this is the case in general. However, for quadratic ideals
this inclusion holds. More precisely, we get the following
relations.

\iproclaim/Iq.3/ Lemma.
Let $T$ and $P$ be as in \Prp/Iq.2/, and set
$$    A_{ij}= P\bigl(i\otimes j\bigr)
            =\sum_{k\ell} P_{ij}^{k\ell}\ k\otimes\ell
\quad,$$
for $i,j\in I$. Then
$$ A\1_{k_1k_2}A_{i_1i_2}
      =\sum_{r_1r_2s_1s_2}
          \me< k_1k_2\,r_1r_2, P_1P_3\, T_2T_1T_3T_2\, P_3P_1,
               i_1i_2\,s_1s_2>
        A_{r_1r_2}A\1_{s_1s_2}
\quad,$$
where as in equation \eq(Bt.5) the indices on $T$ and $P$ indicate
the tensor factors in which the respective operators act. In
particular, $\ideal\1\ideal\subset\ideal\ideal\1$, where $\ideal$ is
the associated Wick ideal.
\eproclaim

\proof: Following the proof of \Prp/Iq.2/ we find
$$\align
   k\1A_{i_1i_2}
      &= \me<ks_1s_2, T_1T_2P_1, i_1i_2u> A_{s_1s_2}u\1    \\
   \ell\1k\1A_{i_1i_2}
      &= \me<k\ell r_1r_2, T_2T_3P_2T_1T_2P_1, i_1i_2uv>
             A_{r_1r_2}v\1u\1    \\
   A_{k_1k_2}\1A_{i_1i_2}
      &= \me<k_1k_2 r_1r_2, P_1T_2T_3P_2T_1T_2P_1, i_1i_2uv>
             A_{r_1r_2}v\1u\1
\quad.\\\endalign$$
We now use the relation $T_1T_2P_1=P_2T_1T_2P_1$ from \Prp/Iq.2/,
together with its translate $T_2T_3P_2=P_3T_2T_3P_2$, and
their adjoints, and the commutation property $X_1X_3=X_3X_1$ for
$X=P,T$, to get
$$\align
  P_1(T_2T_3P_2)T_1T_2P_1
    &= P_1P_3 T_2T_3(P_2T_1T_2P_1)
     = (P_1P_3) T_2 (T_3T_1) T_2P_1 \\
    &= P_3(P_1 T_2T_1)T_3T_2P_1
     = (P_1P_3) T_2T_1(P_2T_3T_2)P_1 \\
    &= P_3P_1\ T_2T_1T_3T_2\ P_3P_1
\quad,\\\endalign$$
where the brackets indicate to which part of the product a rule is
applied in the next step.
Then the result follows by introducing the definition of
$A\1_{s_1s_2}$
\QED

Hence the generators of a quadratic Wick ideal satisfy {\em Wick
commutation relations without constant term}.
When $T$ satisfies the braid condition, the expression
$T_2T_1T_3T_2$, which determines the commutation relations of the
$A_{ij}$, has a simple intuitive interpretation: it describes the
braiding operation on four strands, consisting of taking the strands
two and two together, and braiding these ``cables'' once. In
particular, the coefficients for the algebra $\ideal\ideal\1$
themselves satisfy the braid condition.

Wick relations without constant term behave quite differently from
the Wick relations \eq(H.1). Being homogeneous in the generators,
they always have $A_{ij}=0$ as a solution. On the other hand, if
there is any non-trivial bounded solution, any multiple is also a
solution, so there cannot be a universal bounded \rep. For the
existence of bounded \rep s it is instructive to consider the case
of a single generator, i.e.\ $x\1x=\lambda xx\1$. Since the spectra
of $xx\1$ and $x\1x$ are the same (apart from zero), the spectrum of
$x\1x$ must be invariant under multiplication by integer powers of
$\lambda$. Hence either $x=0$ is the only bounded \rep, or
$\lambda=1$, and $x$ is normal. When we look for bounded \rep s of a
\Wa\ which admits a quadratic ideal, we should therefore look first
for the bounded \rep s of the relations in \Lem/Iq.3/. In the
following section we will study cases in which the ideal is
automatically annihilated in any bounded \rep.

We close this section with a brief review of the Wick ideals in the
examples of Section \sE.
In {\bf\Exa/E.1/}, the $q$-relations and, more generally, for small
$T$ ({\bf\Exa/E.8/}) we typically have the isomorphism of the
universal bounded \rep\ with the Cuntz-Toeplitz algebra. Thus an
\irrep\ is either the Fock \rep\ or a \rep\ of the Cuntz algebra,
which is simple, and hence tolerates no non-trivial ideals. Hence in
these examples we do not have any Wick ideals. Quadratic ideals
cannot exist while $T>-\idty$.

In {\bf\Exa/E.2/}, the Temperley-Lieb-Wick relations, we can choose
$q=-d^{-1}$, so that $-1$ is an eigenvalue of $T$, and the linear
condition \eq(Iq.lin) is satisfied. However, the quadratic condition
\eq(Iq.qua) fails, so these two are independent.

For the twisted canonical commutation relations ({\bf\Exa/E.3/}) we
have a nested sequence of quadratic Wick ideals $\ideal_n$,
generated by $(ij-\mu ji)$, for $i,j\leq n$. In the fermionic case
$\ideal_n$ is generated by $(ij+\mu ji)$, and $i^2$  for
$i,j\leq n$. Wick ideals of higher order, which are contained in
$\ideal_d$ will also occur in the following section.

In {\bf\Exa/E.6/}, the Clifford algebras, we saw that the elements
$\theta_{ij}:=(ij+ji)$ are central. Hence for any choice of a
symmetric, $\Cx$-valued matrix $\widehat\theta_{ij}$, any set of
elements $ij+ji-\widehat\theta_{ij}\idty$ defines a Wick ideal.
Except for $\widehat\theta=0$ this ideal is not quadratic, and since
the generators are not of homogeneous degree, the grading loses its
meaning in the quotient algebra, and becomes replaced by a
$\Ir_2$-grading.

Finally, {\bf\Exa/E.6A/} was constructed in order to get the largest
possible quadratic Wick ideal, containing all monomials of degree
two and higher.

\subhead \ssecno\sIt. Twisted canonical (anti-)commutation relations
\endsubhead

Of the examples from Section \sE\ discussed above, the twisted
canonical (anti-)commu\-tation relations (\Exa/E.3/) have the most
complex Wick ideal structure. In this subsection we will use this
structure to obtain information about the \rep\ theory of these
relations. Since $T$ satisfies the braid relations \cit{Baez} we can
consider, by \Cor/Iq.2b/, the largest quadratic Wick ideal, given by
the ``$-1$''-eigenspace of $T$. In the Bosonic case this ideal is
generated by all $ij-\mu ji$ for $i>j$, and in the Fermionic case
by all $i^2$, and all $ij+\mu ji$ for $i>j$. Pusz and Woronowicz
have determined \cit{PWor,Pusz} all \irrep s of the relations for
which this ideal is annihilated. Therefore, we focus on what can be
said without making that assumption.

We consider the Bosonic case first. The first part of the following
Theorem is based on the idea of considering the \rep s of the
maximal quadratic Wick ideal as an algebra in its own right. As is
often the case for the \rep\ of \Wa\  relations without
constant term, we find that only the zero \rep\ is bounded. The
further classification of bounded \rep s then proceeds along the
lines set by Pusz and Woronowicz \cit{PWor} in their classification
of the \rep s in which the Wick ideal is annihilated. Using the
notion of coherent \rep s their classification can be made more
transparent, which helped us to detect an omission in their
classification (in the notation of the Theorem they have only
$\alpha=1$). Therefore we include a brief, corrected sketch of the
classification of the bounded \rep s as it might have been given in
\cit{PWor}. We emphasize, however, that there are also {\em
unbounded} \rep s of the relations, which are classified in
\cit{PWor}, but which are {\em not} coherent in our sense.

\iproclaim/It.1/ Theorem.
Let $\pi$ be a bounded \rep\ of the twisted canonical commutation
relations \eq(E.3). Then
$$   \pi(ij)=\mu \pi(ji)
\quad,$$
for all $i>j\in I$. Moreover, if $\pi$ is irreducible, it is
coherent, i.e.\ there are a cyclic unit vector $\Omega$, and
$\pfi_i\in\Cx$ such that $\pi(i\1)\Omega=\pfi_i\Omega$. One has
either $\pfi_i=0$ for all $i$ (the Fock \rep), or there is an index
$k$, and a phase $\alpha\in\Cx$, $\abs\alpha=1$, such that
$$ \pfi_i=\cases \alpha\,(1-\mu^2)^{-1/2}& i=k\\
                    0                    & i\neq k\endcases
\quad.\tag\deq()$$
Representations with different $k$ or $\alpha$ are inequivalent.
\eproclaim

\proof:
The standard form of the generators given in \Lem/Iq.3/ is
$A_{ij}=c(ij-\mu ji)$, with an irrelevant overall factor
$c=-\mu/(1+\mu^2)$, which we will drop in the sequel.
A tedious but straightforward computation gives, for $i>j$, the
following explicit version of the commutation relations found in
\Lem/Iq.3/:
$$\align
   A\1_{ij}A_{ij}
       &=\mu^6\aaa ij
         -\mu(1-\mu^2)         \sum_{r<j}    \aaa ir  \\
       & -\mu^4(1-\mu^2)       \sum_{j<s<i}  \aaa sj
         +\mu^2(1-\mu^2)       \sum_{r<j<s<i}\aaa sr   \\
       & -\mu^4(1-\mu^2)       \sum_{r<j}    \aaa jr   \\
       & +(1-\mu^2)^2(1+\mu^2) \sum_{r<s<j}  \aaa sr
\quad.\\\endalign$$
Here we have arranged the terms $A_{nm}A\1_{nm}$ on the right hand
side in descending lexicographic order with respect to the pairs
$(n,m)$ with $n>m$. Note that only pairs $(n,m)$ which are smaller
than $(i,j)$ in this order contribute to this expression. Since for
a bounded operator $A\1_{ij}A_{ij}=\mu^6 A_{ij}A\1_{ij}$ implies
$A_{ij}=0$, the result follows by induction on lexicographic order.

We now give a sketch of the further classification of bounded \irrep
s. Consider the operators $E_i=\pi(ii\1)$. Then, because of the
relation just proven, these operators commute. Moreover,
$\pi(i\1)E_j=\tau_i(E)_j\pi(i\1)$, where $\tau_i:\Rl^d\to\Rl^d$ is
an affine transformation:
$$ \tau_i(E)_j=\cases
        E_j & i<j\\
        1+\mu^2 E_i-(1-\mu^2)\sum_{k<i}E_k   &i=j \\
        \mu^2 E_j &i>j \endcases
\quad.$$
Let $f:\Rl^d\to\Rl$ be a measurable function with the property
$f\circ\tau_i=f$ for all $i$. Then $f(E)$ (evaluated in the joint
functional calculus of the $E_i$) commutes with all generators and
their adjoints, hence is a multiple of the identity. It follows that
the joint spectrum $\sigma(E)\subset\Rl^d$ of the $E_i$ consists of
a single orbit of the $\tau_i$. Let $\lambda\in\sigma(E)$, and let
$\Phi$ be a vector in the corresponding eigenspace:
$E_i\Phi=\lambda_i\Phi$. Then $\pi(i\1)^n\Phi$ is either zero or an
eigenvector for the eigenvalues $\tau_i^{-n}(\lambda)$.

Consider the iteration of $\tau_1^{-1}$. Then, as a first case, we
may have $\tau_1(\lambda)=\lambda$, i.e.\
$\lambda_1=(1-\mu^2)^{-1}$, and $\lambda_i=0$ for $i>1$.
If $\lambda$ is not a fixed point of $\tau_1$, the iteration of
$\tau_1^{-1}$ produces an unbounded sequence, hence, eventually, we
must have $\pi(1\1)^n\Phi=0$ (Here ``$1$'' denotes the first
generator). We may now replace $\Phi$ by the last non-zero vector
$\pi(1\1)^n\Phi$, and repeat the argument for the iteration of
$\tau_2$, with the additional information that $\lambda_1=0$.
By iteration we obtain either an eigenvector $\Phi$ of $E$ such that
$\lambda_i=0$ for all $i$, which clearly is a Fock vector. Or,
otherwise, we find $k$, and an eigenvector $\Phi$ of $E$
such that $\lambda_k=(1-\mu^2)^{-1}$, and $\lambda_i=0$ for
$i\neq k$.

Let $\NN$ be the subspace of eigenvectors $\Phi$ of the latter type.
Then since $\lambda$ is a fixed point of $\tau_k$, we have that
$\pi(k)\NN\subset\NN$, and $\pi(k\1)\NN\subset\NN$. Moreover, the
restriction of $\pi(k)$ to $\NN$ is unitary up the factor
$(1-\mu^2)^{1/2}$. Let $U$ be a unitary operator in $\NN$ commuting
with the restriction of $\pi(k)$. We extend $U$ to an operator on
the cyclic subspace generated by $\NN$, i.e.\ the whole \rep\ space,
by setting
$$U\pi(X)\Phi=\pi(X)U\Phi
\quad,$$ for
$\Phi\in\NN$, and $X\in\wick(T)$. We claim that $U$ is
well-defined, and unitary.  In
\ppt{\break}%
$\bra \pi(X_1)U\Phi_1, \pi(X_2)U\Phi_2>$ we can Wick order the
expression $X_1\1X_2$. From the Wick ordered form of $X_1\1X_2$ we
can drop all terms containing generators $i\neq k$, because of the
commutation relations between the $\pi(i)$, and because
$\pi(i\1)\Phi_\nu=0$, $\nu=1,2$. This leaves a polynomial in
$\pi(k)$ and $\pi(k\1)$, which commutes with $U$. Hence
$\bra \pi(X_1)U\Phi_1, \pi(X_2)U\Phi_2>
  =\bra \pi(X_1)\Phi_1, \pi(X_2)\Phi_2>$, from which it is clear
that $U$ is unitary. Moreover, by definition, $U$ commutes with all
$\pi(X)$, and must be a multiple of the identity by irreducibility.
It follows that $\NN$ is one-dimensional, and hence that
$\pi(k)\Phi=z\Phi$, for some $z\in\Cx$, and $\Phi\in\NN$. From
$\pi(kk\1)\Phi=\lambda_k\Phi=(1-\mu^2)^{-1}\Phi$ we get
$z=\alpha(1-\mu^2)^{-1/2}$, and any unit vector in $\NN$ is the
desired coherent vector. The inequivalence of the \rep s with
different $k$ or $\alpha$ is evident from this construction.
\QED

In the above proof we have not established the positivity of the
coherent \rep s. This can be done in a simple way, which gives
additional insight into this structure: let $\pi'$ denote the Fock
\rep\ of the relations with $k-1$ generators, with Fock vector
$\Omega$. Let $N$ denote the number operator of this \rep, defined
by $N\pi'(i_1i_2,\dots,i_n)\Omega=n\pi'(i_1i_2,\dots,i_n)\Omega$. Then
we set
$$ \pi(i)=\cases
         \pi'(i) & i<k\\
         \alpha (1-\mu^2)^{-1/2}\, \mu^N    &i=k\\
         0       & i>k \endcases
\quad.\tag\deq()$$
The relations for $\pi(ij)$, and $\pi(ij\1)$, with $i\neq j$ readily
follow from the definition of the number operator. The only
non-trivial relation to be checked is for $\pi(kk\1)$, which reduces
to the identity
$$ \sum_{i<k}\pi(ii\1)=\frac{1-\mu^{2N}}{1-\mu^2}
\quad,\tag\deq()$$
which was shown in \cit{PWor}, equations 2.34 and 2.35. Clearly,
this \rep\ is positive, because the Fock \rep\ $\pi'$ is positive,
and it is immediately clear that it is the coherent \rep\ with the
parameters specified in the Theorem.

The Fermionic case is more involved. Since the first generator will
play a special r\^ole, and ``$1$'' is a confusing notation for the
first generator, we use the notational convention introduced in
equation \eq(H.4), and denote the generators in the \rep\ ``$a\1$''
under consideration by $a\1_i=a\1(i)$, and their adjoints by
$a_i=a\1(i\1)$, with $i=1,\ldots,d\leq\infty$. The relations
\eq(E.4) then become
$$
a_ia\1_j=\cases -\mu\ a\1_ja_i   &\text{ for $i\neq j$}\\
                  \idty - a\1_ia_i
                      -(1-\mu^2)\sum_{k<i}a\1_ka_k
                                 &\text{ for $i=j$}.\\
\endcases\tag\deq(tcar)$$
We will refer to these relations as the \mCAR.
It is clear that each $a_i$ is bounded, and we need not distinguish
the abstract \Wa ic adjoint ``$\1$'' from the operator adjoint $^*$.
As in the Bosonic case, we have a nested sequence of quadratic Wick
ideals $\ideal_n$, but this time generated by $i^2$, and $ij+\mu ji$
for $1\leq i,j\leq n$.

The following Proposition is the basic tool for investigating these
nested ideals in a general \rep: it allows us to reduce the case
that $\ideal_1$ is annihilated to the study of a another \rep\ of
the \mCAR\ with one generator less. Applied inductively, it reduces
the study of \rep s in which $\ideal_n$ vanishes to the study of the
\mCAR\ with $n$ generators less. In particular, the \rep s in which
the maximal quadratic Wick ideal $\ideal_d$ is annihilated, i.e.\
those \rep s considered by \cit{Pusz}, are analyzed completely by
the following Proposition, and are given as
$$ a_i= \underbrace{\pmatrix -\mu&0\\0&1\endpmatrix \otimes\dots\otimes
                    \pmatrix -\mu&0\\0&1\endpmatrix}
                                 _{(i-1)\text{ factors}}
         \otimes\pmatrix 0&0\\1&0\endpmatrix
         \otimes \underbrace{\idty\otimes\dots\otimes\idty}
                                 _{(d-i)\text{ factors}}
\quad.\tag\deq()$$
This is precisely the $\mu$-deformation of the standard technique
for analyzing the canonical anticommutation relations with $\mu=1$
\cit{BraRo}. In the general case we get a decomposition of an
arbitrary \rep\ $\pi$ of the \mCAR\ as $\pi=\pi_0\oplus\pi_1$ with
$\pi_0(\ideal_1)=\set0$, or, equivalently $\pi_0(a_1^2)=0$. A
crucial property needed for this decomposition is the normality of
$a_1^2$, which follows from the \mCAR. Unless it vanishes, $a_1$
will not be normal. For the general structure theory of (single)
operators with normal square, we refer to \cit{Radjavi}.

\iproclaim/It.2/ Proposition.
Let $a_i$, $i=1,\ldots,d$ be operators on a Hilbert space $\R$
satisfying the \mCAR\ \eq(tcar).
\roster
\item[1]
Then $a_1^2$ is normal, and
both summands in the direct sum decomposition $\R=\R_0\oplus\R_1$, with
$\R_0=\ker a_1^2$, and $\R_1=\Br{a_1^2\R}$,
are invariant under all $a_i,a_i^*$.
\item[2]
$\R_0$ can be decomposed into $\R_0=\Cx^2\otimes\R_0'$ such that
$$\align
   a_1&=\pmatrix 0&0\\1&0\endpmatrix   \otimes\idty \\
   a_i&=\pmatrix -\mu&0\\0&1\endpmatrix\otimes \widetilde a_{i-1}
\quad,\quad\hbox{for $i>1$}, \endalign$$
where the operators $\widetilde a_i$, $i=1,\ldots, d-1$
again satisfy \eq(tcar).
\endroster
\eproclaim

\proof:
(1)
Since $a_1^*a_1+a_1a_1^*=\idty$, $a_1$ determines a \rep\ of the
Clifford algebra discussed above in Section \sIq, and $N=a_1^2$ is
normal. Then for any $j>1$:
$$ N a_j^*=\mu^2 a_j^*N
\quad.$$
By the Fuglede-Putnam-Rosenblum Theorem \cit{Rudin},
$N^*a_j^*=\mu^2a_j^*N^*$, or $a_jN=\mu^2Na_j$. Hence the kernel of $N$
is an invariant subspace for all generators and their adjoints. The
orthogonal complement of the kernel is then also an invariant
subspace.

(2) In the following we may assume that $\R_1=\set0$, and $\R_0=\R$.
Since $a_1^2=0$, $a_1$ and its adjoint generate a copy of the
$2\times2$ matrix algebra, and we can decompose the space as
$\R=\Cx^2\otimes\R'$, with $a_1$ of the form given in the
Proposition. What we have to show is the form for the $a_i$, $i>1$.
We can write $a_i=\pmatrix C_i&B_i\\ D_i&A_i\endpmatrix$, with
bounded operators $A_i,B_i,C_i,D_i$ on $\R'$. Then
$$ 0=a_ia_1^*+\mu a_1^*a_i
        = \pmatrix \mu D_i& C_i+\mu A_i\\ 0& D_i\endpmatrix
\quad.$$
Hence $D_i=0$, and $C_i=-\mu A_i$. The diagonal block matrix
components of the Wick relation for $a_ia_i^*$ ($i>1$) are
$$\align
    \mu^2(A_iA_i^*+A_i^*A_i) + B_iB_i^*
         &=\idty-(1-\mu^2)\mu^2\sum_{1<k<i}A_k^*A_k
             -(1-\mu^2)\idty\\
    A_iA_i^*+A_i^*A_i + B_i^*B_i
         &=\idty-(1-\mu^2)\sum_{1<k<i}(B_k^*B_k+A_k^*A_k )
\quad,\\\hbox{and\ }\qquad
    B_iB_i^*-\mu^2B_i^*B_i
         &=\mu^2(1-\mu^2)\sum_{1<k<i} B_k^*B_k
\quad,\\\endalign$$
where the third equation is the first minus $\mu^2$ times the
second. For $k=2$ we get $B_2B_2^*=\mu^2B_2^*B_2$, hence $B_2=0$ by
boundedness. Proceeding by induction we find that $B_i=0$ for all
$i>1$.
Hence $a_i=\pmatrix -\mu A_i&0\\0&A_i\endpmatrix$. Substituting
$B_i=0$ in the above equation we find that the $A_i$ with $i>1$
satisfy the Wick relation for $A_iA_i^*$. Since the block matrix for
$a_i$ is diagonal it is clear that the relations for $a_ia_j^*$ with
$i,j>1$ are satisfied if and only if the relations for $A_iA_j^*$
hold.
\QED

It remains to analyze the \rep s of the type $\pi_1$.
As a first step we show that any \irrep\ of this type is coherent,
with a vector $\Omega$ satisfying $a_1\Omega=\alpha\Omega\neq0$,
and $a_i\Omega=0$ for $i>1$. Combined with the previous Proposition
we get the following statement.

\iproclaim/It.3/ Theorem.
\roster
\item[1]
Every \irrep\ of the \mCAR\ \eq(tcar) is coherent.
\item[2]
For every coherent \rep s there are an integer $r$ and
$\alpha\in\Cx$ with $\abs\alpha\leq\mu^{r-1}/\sqrt2$ such that the
cyclic vector $\Omega$ is characterized by
$$ a_j\ \Omega=\alpha\ \delta_{jr}\ \Omega
\quad.\tag\deq(It.3)$$
\endroster
\eproclaim

\proof:
By \Prp/It.2/, an \irrep\ has either $a_1^2=0$, or
$a_1^2$ non-singular. In the first case the explicit form for $\R_0$
makes clear that $a$ is irreducible iff $\widetilde a$ is, and if
$\widetilde a$ is coherent with cyclic vector $\Omega$, then $a$ is
coherent with cyclic vector
$\widetilde\Omega
    ={\left(\smallmatrix 1\\0\endsmallmatrix\right)}\otimes\Omega$.
Moreover, $\widetilde a_i\widetilde\Omega
      =\widetilde\alpha_i\widetilde\Omega$ implies
$a_i\Omega=\alpha_i\Omega$  with $\alpha_1=0$, and
$\alpha_i=-\mu\widetilde\alpha_{i-1}$ for $i>1$.
Hence, by induction on the number of generators, we may assume in both
parts of the Theorem that $\R_0=\set0$, i.e.\ that $N=a_1^2$ is
non-singular, and hence $r=1$.

We have seen in the proof of \Prp/It.2/ that $a_jN=\mu^2Na_j$, and
$a_jN^*=\mu^2N^*a_j$. Hence by the functional calculus
$$   a_j\ f(N)=f(\mu^2N)\ a_j
\quad,$$
where $j>1$, and $f$ is any measurable real-valued function on the
complex plane. In particular, let $f$ be a real-valued function with
the scaling invariance property $f(x)=f(\mu^2x)$, and set $F=f(N)$.
Then $F$ commutes with all generators, including $a_1$. Hence, in an
\irrep, it must be a multiple of the identity. It follows that the
spectrum of $N$ must be contained in just one of the sets
$\set{\mu^{2n}\alpha^2\stt n\in\Ir}$, where $\alpha\in\Cx$. Since
$N$ is bounded, we can take $\alpha^2$ as the element of largest
modulus, and obtain a spectrum contained in
$\set{\mu^{2n}\alpha^2\stt n\geq0}$.

Let $\pfi$ be in the eigenspace $\NN$ of $N$ for the largest
eigenvalue $\alpha^2$. Then $N^*\pfi=\Br\alpha^*\pfi$, and, for
$j>1$: $N^*\,a_j\pfi=\mu^{-2}\Br\alpha^2\,a_j\pfi$. Since $a_j\pfi$
cannot be eigenvector of $N$ for an eigenvalue larger than
$\norm{N}$, we must have $a_j\pfi=0$. Since $\NN$ is also
invariant under $a_1$, either $(a_1-\alpha)\NN$ or $(a_1+\alpha)\NN$
is non-zero, hence $\NN$ contains an eigenvector of $a_1$ with
eigenvalue $\pm\alpha$. Renaming $\alpha\mapsto-\alpha$, if
necessary, we obtain the vector described in (2). The bound on
$\alpha$ comes from the \rep\ theory of the Clifford algebra
(see \Exa/E.6/, and \cit{QCS}).

The cyclic sub\rep\ generated by this vector is the coherent \rep,
and since we have started from an \irrep\ this sub\rep\ has to be
the given one.
\QED

Note that the Theorem does {\em not} claim that all the
coherent \rep s with specified $r,\alpha$ are positive. We strongly
suspect that this is automatically the case, but were not able to
prove it in full generality. However, some information about the
coherent \rep s can be obtained with the methods we have introduced.
{}From the proof of \Prp/It.2/, we know that in every positive \rep\ of
\eq(tcar) the relation $a_1^{*2}a_j^*=\mu^2 a_j^*a_1^{*2}$ holds.
One readily verifies that this relation generates a Wick ideal.
An even larger Wick ideal, which is also automatically annihilated
is given in the following Proposition.

\iproclaim/It.4/ Proposition.
The \Wa\ determined by the \mCAR\ \eq(tcar) has a Wick ideal,
generated by the elements
$$\matrix\format\r&\l&\l\\
  \mo[1,1,i ] &- \mu^2 \mo[i,1,1]                      &\quad1<i \cr
  \mo[1,i,i] &+
         \mu^{-1}\mo[i, 1, i] - \mu \mo[i,1,i] -
         \mo[i,i,1]                                     &\quad1<i\cr
  \mo[1,i,j] &+
        \mu^{-1}\mo[i, 1, j] - \mu^2 \mo[j, 1, i] -
        \mu \mo[j, i, 1]                             &\quad 1<i<j\cr
  \mo[1,j,i] &-
        \mu^2 \mo[i, 1, j] - \mu \mo[i, j, 1] -
        ({\tsize-\mu^{-1} + \mu - \mu^3}) \mo[j, 1, i] -
        (1 - \mu^2) \mo[j, i, 1]                     &\quad 1<i<j\cr
\endmatrix$$
This ideal is annihilated in every positive \rep. Consequently the
coherent \rep\ spaces are spanned by vectors of the form
$$   a_{i_1}^*\cdots a_{i_n}^*\Phi_i
\quad,\tag\deq(It.repsp)$$
with $i,i_\alpha>1$, and
$$ \Phi_i\in\set{a_i^*\Omega,\ a_i^*a_1^*\Omega,\
                 a_1^*a_i^*\Omega,\ a_1^*a_i^*a_1^*\Omega }
\quad.$$
\eproclaim

\proof: The computations necessary to verify that these elements
generate a Wick ideal are so tedious and so straightforward, that we
performed them in {\smc Mathematica} \cit{Mathematica}. The program
is available by anonymous {\tt ftp} (see the Introduction). That
this ideal is annihilated in every positive \rep\ follows
immediately from \Prp/It.4/ and \Lem/Iq.co/. By using the generators
of the ideal as substitution rules, we can transform any monomial in
the $a_i$ to a form in which all factors $a_1$ are collected on the
left, apart from at most one intervening $a_i$ with $i>1$. Applying
the adjoint of this statement to the coherent vector $\Omega$, we
obtain the second statement.
\QED

Further computations in {\smc Mathematica} showed that the scalar
product in the coherent \rep\ space is indeed strictly positive on
some subspaces with low $n$ in \eq(It.repsp). However, we found no
way to show this for general $n$. It would be particularly
interesting to decide positivity for an infinite number of
generators.

\head    \secno\sD.  Wick algebra relations as differential calculus
\endhead

Non-commutative differential calculus is a fast growing subject, and
we cannot even begin to review it
\cit{Baez,Cuntzb,Manin,Rosenberg,Worona}. There are some natural
links with the structure of \Wa s, however, and this section is
devoted to describing these. In the first two subsections we will
consider two prima facie different approaches to generalizing
commutative differential calculus, one based on a generalization of
the partial derivatives, the other on a generalization of the
differential forms. In both cases Wick type relations appear
naturally, and given these, the two approaches amount to the same
thing. In the third subsection we will discuss two characterizations
of the braid relations for $T$ in terms of the resulting calculus,
based on an idea of \cit{WessZ}.

\subhead \ssecno\sDf. The algebra of differential operators
\endsubhead

Let us consider the tensor algebra $\A\equiv\tA(\HH)$ as a
non-commutative analogue of an algebra of ``functions'' generated by
the coordinate functions $x_i$, $i\in I$. How can we generalize the
commutative differential calculus to this algebra? The first
approach to this problem  is to build an analogue of the algebra
{\em differential operators} on an abelian algebra. This is an
associative algebra with identity, which is generated from the
coordinates, considered as multiplication operators,  together with
the partial derivative operators $\partial_i$ for $i\in I$. The
structure of this algebra is determined by the commutation rule
between coordinates and partial derivatives, and it is natural to
assume
$$ \partial_ix_j=\delta_{ij}\idty
         + \sum_{k,\ell\in I} T_{ij}^{k\ell}\ x_\ell \partial_k
\quad,\tag\deq(D.1)$$
We notice that this is precisely our basic \Wa\  commutation
rule with the substitution $i\mapsto x_i,\ i\1\mapsto\partial_i$.
We should point out that if we think of these relations as defining
a differential calculus, it is not clear why we should have a
\HWa. In the abelian case, we have the canonical commutation
relations, which do define a \HWa. Its basic positive \rep\ is the
Bargmann \rep\ \cit{Bargmann}. From this point of view it may be
suggestive to assume hermiticity, but we will not do so for the
moment.

The view of \Wa s as differential calculi suggests some structures
to investigate, which are useful for the study of \Wa s\
independently of this interpretation. For $f\in\tA(\HH)$, we can Wick
order the product $\partial_i f$. This will produce a constant term
and one containing exactly one partial derivative on the right:
$$ \partial_i f= D_i(f) + \sum_\ell\twist_i^\ell(f)\ \partial_\ell
\quad.\tag\deq(D.2)$$
Here $D_i$ and $\twist_i^\ell$ are linear operators on $\tA(\HH)$.
One immediately gets the ``twisted derivation'' properties
$$\aligned
     D_i(fg)&= (D_if)g +  \sum_\ell\twist_i^\ell(f)\ D_\ell(g) \cr
     \twist_i^\ell(fg)&= \sum_k \twist_i^k(f)\twist_k^\ell(g)
\quad,\cr\endaligned\tag\deq(D.3)$$
The second relation can be stated by saying that the {\em twist}
$\twist:\tA(\HH)\to\Mat_d\bigl(\tA(\HH)\bigr)$ is a homomorphism into
the algebra of $d\times d$-matrices over $\tA(\HH)$, where $d=\abs I$
is the number of generators. We can also summarize the two relations
\eq(D.3) by saying that $\widehat\twist:\tA(\HH)\to\Mat_{d+1}(\tA(\HH))$
is a homomorphism, where the matrix elements of $\widehat\twist$ are
indexed by $I\cup\set0$, and are explicitly defined by
$$  \widehat\twist_i^\ell(f)
          =\cases  \twist_i^\ell(f)   & i,\ell\cr
                        D_i(f)        & \ell=0, i\in I \cr
                          0           & i=0 \quad.
\endcases\tag\deq(D.3')$$
In any case, it is clear that $D$ and $\twist$ can be computed by
from \eq(D.3), starting from the initial values
$$\aligned
               D_i\idty=0 \
\qquad&,\qquad D_i(x_j)= \delta_{ij}  \cr
               \twist_i^\ell(\idty)= \delta_{i\ell}
\qquad&,\qquad \twist_i^\ell(x_j)= \sum_k T_{ij}^{\ell k}\ x_k
\quad.\cr\endaligned\tag\deq(D.4)$$

To include the important case with algebraically dependent
coordinates, we pass to a quotient $\tA(\HH)/\ideal$ over a suitable
ideal. We can still use the same commutation rules, provided the
derivative and twist pass to the quotient, i.e.\ provided
$D_i(\ideal)\subset\ideal$, and
$\twist_i^\ell(\ideal)\subset\ideal$. By \Lem/Iq.1/(3), this
is precisely the definition of a {\em Wick ideal}.

\ppt{\vfill\eject}
\subhead \ssecno\sDd. The algebra of differential forms \endsubhead

A second independent route to a non-commutative differential
calculus starts from the {\em differential forms} rather than the
differential operators. The relevant structure here is that of a
{\em graded differential algebra\/} \cit{Connes,Cuntzb}. The
universal differential algebra over an algebra $\A$ (here:
$\tA(\HH)$) is defined as the associative algebra $\Dff(\A)$
generated by the elements $f$ and $\df f$ for $f\in\A$, where
$\df f$ depends linearly on $f$, $\df\idty=0$, and
$\df(fg)=(\df f)g+f\df g$. By applying this rule, every element of
$\Dff(\A)$ can be written in the form $\omega=f_0\df f_1\cdots\df
f_p$. Such elements are defined to have degree $p$, for which we
write $\degree\omega=p$, or $\omega\in\Dff^p(\A)$. There is a unique
extension of $\df$ to a linear operator on $\Dff(\A)$, satisfying
$d^2=0$, and the graded derivation property
$$ \df(\omega_1\omega_2)
      =\df(\omega_1)\omega_2+(-1)^{\degree{\omega_1}}
       \omega_1\,\df(\omega_2)
\quad,\tag\deq(D.5)$$
applied whenever $\omega_1$ is of homogeneous degree.

Note that while we can write every element as a linear combination
of simple forms $f_0\df f_1\cdots\df f_p$, we cannot insist on a
similar expansion in coordinate differentials $\dx i$, since in the
universal algebra $\Dff(\A)$ there is no commutation relation
between $x_i$ and $\dx j$. Let us postulate a relation of this kind in
the general form
$$ x_j\ \dx\ell = \sum_{ik}T_{ij}^{\ell k}\  (\dx i)\ x_k
\quad.\tag\deq(D.6)$$
We will denote by  $\DffT=\DffT(\tA(\HH))$ the quotient of
$\Dff(\tA(\HH))$ by this relation and its differential
$$ \dx j\ \dx\ell + \sum_{ik}T_{ij}^{\ell k}\  \dx i\ \dx k =0
\quad.\tag\deq(D.7)$$
Thus $\DffT$ is the graded differential algebra with the universal
property that for any homomorphism $\eta:\tA(\HH)\to\Dff'$ into
another graded differential algebra $\Dff'$, with the images
$x_i'=\eta(x_i)$ satisfying \eq(D.6), there is a unique graded
homomorphism $\hat\eta:\DffT\to\Dff'$ extending $\eta$, and
satisfying $\df\hat\eta=\hat\eta\df$.

By induction on the degree of a polynomial $f\in\tA(\HH)$, we then find
in $\DffT$ relations of the form
$$\aligned
          \df f&=\sum_i \dx i\ D_i(f)  \cr
    f\ \dx\ell &= \sum_i (\dx i)\ \twist_i^\ell(f)
\quad,\endaligned\tag\deq(D.8)$$
where $D_i$ and $\twist_i^\ell$ are the linear operators on
$\tA(\HH)$, which are inductively defined by equations \eq(D.4) and
\eq(D.3). It is clear that we can invert this construction, and define
the algebra of differential forms $\DffT$ starting from the
algebra of differential operators $x_i, \partial_i$, so these two
approaches are essentially equivalent.

A consequence of equation \eq(D.7) is that, for many of the \Wa s
considered in Section \sE, the differential calculus collapses after
the $1$-forms, i.e.\ all forms of degree $\geq2$ vanish. A more
precise description of the space of $p$-forms is given by the
following Lemma.

\iproclaim/Dd.fcs/ Lemma.
Let $\Dff^p_0$ denote the space of $p$-forms with constant
coefficients, i.e.\ the linear subspace of  $\DffT$ generated by all
expressions $\dx{i_1}\cdots\dx{i_p}$. Then $\Dff^p_0$ is canonically
isomorphic to
$$  \bigcap_{r=1}^{p-2}\HH^{\otimes (r-1)}\otimes(\ker(\idty+T^*))
                    \otimes\HH^{\otimes(p-r-1)}
      \quad\subset\quad \HH^{\otimes p}
\quad.\tag\deq(Dd.fcs)$$
\eproclaim

\proof:
$\Dff^p_0$ is the quotient of the space of linear combinations
$\sum\Phi(i_1,\dots,i_p)\dx{i_1}\cdots\dx{i_p}$ by the subspace
generated by the coefficients $\Phi$ of the form
$$\sum_{k\ell}\bra i_ri_{r+1}\abs{\idty+T} k\ell>
       \Psi(i_1,\dots,i_r,k,\ell,i_{r+2},\dots,i_p)
   =\bigl((\idty+T_r)\Psi\bigr)(i_1,\dots,i_p)
\quad.$$
This subspace is the span of the ranges of the operators
$(\idty+T_r)$. Since in a Hilbert space we can identify the
quotient by a subspace with its orthogonal complement, we get
$$ \Dff^p_0\cong\Bigl(\bigcup_r(\idty+T_r)
                   \HH^{\otimes p}\Bigr)^\perp
           =\bigcap_r\bigl((\idty+T_r)\HH^{\otimes p}\bigr)^\perp
\quad.$$
\QED

Let us consider the space $\Dff^p_0$ in the examples. In {\bf
\Exa/E.1/} the kernel of $\idty+T$ is zero, so $\Dff^p_0=\set0$ for
all $p\geq2$, and the same conclusion holds for any \Wa\ with small
$T$ ({\bf \Exa/E.8/}). In {\bf\Exa/E.2/} we may choose $q=-d^{-1}$,
and get
$$ \dx i\,\dx j=\delta_{ij}\ {1\over d}\sum_k \dx k\,\dx k
               \equiv \delta_{ij}\ \omega
\quad.\tag\deq(Dd.tl)$$
On the other hand, one easily verifies that
$\omega\,\dx i=d^{-1}\dx i\,\omega=d^{-2}\omega\,\dx i$.
Hence $\Dff^2_0=\Cx\omega$, and $\Dff^p_0=\set0$ for $p\geq3$.

For  the twisted canonical commutation relations, i.e.\ the Bosonic
case of {\bf \Exa/E.3/}, the kernel of $\idty+T$ is generated by the
vectors $\ket{ij}-\mu\ket{ji}$, with $i<j$. In particular, any
form in which some $\dx i$ appears twice vanishes. As in the
commutative case, which corresponds to $\mu=0$, we get
$\dim\Dff^p_0={d\choose p}$, where $d$ is the number of generators.
In the Fermionic case (twisted anti-commutation relations) we get
the same combinatorial problem as in the determination of the
dimension of the Bosonic Fock space over a $d$-dimensional
one-particle space. Thus
$$ \dim \Dff^p_0=(-1)^p{-d\choose p}={d+p-1\choose p}
\quad.\tag\deq(Dd.tcar)$$
Since the Clifford algebras ({\bf \Exa/E.6/}) are the special case with
$\mu=1$, this formula also holds in that case.
These examples follow a general pattern: the sign in the graded
derivation property \eq(D.5) forces the differentials of commuting
variables to anti-commute, and conversely. Moreover, this duality
persists under $\mu$-deformation.

One consequence of \eq(Dd.tcar) is that $\Dff^p_0\neq\set0$ for all
$p$. A more trivial example with this property is {\bf\Exa/E.6A/},
where $\dim\Dff^p_0=d^p$ even grows exponentially. However, if we
replace $\ker(\idty+T^*)$ by a generic subspace of $\HH\otimes\HH$
the intersection in \eq(Dd.fcs) vanishes for large $p$. The
geometric property of this subspace making $\Dff^p_0\neq\set0$ for
all $p$ has been investigated recently as a property of finite range
interactions of quantum spin chains, where it leads to a new class
of exactly solvable ground state problems \cit{FCL,FCS}. A paradigm
of a subspace with this property \cit{AKLT} is the sum of the
spin-$0$ and spin-$1$ subspaces of a spin-$1$ chain. This
corresponds to a \Wa\ on three generators with relations
$$ i\1j= \delta_{ij}\Bigl(1-\frac{\lambda-2}2 \sum_kkk\1\Bigr)
           +\frac\lambda2 ji\1 -\frac\lambda3 ij\1
\quad.\tag\deq(Dd.aklt)$$
In that case we have $\dim\Dff^p_0=4$ for all $p\geq2$.

\subhead \ssecno\sDb. Differential calculus with braid relations
\endsubhead

Wess and Zumino \cit{WessZ} have proposed to ``complete the
algebra'' generated by the three kinds of objects $x_i, \partial_i$,
and $\dx i$ by introducing commutation relations between
$\partial_i$ and $\dx i$. This would indeed be natural if we want to
upgrade the partial derivatives $\partial_i$ to ``{\em covariant
derivatives}'' acting also on differential forms. They make the
ansatz
$$\partial_i\dx j= \sum_{k,\ell} S_{ij}^{k\ell}\
                   \dx\ell\ \partial_k
\quad,\tag\deq(Db.wz)$$
where the $S_{ij}^{k\ell}$ are complex coefficients.
If we apply relations \eq(D.1), \eq(D.6), and \eq(Db.wz) to
$\partial_i\,x_j\, \dx k$  we have two different routes for applying
the rules which agree if and only if
\eqgroup(Db.b1)$$\align
    \sum_\ell \dx\ell &\me<i\ell, ST-\idty, jk>=0
\tag{\deq\lasteq a(Db.b1a)}\\
    \sum_{\ell nm} \dx\ell\, x_m \partial_m
       &\me<i\ell n ,T_2T_1S_2-S_1T_2T_1, jkm> =0
\quad,\tag{\deq\lasteq b(Db.b1b)}\endalign$$
where $S:\HH\otimes\HH\to\HH\otimes\HH$ has the matrix elements
$$    \me<i\ell,S,jk>= S_{ij}^{k\ell}
\quad.\tag\deq(D.10)$$
Now if we want the differentials $\dx\ell$, as well as the
combinations $\dx\ell\, x_m \partial_m $ to be linearly independent
in the resulting algebra, we have to take $ST=\idty$ from
\eq(Db.b1a), and taking the number of generators finite, we also get
$TS=\idty$. We emphasize that the existence of this inverse has
nothing to do with the possibility of solving the Wick relations for
the right hand side, thereby defining an ``anti-Wick ordering''.
This transformation would be connected with the existence $\Tt^{-1}$,
with $\Tt$ from \eq(H.4) rather than $S=T^{-1}$ with $T$ from
\eq(Bt.2). It is easy to find examples, where $\Tt$ is invertible,
while $T$ is not, and conversely.

Writing the expression in \eq(Db.b1b) as
$S_1(T_1T_2T_1-T_2T_1T_2)S_2=0$, we get the braid condition
$T_1T_2T_1=T_2T_1T_2$ for $T$. The verification that the braid
condition indeed suffices to define an algebra from relations
\eq(D.1), \eq(D.6), and \eq(Db.wz) can be done using the Diamond
Lemma \cit{Bergmann}, in the same way as in the proof of
\Thm\sDb.2/Db.2/ below.

If we read $\partial_i$ as $x_i\1$, it is natural to try to extend
the algebra still further, by including the adjoint differentials
$\dxs i$. We would thus arrive at a structure containing both
operations, involution as well as differential.
It turns out that the development in this section has put us on the
wrong track for doing this. Roughly, the reason is the following:
since we have a relation transforming $x\1\,\dx{}$ to $\dx{}\,x\1$,
with coefficients $S$, its adjoint should transform $x\,\dxs{}$ to
$\dxs{}\,x$ with coefficients $S^*$. Taking the differential of
either relation we get a rule transforming $\dxs{}\,\dx{}$ to
$\dx{}\,\dxs{}$, and the two possibilities match only when $S=S^*$.
But when we apply these rules to the differential $d(x_i^*x_j)$, and
compare with the basic Wick relation, we find that we should have
$S=T$, rather than $S=T^{-1}$. Nevertheless, the braid condition is
sufficient to define a structure allowing both involution and
differentials:

\iproclaim/Db.1/ Definition.
Let $I$ be finite, and let $T_{ij}^{k\ell}$, $i,j,k,\ell\in I$, be
the coefficients of a \HWa\ $\wick(T)$.
Then a {\bf Wick differential *-algebra} over $T$ is an associative
algebra $\Dffwick$ with the following properties:
\roster
\item[1]
$\Dffwick$ is generated by elements ``$f$'', and ``$\df f$'', for
$f\in\wick(T)$, where ``$\df f$'' depends linearly on $f$.
\item[2]
$\df\idty=0$, and $\df(fg)=(\df f)g+f\df g$.
\item[3]
$\Dffwick$ has an involution $\dagger$ extending that of $\wick(T)$
such that $\df(f\1)=(\df f)\1$.
\item[4]
There are complex coefficients $S_{ij}^{k\ell},R_{ij}^{k\ell}$ such
that
\eqgroup(Db.dw)$$\align
      x_i\dx j&= \textstyle \sum_{k\ell}\
                   R_{ij}^{k\ell}\   \dx\ell\, x_k
\tag{\deq\lasteq a(Db.R)}\\
    x_i\1\dx j&= \textstyle \sum_{k\ell}\
                   S_{ij}^{k\ell}\   \dx\ell\,x_k\1
\quad.\tag{\deq\lasteq b(Db.S)}\endalign$$
\item[5]
Monomials of the form
$$ \dx{i_1}\cdots\dx{i_m}\ x_{j_1}\cdots x_{j_n}\
   x\1_{k_p}\cdots x\1_{k_1}\ \dxs{\ell_q}\cdots \dxs{\ell_1}\
\quad,\tag\deq(Db.worder)$$
with $i_\alpha,j_\alpha,k_\alpha,\ell_\alpha\in I$, and
$m,n,p,q\geq0$, span $\Dffwick$.
\item[6]
The monomials (5) are linearly independent modulo the relations
$$ \dx i\dx j= - \textstyle \sum_{k\ell}\
                   R_{ij}^{k\ell}\   \dx\ell\, \dx k
\quad,\quad i,j\in I \quad,$$
and their adjoints.
\endroster
\eproclaim

\iproclaim/Db.2/ Theorem.
Let $T_{ij}^{k\ell}$ be as in \Def/Db.1/, $\HH=\Cx^I$, and define the
operators $T,R,S:\HH\otimes\HH\to\HH\otimes\HH$ by
$$
 T_{ij}^{k\ell}=\me<i\ell,T,jk>
\quad,\quad
 S_{ij}^{k\ell}=\me<i\ell,S,jk>
\quad,\quad\hbox{and}\quad
 R_{ij}^{k\ell}=\me<\ell k,R,ij>
\quad.$$
Then a Wick differential *-algebra over $T$ exists if and only if
$T$ is invertible, and satisfies the braid relation $T_1T_2T_1=T_2T_1T_2$.
In this case, the $S$ and $R$ must be chosen to be $S=T$, and
$R=T^{-1}$, and $\Dffwick\equiv\Dffwick(T)$ is uniquely determined
by $T$.
\eproclaim

\proof:
By taking adjoints and differentials of the relations \eq(Db.dw)
we get the following system of equations:
\eqgroup(Db.rel)$$\align
      x_i\1x_j&= \delta_{ij}\idty+ \textstyle \sum_{k\ell}\
                    \me<i\ell,T,jk>\, x_\ell\,x_k\1
\tag{\deq\lasteq a(Db.x+x)}\\
      x_i\dx j&= \hskip30pt\textstyle \sum_{k\ell}\
                     \me<\ell k,R,ij>\ \dx\ell\, x_k
\tag{\deq\lasteq b(Db.xdx)}\\
    x_i\1\dx j&= \hskip30pt\textstyle \sum_{k\ell}\
                     \me<i\ell,S,jk>\ \dx\ell\,x_k\1
\tag{\deq\lasteq c(Db.x+dx)}\\
   \dxs ix_j&= \hskip30pt\textstyle \sum_{k\ell}\
                      \me<i\ell,S^*,jk>\,x_\ell\,  \dxs k
\tag{\deq\lasteq d(Db.dx+x)}\\
 \dxs i x_j\1&= \hskip30pt\textstyle \sum_{k\ell}\
                       \me<ji,R^*,k\ell>\, x_\ell\1\,\dxs k
\tag{\deq\lasteq e(Db.dx+x+)}\\
 \dxs i \dx j&= \phantom{\delta_{ij}\idty }\textstyle -\sum_{k\ell}\
                       \me<i\ell,S^*,jk>\, \dx\ell\, \dxs k
\tag{\deq\lasteq f(Db.dx+dx)}\\
   \dx i\dx j&= \phantom{\delta_{ij}\idty }\textstyle -\sum_{k\ell}\
                       \me<\ell k,R,ij>\ \dx\ell\, \dx k\,
\tag{\deq\lasteq g(Db.dxdx)}\\
   \dxs i\dxs j&= \phantom{\delta_{ij}\idty }\textstyle -\sum_{k\ell}\
                       \me<ij,R^*,k\ell>\, \dxs\ell\, \dxs k
\tag{\deq\lasteq h(Db.dx+dx+)}\\
\endalign$$
We can obtain \eq(Db.dx+dx) as the differential of either
\eq(Db.x+dx) or \eq(Db.dx+x). These two have to agree, so that we
must have $S=S^*$. Taking the differential of \eq(Db.x+x), and
substituting \eq(Db.dx+x) and \eq(Db.x+dx), we find that $S=T$.
Thus the system \eq(Db.rel) becomes closed under both adjoints and
differentials.

The crucial part of the proof is to apply the Diamond Lemma
\cit{Bergmann} to this system, omitting the last two. The ``Wick
ordered form'' in which no further substitutions are possible is the
form given in \eq(Db.worder). It is obvious that the substitution
process terminates. We can apply any of the rules, whenever we find
two factors $x,\,x\1,\,\dx{},\,\dxs{}$ in the ``wrong'' order.  If
we find two occurrences of this in the same monomial, in
two non-overlapping pairs of factors, it is clear that either of the
two choices leads to the same result. The only situations we have to
check therefore, are those of three adjacent factors, which are in
the opposite order. Since we have only three types of factors, there
are only four types of such triples. Of these two are adjoints of
the other pair, which leaves us with checking the consistency of the
two possible reduction paths for the two types of products
$x_i\1\,x_j\,\dx k$, and $\dxs i\,x_j\,\dx k$. Now the first amounts
to two conditions, one arising from the first term in \eq(Db.x+x),
and one from the second. These are $RT=\idty$, and the braid
relation, respectively. The condition for the second type of triple
follows from these.

Hence the algebra defined by relations \eq\lasteq.a-f() has the
required properties. It remains to be checked that \eq\lasteq.g,h()
are compatible with this, i.e.\ that it makes no difference whether
we apply these relations before of after Wick ordering. For
$x\,\dx{}\,\dx{}\to \dx{}\,\dx{}\,x$ this follows from the braid
relation for $R$. Since the coefficients for the commutation of
coordinates and differentials are the same as for the commutation of
two differentials, the compatibility of \eq(Db.dxdx) with the
transformation $x\1\,\dx{}\,\dx{}\to \dx{}\,\dx{}\,x\1$ follows by
virtually the same computation as for $x\1\,x\,\dx{}\to
\dx{}\,x\,x\1$. For the same reason we need not check compatibility
for $\dxs{}\,\dx{}\,\dx{}$, and $\dx{}\,\dx{}\,\dx{}$, and the
compatibility for relation \eq(Db.dx+dx+) follows by taking
adjoints.
\QED

Since $\Dffwick$ is an involutive algebra we may once again ask
for {\em positive} \rep s. Rather than aiming at the development of
a general theory, we state some basic features of the \rep\ theory
of $\Dffwick(T)$. Note that, in spite of the difference between
\eq(D.6) and \eq(Db.R) with $R=T^{-1}$, \Lem/Dd.fcs/ applies without
change, because $(\idty+T)$ and $(\idty+T^{-1})$ have the same
kernel.

\iproclaim/Db.3/ Proposition.
Let $\Dffwick(T)$ be a differential \Wa\  over $T$.
\roster
\item
Then if $\Dff^p_0=\set0$ for some $p\geq2$, we have
$\pi(\dx i)=\pi(\dxs i)=0$ for all $i$.
\item
There is no universal bounded \rep\ of $\Dffwick(T)$.
\endroster
\eproclaim

\proof:
(1) Let $\pi$ be a positive \rep\ of $\Dffwick(T)$, and let $q$ be the
largest integer such that there is a $q$-form $\omega\in\Dff^p_0$
with $\pi(\omega)\neq0$. Clearly, $q< p<\infty$. On the other hand,
let $\omega\in\Dff^{q'}_0$ with $q/2<q'\leq q$. Then
$\pi\bigl((\omega\1\omega)^2\bigr)=0$, because after Wick ordering
this expression contains a form of degree $2q'>q$. Hence in
$\pi\bigl(1+\lambda\omega\1\omega\bigr)^2$ the term quadratic in
$\lambda$ vanishes, and, by positivity, the linear term has to
vanish, as well. Hence $\pi(\omega)=0$, contradicting the minimality
of $q$, unless $q=0$.

(2) Let $\pi$ be a bounded \rep. Then we can obtain a new one,
$\tilde\pi$ by setting $\tilde\pi(x_i)=\pi(x_i)$, and $\tilde\pi(\dx
i)=\lambda\pi(\dx i)$, with $\lambda\in\Cx$ arbitrary. Hence there
cannot be a universal bound on $\norm{\dx i}$, which holds in every
bounded \rep.
\QED

Note that the first item excludes the twisted canonical commutation
relations (\Exa/E.3/), by the discussion at the end of Section \sDd.
On the other hand, the differential \Wa\  of the twisted
canonical anti-commutation relations (\mCAR) seem to have
non-trivial \rep s. This is certainly the case for Clifford algebras
(\Exa/E.6/), as the following example shows: consider the Fock \rep\
of the (untwisted) ca\-no\-nical anti-commutation relations, whose
scalar product is clearly positive. Denote by $Z$ the operator
defined by $Z\Omega=\Omega$, and $Zx_i=-x_iZ$. Then we can define
$\dx i=\xi_i Z$, where $\xi_i\in\Cx$ is an arbitrary constant, and
it is immediately verified that the relations \eq(Db.dw) and
\eq(Db.rel) hold.

\head \secno\sG. Gauge automorphisms and their KMS states\endhead

{}From the universal property of $\wick(T)$ we see that any \Wa\ has a
$\1$-auto\-mor\-phism group $\gauge_t:\wick(T)\to\wick(T)$ defined by
$$ \gauge_t(k)=e^{it} k
\quad,\quad\hbox{for $t\in\Rl$, and $k\in I$,}
\tag\deq(G.1)$$
which we call the {\em gauge group} of $\wick(T)$.
That $\gauge_t$ is a $\1$-automorphism is equivalent to
$\gauge_t(k\1)=e^{-it}k\1$.
The pattern for defining these automorphisms
\Wa s sometimes applies also to other linear transformations of the
generators, depending on $T$.
The only thing one needs to check is
the invariance of the elements
$i\1\otimes j-\delta_{ij}\idty-\sum T_{ij}^{k\ell}\ell\otimes k\1$
in the tensor algebra $\tA(\HH,\HH\1)$. The definition of the
corresponding automorphism is then clear form the universal
property. What makes the gauge transformations special is that this
scheme works for {\em all} \Wa s.
The gauge automorphisms also extend to the universal bounded \rep\
$\wicc(T)$, if the latter exists, and we will then denote the
corresponding C*-automorphism group by $\gauge_t$, as well.

The gauge group defines on $\wick(T)$ a natural $\Ir$-grading by
$$ \degree(X)=n  \quad\iff\quad  \gauge_t(X)=e^{int}X
\quad.\tag\deq(G.2)$$
When we use the notation $\degree(X)$ it is implied that the element
$X$ is indeed of homogeneous degree. On a monomial in the generators
the degree is computed simply by counting the number of unstarred
generators $k$ and subtracting the number of starred ones $k\1$ in
the monomial. The degree zero, or gauge invariant part of $\wick(T)$
(resp.\ $\wicc(T)$) is denoted by $\wick(T)^\gauge$ (resp.\
$\wicc(T)^\gauge$). There is a faithful conditional expectation from
$\wicc(T)$ onto $\wicc(T)^\gauge$, defined by
$$ X\, \longmapsto
        {1\over2\pi}\int_{-\pi}^\pi dt\ \gauge_t(X)
\quad.\tag\deq(G.3)$$

\iproclaim/G.1/ Definition.
Let $\wick(T)$ be a \Wa, and $\lambda\in\Rl$. Then a {\bf gauge
KMS-functional} at fugacity $\lambda$ is a linear functional
$\kmsl$ on $\wick(T)$ such that $\kmsl(\idty)=1$,
$$  \kmsl\bigl(kX\bigr)
        =\lambda\,\kmsl\bigl(Xk\bigr)
\quad,\quad\hbox{and}\quad
    \kmsl\bigl(Xk\1\bigr)
        =\lambda\,\kmsl\bigl(k\1X\bigr)
\quad.$$
If this functional extends to $\wicc(T)$, it will be denoted by the
same letter.
\eproclaim

On the C*-algebra $\wicc(T)$, states with this property are indeed
precisely the KMS-states for the gauge group. To see this, note that
the polynomials in the generators form a dense subalgebra of
$\wicc(T)$ on which the gauge group is even analytic with
$\gauge_{i\beta}(X)=\exp\bigl(-\beta\degree(X)\bigr)X$. Now one of
the equivalent versions of the KMS-condition at inverse temperature
$\beta$ (e.g.\ 5.3.1 in \cit{BraRo}) requires the existence of such
a dense subalgebra of analytic elements, and the condition
$\omega\bigl(A\gauge_{i\beta}(B)\bigr)
      =\omega(BA)$, for $A,B$ in that subalgebra.
Hence the KMS-condition for is equivalent to
$$     \omega(XA)=e^{-\beta\degree(X)}\omega(AX)
\quad,\tag\deq(G.4)$$
for $A\in\wicc(T)$, and $X$ a homogeneous polynomial. By induction
on the degree of $X$ this is easily seen to be equivalent to the
definition of $\kmsl$, with $\lambda=e^{-\beta}$.
(The expression ``fugacity'' for the parameter $\lambda$ is borrowed
from statistical mechanics: in the language of physics the gauge
group is generated by the ``number operator'',  $-\beta$ is called
the ``chemical potential'', and ``fugacity'' for $\lambda$ is
standard terminology. The case $\lambda=0$ defines the Fock state,
which is the unique ground state \cit{BraRo} for $\gauge_t$,
formally corresponding to $\beta=\infty$, or temperature
$1/\beta=0$). Negative $\lambda$ are excluded by the following Lemma.

The most striking consequence of the \Wa\ structure for the
gauge-KMS states is that we obtain an explicit prescription for
computing $\kmsl$:

\iproclaim/G.2/ Lemma.
Let $\wick(T)$ be a \Wa, and $\lambda\in\Rl$, and let
$\kmsl$ be a gauge KMS-functional at fugacity $\lambda$. Then
\roster
\item[1]
For $\lambda\neq1$, $\degree(X)\neq0$ implies $\kmsl(X)=0$.
\item[2]
If $\kmsl(X\1X)\geq0$ for all $X\in\wick(T)$, we must have
$\lambda\geq0$.
\item[3]
Suppose that $\Norm\big{\Tt}\leq1$, and
$\abs\lambda<\Norm\big{\Tt}^{-1}$ where
$\Tt:\HH\1\otimes\HH\to\HH\otimes\HH\1$ is the operator introduced in
equation \eq(H.4). Then $\kmsl$ exists and is uniquely determined.
\item[4]
Suppose that $\Norm\big{\Tt}\leq1$. Then, for any $X\in\wick(T)$,
$\lambda\mapsto \kmsl(X)$ is an analytic function for
$\abs\lambda<\Norm\big{\Tt}^{-1}$.
\endroster
\eproclaim

\proof:
(1) $\kmsl(X)=\kmsl(X\idty)=\lambda^{\degree(X)}\kmsl(\idty X)$,
hence $\lambda^{\degree(X)}=1$, whenever $\kmsl(X)\neq0$.

(2) When $\kmsl$ is positive,
$\kmsl(kk\1)=\lambda\kmsl(k\1k)$ implies that either $\lambda\geq0$
or $\kmsl(k\1k)=0$ for all $k\in I$. The latter condition would
imply $\kmsl(XkY)=\lambda^{-\degree(Y)}\kmsl(YXk)=0$, so all $k$
would be annihilated by the GNS-\rep\ with respect to $\kmsl$, which
contradicts the Wick relations \eq(H.1).

(3) The prescription for evaluating $\kmsl$ on all Wick ordered
monomials of the form $X=i_1\cdots i_nj\1_m\cdots j\1_1$ is the
following: we exchange the two groups of generators, to get
$\kmsl(X)=\lambda^n\kmsl(\widetilde X)$, where $\widetilde
X=j_m\1\cdots j_1\1i_1\cdots i_n$. Note that $\lambda^n=\lambda^m$,
or $\kmsl(X)=0$ by (1). Now we Wick order $\widetilde X$. This
produces a linear combination of Wick ordered monomials with degrees
$n'\leq n, m'\leq m$. We collect all terms of the same degree as
$X$. This leading term is computed by using the Wick relations
without constant term, and is described by an operator
$$ \Tt^{(n,m)}:(\HH\1)^{\otimes m}\otimes\HH^{\otimes n}
      \longrightarrow \HH^{\otimes n}\otimes(\HH\1)^{\otimes m}
\quad,$$
inductively defined from $\Tt^{(1,1)}=\Tt$ from equation \eq(H.4).
$\Tt^{(n,m)}$ is the product of $nm$ copies of $\Tt$, acting in
appropriate tensor factors. Hence
$\Norm\big{\Tt^{(n,m)}}\leq\Norm\big{\Tt}^{nm}$. Now
$$ \kmsl(X-\lambda^n\Tt^{n,m}\widetilde X)
    =\kmsl(\hbox{terms of lower degree)}
$$
is an inductive formula for $\kmsl$ provided the operator
$X\mapsto \bigl(X-\lambda^n\Tt^{n,m}\widetilde X\bigr)$ is
invertible on $\HH^{\otimes n}\otimes(\HH\1)^{\otimes m}$. Clearly, it
is sufficient for this that
$\abs\lambda^n \Norm\big{\Tt^{(n,m)}}
      <(\abs\lambda\Norm\big{\Tt})^n\Norm\big{\Tt}^m<1$,
which is in turn implied by the assumption made in (3).

(4) is immediate from (3) and the analyticity of the inverse
operator used in (3).
\QED

For special choices of $T$ the uniqueness result can be improved.
What we have used in the proof is only the invertibility of the
operator $X\mapsto \bigl(X-\lambda^n\Tt^{n,m}\widetilde X\bigr)$.
For the $q$-relations (\Exa/E.1/) this is simply
$(1-\lambda^nq^{nm})\idty$. Hence in this case uniqueness holds
whenever $\lambda\neq q^{-m}$ for all $m$. This includes the value
$\lambda=1$, for which $\kmsl$ is a trace, and almost all
$\lambda>1$ (for which the ``temperature'' $\beta=-\ln\lambda$ is
negative). On the other hand, we know that the algebra $\wicc(T)$ in
this case has no tracial states. Hence the functional $\kms1$ cannot
be positive. The following Theorem collects the basic results on
the positivity of $\kmsl$.

\iproclaim/G.3/ Theorem.
Let $\wick(T)$ be a \Wa\ with $d=\abs I<\infty$ generators,
and suppose that the Fock state $\cst0=\kms0$ is positive. Let $P_n$
be the operators on $\HH^{\otimes n}$ defined in \eq(Bt.4).
Then the following conditions are equivalent:
\roster
\item[1]
$\kmsl$ is positive, and normal with respect to the Fock \rep.
\item[2]
$\sum_{n=0}^\infty \lambda^n\rank P_n<\infty$.
\endroster
In particular, $\kmsl$ is positive for $0\leq\lambda\leq1/d$, and if
$P_n>0$ for all $n$, $\kms{1/d}$ is not normal with respect to the
Fock \rep.
\eproclaim

\proof:
In the Fock \rep\ the gauge group is implemented by the number
operator $N$, defined by $(N-n)i_1\cdots i_n \idty=0$. Then if the
operator $\lambda^N$ is trace class, it defines the density matrix
of a KMS state $\kmsl$ for $\gauge_t$. The converse holds, because
the Fock state is pure, and hence the Fock \rep\ is irreducible.
Hence (1) is equivalent to
$$ \tr\lambda^N
     =\sum_{n=0}^\infty \lambda^n\rank P_n<\infty
\quad.$$
Since $\rank P_n\leq d^n$, this criterion holds when
$\lambda d<1$. Moreover, when $P_n>0$, $\rank P_n= d^n$, and the
criterion fails for $\lambda d=1$.
\QED

With the help of this Theorem we can now discuss the existence
of gauge KMS-states for the examples in Section \sE.
For the generalized $q$-relations (\Exa/E.5/, \eq(E.7)) we have
$\Norm\big{\Tt}=\max_{ij}\abs{q_{ij}}$. Hence for constant
$q_{ij}\equiv q$, and $\abs q<1$ we have $\Norm\big{\Tt}<1$, and
we have unique KMS states for $0\leq\lambda\leq1/d$. For $q=0$ it is
easy to see (Example 5.3.27 in \cit{BraRo}) that there cannot be a
state $\kmsl$ with $\lambda>1/d$. More generally, we know for
any sufficiently small $T$, i.e.\ for $T$ satisfying the hypothesis
of \Thm/Bs.6/, the exact range of the parameter $\lambda$
for which $\kmsl\geq0$. This follows readily from the isomorphism
$\wicc(T)\cong\wicc(0)$ established in that Theorem, and the
observation that the isomorphism intertwines the respective gauge
groups. The growth $\rank P_n=d^n$ is typical for all these
examples. It is also responsible for the anomalous thermodynamical
behaviour of the ``free quon gas'' \cit{QTD}.

For the twisted canonical (anti-) commutation relations the
dimension of the $n$-particle Fock space grows exactly as for the
untwisted counterparts. This agreement of the ``Poincar\'e
series'' \cit{WessZ} with the undeformed case, is, in fact, one of
the motivations for studying just this type of deformation.
Thus we have
$$  \tr(\lambda^N)=\cases  (1-\lambda)^{-d} &\text{ for \eq(E.3)}\\
                           (1+\lambda)^d    &\text{ for \eq(E.4).}\\
\endcases$$
Hence the range of $\lambda$ is $\bracks{0,1}$ in the Bosonic case,
and $\bracks{0,\infty}$ in the Fermionic case ($\lambda=\infty$
defines a ``ceiling state'' \cit{BraRo}).

The state $\kms1$ in the Bosonic case, which is positive by
continuity of $\kmsl$ in $\lambda$, is a trace. Indeed, the algebra
$\wicc(T)$ in this case has a tracial state. In the corresponding
GNS-\rep\ all generators are zero, except the first, which is
unitary up to a factor.

\head Acknowledgement \endhead

We are indebted to Marek \Bozejko\ and Roland Speicher for making
their unpublished work \cit{BoSpein} available to us.
PJ was supported in part by the U.S. NSF, NATO, a University of Iowa
Faculty Scholar Award, and by the U.I.--Oakdale Institute for
Advanced Studies.
PJ worked at the end of the project at the University of Oslo, with
support from the Norwegian Research council.
RW was supported by a fellowship and a travel grant from the
Deutsche Forschungsgemeinschaft (Bonn).

\ppt{\parskip=0pt}
\enddocument
\bye